\documentclass[a4paper,12pt]{article}

\usepackage{amsmath,amssymb,amsthm,amscd}

\newcommand{\Fg}{\mathfrak{g}}
\newcommand{\Fh}{\mathfrak{h}}

\newcommand{\BB}{\mathbb{B}}
\newcommand{\BP}{\mathbb{P}}
\newcommand{\BR}{\mathbb{R}}
\newcommand{\BZ}{\mathbb{Z}}
\newcommand{\BQ}{\mathbb{Q}}

\newcommand{\CB}{\mathcal{B}}

\newcommand{\wt}{\mathop{\rm wt}\nolimits}
\newcommand{\cl}{\mathop{\rm cl}\nolimits}
\newcommand{\rr}{\Delta^{\mathrm{re}}}
\newcommand{\prr}{\Delta^{\mathrm{re}}_{+}}
\newcommand{\id}{\mathop{\rm id}\nolimits}
\newcommand{\ext}{\mathop{\rm ext}\nolimits}
\newcommand{\res}{\mathop{\rm res}\nolimits}

\newcommand{\Hom}{\mathop{\rm Hom}\nolimits}
\newcommand{\Bij}{\mathop{\rm Bij}\nolimits}
\newcommand{\GL}{\mathop{\rm GL}\nolimits}
\newcommand{\Deg}{\mathop{\rm Deg}\nolimits}

\newcommand{\vpi}{\varpi}
\newcommand{\ve}{\varepsilon}
\newcommand{\vp}{\varphi}

\newcommand{\bi}{{\bf i}}
\newcommand{\bzero}{{\bf 0}}

\newcommand{\vsp}{\vspace{3mm}}

\newcommand{\fin}[1]{\overset{\circ}{#1}}
\newcommand{\ti}[1]{\widetilde{#1}}
\newcommand{\ud}[1]{\underline{#1}}

\makeatletter
\@addtoreset{equation}{subsection}

\renewcommand\section{\@startsection{section}{1}{0pt}
{-3.5ex plus -1ex minus -.2ex}{1.0ex plus .2ex}{\large\bf}}
\renewcommand\subsection{\@startsection{subsection}{1}{0pt}
{2.5ex plus 1ex minus .2ex}{-1em}{\bf}}
\makeatother

\theoremstyle{plain}
\newtheorem{thm}{Theorem}[subsection]
\newtheorem{lem}[thm]{Lemma}
\newtheorem{prop}[thm]{Proposition}
\newtheorem{cor}[thm]{Corollary}

\newtheorem{ithm}{Theorem}
\newtheorem{icor}[ithm]{Corollary}

\newtheorem*{claim}{Claim}
\newtheorem*{sublem}{Sublemma}

\theoremstyle{definition}
\newtheorem{dfn}[thm]{Definition}

\theoremstyle{remark}
\newtheorem{rem}[thm]{Remark}

\begin{document}
\setlength{\baselineskip}{18pt}

\title{{\Large\bf 
Lakshmibai-Seshadri paths of level-zero weight shape \\
and one-dimensional sums associated to \\
level-zero fundamental representations}}
\author{
 Satoshi Naito \\ 
 \small Institute of Mathematics, University of Tsukuba, \\
 \small Tsukuba, Ibaraki 305-8571, Japan \ 
 (e-mail: {\tt naito@math.tsukuba.ac.jp})
 \\[2mm] and \\[2mm]
 Daisuke Sagaki \\ 
 \small Institute of Mathematics, University of Tsukuba, \\
 \small Tsukuba, Ibaraki 305-8571, Japan \ 
 (e-mail: {\tt sagaki@math.tsukuba.ac.jp})
}
\date{}
\maketitle

%=======================%
%     START ABSTRACT    %
%=======================%
%
%
\begin{abstract}
We give interpretations of energy functions and 
(classically restricted) one-dimensional sums 
associated to tensor products of level-zero fundamental representations 
of quantum affine algebras in terms of 
Lakshmibai-Seshadri paths of level-zero weight shape.
\end{abstract}
%
%
%
%=========================%
%     START SECTION 01    %
%=========================%
%
\section{Introduction.}
\label{sec:intro}
A one-dimensional (configuration) sum (1dsum for short) is 
a weighted sum over certain one-dimensional configurations, 
where the weights are given by an ``energy function'', 
and arose from the studies of solvable lattice models 
in statistical mechanics by Baxter's corner 
transfer matrix method. 
However, the crystal basis theory of Kashiwara provided 
an intrinsic definition of a 1dsum, and a conceptual proof of 
the fact that in the infinite lattice size limit, a 1dsum tends to 
the character of a highest weight module over an affine Lie algebra
(see \cite{KMN}).
The purpose of this paper is to give interpretations of energy functions and 
hence of classically restricted 1dsums associated 
to tensor products of certain finite-dimensional irreducible modules,
called level-zero fundamental representations,
over a quantum affine algebra
via Lakshmibai-Seshadri paths (LS paths for short) 
of level-zero weight shape.

Let us explain our results more precisely. 
Let $\Fg$ be an affine Lie algebra over $\BQ$
with Cartan subalgebra $\Fh$, 
integral weight lattice $P \subset \Fh^{\ast}$, 
simple roots $\bigl\{\alpha_{j}\bigr\}_{j \in I} \subset \Fh^{\ast}$, 
null root $\delta=\sum_{j \in I} a_{j}\alpha_{j} \in \Fh^{\ast}$, 
Weyl group $W$, and level-zero fundamental weights
$\bigl\{\vpi_{i}\bigr\}_{i \in I_{0}} \subset P$,
where $I$ is an index set with a distinguished index ``0'',
and $I_{0}:=I \setminus \{0\}$.
Also, let $\Fg_{I_{0}}$ be
the canonical finite-dimensional Lie subalgebra
of $\Fg$ corresponding to the subset $I_{0} \subset I$. 
We denote by $U_{q}'(\Fg)$ the quantum affine algebra 
with weight lattice $P_{\cl}:=\cl(P) \subset \Fh^{\ast}/\BQ\delta$, 
where $\cl:\Fh^{\ast} \twoheadrightarrow \Fh^{\ast}/\BQ\delta$ 
is the canonical projection, and by $U_{q}'(\Fg)_{I_{0}}$ 
the subalgebra of $U_{q}'(\Fg)$ corresponding to 
the subset $I_{0} \subset I$.
In \cite{HKOTY}, \cite{HKOTT}, 
to equate with fermionic formulas $M$ originated from 
the Bethe Ansatz in solvable lattice models, 
they introduced a special kind of 
classically restricted 1dsums $X$ as a specific $q$-analogue of the 
multiplicities of an irreducible $U_{q}'(\Fg)_{I_{0}}$-module in 
tensor products of the so-called Kirillov-Reshetikhin modules
(KR modules for short) $W_{s}^{(i)}$, $i \in I_{0}$, 
$s \in \BZ_{\ge 1}$, over $U_{q}'(\Fg)$, 
under the assumption that the KR modules $W_{s}^{(i)}$, $i \in I_{0}$,
$s \in \BZ_{\geq 1}$, have
``simple'' crystal bases $\CB^{i,s}$, called KR crystals.
Here the KR modules $W_{s}^{(i)}$, $i \in I_{0}$, 
$s \in \BZ_{\ge 1}$, are finite-dimensional irreducible 
$U_{q}'(\Fg)$-modules having the conjectural irreducible 
decomposition as a $U_{q}'(\Fg)_{I_{0}}$-module.

Soon afterward, in \cite{Kas-lz}, 
Kashiwara constructed finite-dimensional irreducible 
$U_{q}'(\Fg)$-modules $W(\vpi_{i})$, $i \in I_{0}$,
called level-zero fundamental representations, and proved that they have
simple crystal bases $\CB(\vpi_{i})$.
In a series of papers \cite{NSp1}--\cite{NSls}, we studied the crystal $\BB(\lambda)$
of all LS paths of shape $\lambda$ for a level-zero integral weight
$\lambda \in P^{0}_{+}:=\sum_{i \in I_{0}} \BZ_{\ge 0} \vpi_{i}$ 
of the form $\lambda = \sum_{i \in I_{0}} \lambda^{(i)} \vpi_{i}$ 
with $\lambda_{(i)} \in \BZ_{\ge 0}$, and also 
the associated $P_{\cl}$-weighted crystal 
($P_{\cl}$-crystal for short) $\BB(\lambda)_{\cl}$.
The main results of \cite{NSp1} and \cite{NSp2} 
show that the $P_{\cl}$-crystal $\BB(\vpi_{i})_{\cl}$ is 
isomorphic to the crystal basis $\CB(\vpi_{i})$ of 
$W(\vpi_{i})$ for each $i \in I_{0}$.

We should mention that through enough evidence 
(see, for example, \cite{Kas-rims} and \cite{FL}), 
it is confirmed that the level-zero fundamental representation $W(\vpi_{i})$ 
is indeed the KR module $W_{1}^{(i)}$, and hence the crystal basis 
$\CB(\vpi_{i})$ is indeed the KR crystal $\CB^{i,1}$ for every $i \in I_{0}$.
In this paper, following the definition in \cite{HKOTY} and \cite{HKOTT} 
of classically restricted 1dsums $X$ associated 
to tensor products of KR crystals 
$\CB^{i,s}$, $i \in I_{0}$, $s \in \BZ_{\ge 1}$, 
we define classically restricted 1dsums associated 
to tensor products of the simple crystals 
$\BB(\vpi_{i})_{\cl}$ ($\simeq \CB(\vpi_{i})$), $i \in I_{0}$, 
in place of $\CB^{i,1}$, 
$i \in I_{0}$, as follows.
Let $\bi = (i_{1}, i_{2}, \ldots, i_{n})$ be 
a sequence of elements of $I_{0}$ (with repetitions allowed), 
and set 
$\BB_{\bi} := 
 \BB(\vpi_{i_{1}})_{\cl} \otimes 
 \BB(\vpi_{i_{2}})_{\cl} \otimes \cdots \otimes
 \BB(\vpi_{i_{n}})_{\cl}$.
Then, for an element 
$\mu \in \cl(P^{0}_{+})=
 \sum_{i \in I_{0}}\BZ_{\ge 0} \cl(\vpi_{i})$, 
the classically restricted 1dsum 
$X(\BB_{\bi},\mu; q)$ is defined by:
\begin{equation*}
X(\BB_{\bi},\mu; q)=\sum_{
 \begin{subarray}{c}
 b \in \BB_{\bi} \\[1mm]
 e_{j}b=\bzero \ (j \in I_{0}) \\[1mm]
 \wt b = \mu
 \end{subarray}} 
 q^{D_{\bi}(b)},
\end{equation*}
where $D_{\bi}:\BB_{\bi} \rightarrow \BZ$ is the energy function 
associated to the tensor product $P_{\cl}$-crystal 
$\BB_{\bi}$ (see \S\ref{subsec:main} for its definition), and 
$e_{j}$, $j \in I_{0}$, are 
the (raising) root operators on $\BB_{\bi}$ for 
the canonical Lie subalgebra $\Fg_{I_{0}}$ of $\Fg$. 

Now, for a level-zero integral weight 
$\lambda \in P$ of the form 
$\lambda =\sum_{i \in I_{0}} \lambda^{(i)} \vpi_{i}$ 
with $\lambda^{(i)} \in \BZ_{\ge 0}$, we define a function 
$\Deg_{\lambda}:\BB(\lambda)_{\cl} \rightarrow \BZ_{\le 0}$ 
on the crystal $\BB(\lambda)_{\cl}$ as follows. 
The $P_{\cl}$-crystal $\BB(\lambda)_{\cl}$ is,
by definition, the set of all paths 
$[0,1] \rightarrow \Fh_{\BR}^{\ast}/\BR \delta$ of the form 
$\cl(\pi)=\cl \circ \pi$, where $\pi:[0,1] \rightarrow 
\Fh_{\BR}^{\ast}:= \BR \otimes_{\BQ} \Fh^{\ast}$ 
is an LS path of shape $\lambda$ and
$\cl:\Fh_{\BR}^{\ast} \rightarrow \Fh_{\BR}^{\ast}/\BR\delta$ is 
the canonical projection. 
Let $\eta \in \BB(\lambda)_{\cl}$.
Then there exists a unique LS path $\pi_{\eta} \in \BB_{0}(\lambda)$ 
such that $\cl(\pi_{\eta})=\eta$ and such that the ``initial direction'' 
$\nu_{1} \in W\lambda$ of $\pi_{\eta}$ lies 
in the set $\lambda-\fin{Q}_{+}$, 
where $\BB_{0}(\lambda)$ denotes 
the connected component of $\BB(\lambda)$ 
containing the straight line 
$\pi_{\lambda}(t)=t\lambda$, $t \in [0,1]$, and 
$\fin{Q}_{+}:=\sum_{j \in I_{0}} \BZ_{\ge 0} \alpha_{j}$. 
Note that $\pi_{\eta}(1) \in P$ is of the form 
$\lambda-a_{0}^{-1}\beta+a_{0}^{-1}k\delta$, 
with $\beta \in \fin{Q}_{+}$ and $k \in \BZ_{\ge 0}$. 
Thus the degree $\Deg_{\lambda}(\eta) \in \BZ_{\le 0}$ of 
$\eta \in \BB(\lambda)_{\cl}$ is defined to be
the nonpositive integer $-k \in \BZ_{\le 0}$. 

We recall from \cite{NSz} that for every sequence 
$\bi = (i_{1},\,i_{2},\,\dots,\,i_{n})$ 
of elements of $I_{0}$, the tensor product $P_{\cl}$-crystal 
$\BB_{\bi}$ is isomorphic to the $P_{\cl}$-crystal 
$\BB(\lambda)_{\cl}$, 
where $\lambda := \sum_{k=1}^{n} \vpi_{i_{k}} \in P^{0}_{+}$.
Thus, there exists a (unique) isomorphism of $P_{\cl}$-crystals
\begin{equation*}
\Psi_{\bi}:
 \BB(\lambda)_{\cl} \stackrel{\sim}{\rightarrow} 
 \BB_{\bi} = 
 \BB(\vpi_{i_{1}})_{\cl} \otimes 
 \BB(\vpi_{i_{2}})_{\cl} \otimes \cdots \otimes 
 \BB(\vpi_{i_{n}})_{\cl}.
\end{equation*}
Our main result of this paper is the following
description of the energy function 
$D_{\bi}:\BB_{\bi} \rightarrow \BZ$ 
associated to the $P_{\cl}$-crystal 
$\BB_{\bi}$ in terms of the function 
$\Deg_{\lambda}:\BB(\lambda)_{\cl} \rightarrow \BZ_{\leq 0}$. 
%
%%%%%%%%%%%%%%%%%
%%% ithm:main %%%
%%%%%%%%%%%%%%%%%
%
\begin{ithm} \label{ithm:main}
Let $\bi = (i_1, i_2, \ldots, i_n)$ 
be an arbitrary sequence of elements of $I_{0}$,
and set 
$\BB_{\bi}:= 
 \BB(\vpi_{i_{1}})_{\cl} \otimes 
 \BB(\vpi_{i_{2}})_{\cl} \otimes \cdots \otimes 
 \BB(\vpi_{i_{n}})_{\cl}$, 
$\lambda := \sum_{k=1}^{n} \vpi_{i_{k}} \in P^{0}_{+}$.
Then, for every $\eta \in \BB(\lambda)_{\cl}$, 
the following equation holds\,{\rm:}
\begin{equation*}
\Deg_{\lambda}(\eta) = D_{\bi}(\Psi_{\bi}(\eta)) - D_{\bi}^{\ext},
\end{equation*}
where the $D_{\bi}^{\ext} \in \BZ$ is 
a specific constant 
{\rm(}see \S\ref{subsec:main} for its explicit definition{\rm)}.  
\end{ithm}
Because the root operators $e_{j}$, $j \in I_{0}$, on $\BB(\lambda)_{\cl}$ 
and those on $\BB_{\bi}$ are compatible with the isomorphism 
$\Psi_{\bi}:\BB(\lambda)_{\cl} \rightarrow \BB_{\bi}$ of 
$P_{\cl}$-crystals, we obtain the following corollary.
%
%%%%%%%%%%%%%%%%%
%%% icor:main %%%
%%%%%%%%%%%%%%%%%
%
\begin{icor} \label{icor:main}
Keep the notation of Theorem~\ref{ithm:main}. 
For every $\mu \in \cl(P^{0}_{+})=
\sum_{i \in I_{0}} \BZ_{\ge 0} \cl(\vpi_{i})$, 
the following equation holds\,{\rm:}
\begin{equation*}
\sum_{
 \begin{subarray}{c}
 \eta \in \BB(\lambda)_{\cl} \\[1mm]
 e_{j}\eta=\bzero \ (j \in I_{0}) \\[1mm]
 \eta(1) = \mu
 \end{subarray}} 
 q^{\Deg_{\lambda}(\eta)}
= q^{-D_{\bi}^{\ext}} X(\BB_{\bi}, \mu; q), 
\end{equation*}
where $e_{j}$, $j \in I_{0}$, are the {\rm(}raising\,{\rm)}
root operators on $\BB(\lambda)_{\cl}$. 
\end{icor}

Now we restrict our attention to the case 
in which $\Fg$ is of type $A_{\ell-1}^{(1)}$. 
Let $\bi = (i_{1},\,i_{2},\,\dots,\,i_{n})$ be 
an arbitrary sequence of elements of 
$I_{0}=\bigl\{1,\,2,\,\dots,\,\ell-1\bigr\}$ such that 
$i_{1} \ge i_{2} \ge \cdots \ge i_{n}$, and set 
$\lambda:=\sum_{k=1}^{n} \vpi_{i_{k}} \in P^{0}_{+}$; 
we denote this sequence $\bi=(i_{1},\,i_{2},\,\dots,\,i_{n})$ 
by $\lambda^{\dagger}$ when we regard it 
as a partition (or a Young diagram). 
Also, in the following, we identify an element 
$\mu=\sum_{i \in I_{0}} \mu^{(i)}\cl(\vpi_{i}) \in 
\cl(P^{0}_{+})$ with the partition
\begin{equation*}
\left(
 \sum_{i=1}^{\ell-1}\mu^{(i)}, \ 
 \sum_{i=2}^{\ell-1}\mu^{(i)}, \ \dots, \ 
 \mu^{(\ell-1)}
\right).
\end{equation*}
Then, from \cite[Corollary~4.3]{NY},
we deduce (see \S\ref{subsec:kostka} for details) that 
for every $\mu \in \cl(P^{0}_{+})=\sum_{i \in I_{0}} \BZ_{\ge 0}\cl(\vpi_{i})$, 
the Kostka-Foulkes polynomial $K_{\mu^{t},\,\lambda^{\dagger}}(q)$
associated to the conjugate (or transpose) $\mu^{t}$ 
of the partition $\mu$ and the partition $\lambda^{\dagger}$ is 
equal to the following:
\begin{equation*}
X(\BB_{\bi},\mu; q^{-1}) =\sum_{
 \begin{subarray}{c}
 b \in \BB_{\bi} \\[1mm]
 e_{j}b=\bzero \ (j \in I_{0}) \\[1mm]
 \wt b = \mu
 \end{subarray}} 
 q^{-D_{\bi}(b)}.
\end{equation*}
Thus we obtain the following expression for 
the Kostka-Foulkes polynomials in terms of LS paths 
(note that the constant $D_{\bi}^{\ext} \in \BZ$ 
vanishes in the case of type $A_{\ell-1}^{(1)}$).
%
%%%%%%%%%%%%%%%%%%%
%%% icor:kostka %%%
%%%%%%%%%%%%%%%%%%%
%
\begin{icor} \label{icor:kostka}
Assume that $\Fg$ is of type $A_{\ell-1}^{(1)}$.
Let $\mu \in \cl(P^{0}_{+})=\sum_{i \in I_{0}} \BZ_{\ge 0}\cl(\vpi_{i})$. 
Then, with the notation above, we have
\begin{equation*}
K_{\mu^{t},\,\lambda^{\dagger}}(q)= 
\sum_{
 \begin{subarray}{c}
 \eta \in \BB(\lambda)_{\cl} \\[1mm]
 e_{j}\eta=\bzero \, (j \in I_{0}) \\[1mm]
 \eta(1)=\mu
 \end{subarray}} 
 q^{-\Deg_{\lambda}(\eta)}.
\end{equation*}
\end{icor}

This paper is organized as follows. 
In \S\ref{sec:review}, 
we first fix our notation for quantum affine algebras. 
Then we briefly review some fundamental facts about 
LS path crystals with weight lattice $P$ or $P_{\cl}$, 
and those about simple $P_{\cl}$-crystals for quantum affine algebras.
In \S\ref{sec:degree},
we define our ``degree functions'' on 
$P_{\cl}$-crystals of LS paths of level-zero weight shape, 
and show some of their basic properties. 
In \S\ref{sec:main},
we first state our main result (Theorem~\ref{ithm:main}) 
describing energy functions associated to 
level-zero fundamental representations. 
Then, we give a proof of it after showing 
a key proposition to our proof. 
Finally, we mention the relation to 
classically restricted 1dsums and the Kostka-Foulkes polynomials 
(Corollaries~\ref{icor:main} and \ref{icor:kostka}).

%=========================%
%     START SECTION 02    %
%=========================%
%
\section{Preliminaries.}
\label{sec:review}

%==============================%
%     START SUBSECTION 0201    %
%==============================%
%
\subsection{Affine Lie algebras and quantum affine algebras.}
\label{subsec:affine}

Let $A=(a_{ij})_{i,j \in I}$ 
be a generalized Cartan matrix of affine type. 
Throughout this paper, we assume that 
the elements of the index set $I$ are numbered 
as in \cite[\S4.8, Tables Aff~1--Aff~3]{Kac}. 
Take a special vertex $0 \in I$ as in these tables, 
and set $I_{0}:=I \setminus \{0\}$. 
Let $\Fg=\Fg(A)$ be the affine Lie algebra 
associated to the Cartan matrix 
$A=(a_{ij})_{i,j \in I}$ of affine type 
over the field $\BQ$ of rational numbers, and 
let $\Fh$ be its Cartan subalgebra. 
Note that 
$\Fh=\bigl(\bigoplus_{j \in I} \BQ h_{j}\bigr) \oplus \BQ d$, 
where 
$\Pi^{\vee}:=\bigl\{h_{j}\bigr\}_{j \in I} \subset \Fh$ is 
the set of simple coroots, and 
$d \in \Fh$ is the scaling element. 
Also, we denote by 
$\Pi:=\bigl\{\alpha_{j}\bigr\}_{j \in I} \subset 
\Fh^{\ast}:=\Hom_{\BQ}(\Fh,\BQ)$ 
the set of simple roots, and by 
$\Lambda_{j} \in \Fh^{\ast}$, $j \in I$, 
the fundamental weights; 
note that $\alpha_{j}(d)=\delta_{j,0}$ and 
$\Lambda_{j}(d)=0$ for $j \in I$. 
Let 
$\delta=\sum_{j \in I} a_{j}\alpha_{j} \in \Fh^{\ast}$ and 
$c=\sum_{j \in I} a^{\vee}_{j} h_{j} \in \Fh$ be 
the null root and 
the canonical central element of $\Fg$, 
respectively. 
Here we should note that 
$a_{0}=2$ if $\Fg$ is of type $A^{(2)}_{2\ell}$, and 
$a_{0}=1$ otherwise. 
We define the Weyl group $W$ of $\Fg$ by: 
$W=\langle r_{j} \mid j \in I\rangle 
\subset \GL(\Fh^{\ast})$, where 
$r_{j} \in \GL(\Fh^{\ast})$ is the simple reflection 
associated to $\alpha_{j}$ for $j \in I$, and then 
define the set $\rr$ of real roots by: $\rr=W\Pi$. 
The set of positive real roots is denoted by 
$\prr \subset \rr$. 
Also, let us denote 
by $(\cdot\,,\,\cdot)$ the (standard) 
bilinear form on $\Fh^{\ast}$ 
normalized so that $(\alpha_{i},\alpha_{j}) \in \BZ$ 
for all $i,\,j \in I$.

We take a dual weight lattice $P^{\vee}$ 
and a weight lattice $P$ as follows:
%
%%%%%%%%%%%%%%%%%%%
%%% eq:lattices %%%
%%%%%%%%%%%%%%%%%%%
%
\begin{equation} \label{eq:lattices}
P^{\vee}=
\left(\bigoplus_{j \in I} \BZ h_{j}\right) \oplus \BZ d \, 
\subset \Fh
\quad \text{and} \quad 
P= 
\left(\bigoplus_{j \in I} \BZ \Lambda_{j}\right) \oplus 
   \BZ a_{0}^{-1}\delta
   \subset \Fh^{\ast}.
\end{equation}
It is clear that $P$ contains 
all the simple roots $\alpha_{j}$, $j \in I$, and that 
$P \cong \Hom_{\BZ}(P^{\vee},\BZ)$. 
The quintuplet
$(A,P,P^{\vee},\Pi,\Pi^{\vee})$ is called 
a Cartan datum for the generalized Cartan matrix 
$A=(a_{ij})_{i,j \in I}$ of affine type 
(see \cite[Definition~2.1.1]{HK}).

Let $\Fg_{I_{0}}$ be the canonical finite-dimensional 
Lie subalgebra of $\Fg$ generated 
by $x_{j}$, $y_{j}$, $j \in I_{0}$, and 
$h_{j}$, $j \in I$, 
where $x_{j}$ (resp., $y_{j}$) is 
a nonzero root vector of $\Fg$ corresponding to 
the simple root $\alpha_{j}$ (resp., $-\alpha_{j}$); 
note that $\Fh_{I_{0}}:=\bigoplus_{j \in I} \BQ h_{j}$ is 
the Cartan subalgebra of $\Fg_{I_{0}}$.
We denote by $\fin{W}$ the subgroup of $W$ generated by 
$r_{j}$, $j \in I_{0}$, which can be thought of as 
the Weyl group of the Lie subalgebra 
$\Fg_{I_{0}} \subset \Fg$. 
Let $w_{0} \in \fin{W}$ denote 
the longest element of $\fin{W}$. 

%%%%%%%%%%%%%%%
%%% dfn:lv0 %%%
%%%%%%%%%%%%%%%
%
\begin{dfn} \label{dfn:lv0}
An integral weight 
$\lambda \in P$ is said to 
be level zero if $\lambda(c)=0$. 
In addition, 
a level-zero integral weight 
$\lambda \in P$ is said to be level-zero dominant 
(resp., strictly level-zero dominant) 
if $\lambda(h_{j}) \ge 0$ 
(resp., $\lambda(h_{j}) > 0$) 
for all $j \in I_{0}$.
\end{dfn}

For each $i \in I_{0}=I \setminus \{0\}$, 
we define a level-zero fundamental weight 
$\vpi_{i} \in P$ by:
\begin{equation}
\vpi_{i}=\Lambda_{i}-a_{i}^{\vee}\Lambda_{0}.
\end{equation}
Note that the $\vpi_{i}$ is actually
a level-zero dominant integral weight; 
in fact, $\vpi(c)=0$ and 
$\vpi_{i}(h_{j})=\delta_{i,j}$ 
for $j \in I_{0}$. 
We set 
\begin{equation}
P_{+}^{0}:=\sum_{i \in I_{0}} \BZ_{\ge 0}\vpi_{i}.
\end{equation}
The next lemma follows immediately from 
\cite[Lemma~2.3.2 and Remark~4.1.1]{NSls} and 
the proof of \cite[Lemma~2.3.3]{NSls} 
by noting the linear independence of 
$\alpha_{j}$, $j \in I_{0}$, and $\delta$.
%
%%%%%%%%%%%%%%%
%%% lem:orb %%%
%%%%%%%%%%%%%%%
%
\begin{lem} \label{lem:orb}
Let $\lambda \in P_{+}^{0}$ be 
a level-zero dominant integral weight.
Then, every element $\nu$ in the $W$-orbit $W\lambda$ of 
$\lambda$ can be written uniquely as\,{\rm:} 
$\nu=\lambda-\beta+kd_{\lambda}\delta$ 
for some $\beta \in \fin{Q}_{+}$ and $k \in \BZ$, 
where $\fin{Q}_{+}:=\sum_{j \in I_{0}}\BZ_{\ge 0}\alpha_{j}$, 
and the positive integer $d_{\lambda} \in \BZ_{> 0}$ is defined 
by\,{\rm:} 
$W\lambda \cap (\lambda+\BZ\delta) = 
 \lambda+\BZ d_{\lambda}\delta$.
Furthermore, for the $\beta \in \fin{Q}_{+}$ above, 
there exists some $w \in \fin{W}$ such that 
$w\lambda=\lambda-\beta$.
\end{lem}

Let $\lambda \in P$ be an integral weight. 
For $\mu,\nu \in W\lambda$, we write $\mu > \nu$ if there exist 
a sequence $\mu=\nu_{0},\,\nu_{1},\,\dots,\,\nu_{n}=\nu$ 
of elements of $W\lambda$ and a sequence 
$\xi_{1},\,\xi_{2},\,\dots,\xi_{n}$ 
of positive real roots such that 
$\nu_{k}=r_{\xi_{k}}(\nu_{k-1})$ and 
$\nu_{k-1}(\xi_{k}^{\vee}) \in \BZ_{< 0}$ 
for all $1 \le k \le n$, where 
$\xi_{k}^{\vee} \in \Fh$ denotes the dual root of 
$\xi_{k} \in \prr$, and 
$r_{\xi_{k}}$ denotes the associated reflection; 
we write $\mu \ge \nu$ if $\mu > \nu$ or $\mu=\nu$. 
%
%%%%%%%%%%%%%%%%%%
%%% rem:bruhat %%%
%%%%%%%%%%%%%%%%%%
%
\begin{rem} \label{rem:bruhat}
Let $\lambda \in P^{0}_{+}$ be 
a level-zero dominant integral weight, and 
let $\nu,\,\nu' \in W\lambda$ be such that $\nu > \nu'$. 
Write $\nu$ and $\nu'$ as 
$\nu=\lambda-\beta+kd_{\lambda}\delta$ and 
$\nu'=\lambda-\beta'+k'd_{\lambda}\delta$ 
for $\beta,\,\beta \in \fin{Q}_{+}$ and $k,\,k' \in \BZ$ 
(see Lemma~\ref{lem:orb}), respectively. 
Then we deduce from \cite[Remark~2.4.3\,(1)]{NSls} that 
either $k < k'$ holds, or $k=k'$ and $\beta-\beta' \in 
\fin{Q}_{+} \setminus \{0\}$ holds.
\end{rem}

Now, let 
$\cl:\Fh^{\ast} \twoheadrightarrow \Fh^{\ast}/\BQ\delta$ 
be the canonical projection from $\Fh^{\ast}$ onto 
$\Fh^{\ast}/\BQ\delta$, and 
define a classical weight lattice 
$P_{\cl}$ and a classical dual weight lattice 
$P_{\cl}^{\vee}$ by:
%
%%%%%%%%%%%%%%%%%
%%% eq:lat-cl %%%
%%%%%%%%%%%%%%%%%
%
\begin{equation} \label{eq:lat-cl}
P_{\cl} = \cl(P) = 
 \bigoplus_{j \in I} \BZ \cl(\Lambda_{j})
\quad \text{and} \quad 
P_{\cl}^{\vee} = 
 \bigoplus_{j \in I} \BZ h_{j} 
 \subset P^{\vee}.
\end{equation}
We see that 
$P_{\cl} \simeq P/(\BQ\delta \cap P)$, and that 
$P_{\cl}$ can be identified with 
$\Hom_{\BZ}(P_{\cl}^{\vee},\BZ)$ as a $\BZ$-module by: 
$(\cl(\lambda))(h)=\lambda(h)$ 
for $\lambda \in P$ and $h \in P_{\cl}^{\vee}$.  
The quintuple $(A, \cl(\Pi), \Pi^{\vee}, 
P_{\cl}, P_{\cl}^{\vee})$ is 
called a classical Cartan datum 
(see \cite[\S10.1]{HK}). 
Note that there exists 
a natural action of the Weyl group $W$ on 
$\Fh^{\ast}/\BQ\delta$ induced 
from the one on $\Fh^{\ast}$, 
since $W\delta=\delta$.
It is obvious that 
$w \circ \cl = \cl \circ w$ for all $w \in W$.
If we set $\Fh^{\ast 0}:=
\bigl\{\lambda \in \Fh^{\ast} \mid \lambda(c)=0\bigr\}$, 
then there exists a (positive definite) symmetric bilinear form 
$(\cdot\,,\,\cdot)_{\cl}$ on 
$\cl(\Fh^{\ast 0})=\Fh^{\ast 0}/\BQ\delta$ 
induced from the restriction to $\Fh^{\ast 0}$ of 
the standard bilinear form $(\cdot\,,\,\cdot)$ 
on $\Fh^{\ast}$, since $(\delta,\Fh^{\ast 0})=\{0\}$. 
%
%
%%%%%%%%%%%%%%%%%%
%%% dfn:lv0-cl %%%
%%%%%%%%%%%%%%%%%%
%
\begin{dfn} \label{dfn:lv0-cl}
An integral weight $\mu \in P_{\cl}$ is said to be level zero 
if $\mu(c)=0$. A level-zero integral weight $\mu \in P_{\cl}$ is 
said to be level-zero dominant 
(resp., strictly level-zero dominant) if 
$\lambda(h_{j}) \ge 0$ (resp., $\lambda(h_{j}) > 0$) 
for all $j \in I_{0}$. 
\end{dfn}
%
%%%%%%%%%%%%%%%%%
%%% rem:orbcl %%%
%%%%%%%%%%%%%%%%%
%
\begin{rem} \label{rem:orbcl}
Let $\lambda \in P^{0}_{+}$ 
be a level-zero dominant integral weight. 
It is easily to verify that 
$\cl(W\lambda)=\fin{W}\cl(\lambda)$ 
(see the proof of \cite[Lemma~2.3.3]{NSls}).
Also, we see that $\cl(\lambda)$ is 
the unique level-zero dominant integral weight in 
$\cl(W\lambda)=\fin{W} \cl(\lambda)$, and that 
$\cl(w_{0}\lambda)=w_{0}\cl(\lambda)$ is 
the unique element of 
$\cl(W\lambda)=\fin{W} \cl(\lambda)$ such that 
$(\cl(w_{0}\lambda))(h_{j})=
 (w_{0}\cl(\lambda))(h_{j}) \le 0$ 
for all $j \in I_{0}$.
\end{rem}

Let $U_{q}'(\Fg)$ be the quantized universal enveloping algebra 
of the affine Lie algebra $\Fg$ with weight lattice $P_{\cl}$ 
over the field $\BQ(q)$ of rational functions in $q$. 
We denote by $x_{j}$, $y_{j}$, 
$j \in I$, and $q^{h}$, $h \in P^{\vee}_{\cl}$, 
the canonical generators of $U_{q}'(\Fg)$, 
where $x_{j}$ (resp., $y_{j}$) 
corresponds to the simple root $\alpha_{j}$ 
(resp., $-\alpha_{j}$) for $j \in I$. 
%
%
%
%==============================%
%     START SUBSECTION 0202    %
%==============================%
%
\subsection{Crystals of LS paths with weight lattice $P$, or $P_{\cl}$.}
\label{subsec:LSpath}
A path (with weight in $P$) is, by definition, 
a piecewise linear, continuous map 
$\pi:[0,1] \rightarrow \Fh^{\ast}_{\BR}:=
\BR \otimes_{\BQ} \Fh^{\ast}$ such that 
$\pi(0)=0$ and $\pi(1) \in P \subset 
 \BR \otimes_{\BZ} P=\Fh^{\ast}_{\BR}$. 
We denote by $\BP$ the set of all paths 
$\pi:[0,1] \rightarrow \Fh^{\ast}_{\BR}$. 
For each $\pi_{1},\,\pi_{2} \in \BP$, we define a path 
$\pi_{1} \pm \pi_{2} \in \BP$ by: 
$(\pi_{1} \pm \pi_{2})(t)=\pi_{1}(t) \pm \pi_{2}(t)$ 
for $t \in [0,1]$. 
For an integral weight $\nu \in P$, 
let $\pi_{\nu}$ denote
the straight line connecting $0 \in P$ 
with $\nu \in P$, i.e., 
$\pi_{\nu}(t):=t\nu$ for $t \in [0,1]$. 

Let $\pi \in \BP$. 
A pair $(\ud{\nu}\,;\,\ud{\sigma})$ of 
a sequence $\ud{\nu}:
 \nu_{1},\,\nu_{2},\,\dots,\,\nu_{s}$ of 
elements of $\Fh^{\ast}_{\BR}$ and 
a sequence $\ud{\sigma}:
 0=\sigma_{0} < \sigma_{1} < \cdots < \sigma_{s}=1$ 
of rational numbers is called an expression of 
$\pi \in \BP$ if the following equation holds:
%
%%%%%%%%%%%%%%%%%%%
%%%%% eq:path %%%%%
%%%%%%%%%%%%%%%%%%%
%
\begin{equation} \label{eq:path}
\pi(t)=\sum_{u'=1}^{u-1}
(\sigma_{u'}-\sigma_{u'-1})\nu_{u^{\prime}}+
(t-\sigma_{u-1})\nu_{u} 
\quad
\text{for \, }
   \sigma_{u-1} \le t \le \sigma_{u}, \ 
   1 \le u \le s. 
\end{equation}
In this case, 
we write $\pi=(\ud{\nu}\,;\,\ud{\sigma})$.
An expression 
$(\nu_{1},\,\nu_{2},\,\dots,\,\nu_{s} \,;\, \ud{\sigma})$ of $\pi$ 
is said to be reduced 
if $\nu_{u} \ne \nu_{u+1}$ for any $u=1,\,2,\,\dots,\,s-1$. 
%
%%%%%%%%%%%%%%%%
%%% rem:expr %%%
%%%%%%%%%%%%%%%%
%
\begin{rem}[{see \cite[Remark~2.5.2]{NSls}}] \label{rem:expr}
Let $\pi \in \BP$.
We easily see that there exists a unique 
reduced expression of $\pi$. 
Also, if 
$(\nu_{1},\,\nu_{2},\,\dots,\,\nu_{s} \,;\, 
  \sigma_{0},\,\sigma_{1},\,\dots,\,\sigma_{s})$ is
an expression of $\pi$, 
then the reduced expression of $\pi$ is obtained 
from this expression 
by ``omitting'' $\nu_{u}$'s 
such that $\nu_{u}=\nu_{u+1}$ and 
corresponding $\sigma_{u}$'s . 
\end{rem}
%
%%%%%%%%%%%%%%%
%%% dfn:dir %%%
%%%%%%%%%%%%%%%
%
\begin{dfn} \label{dfn:dir}
Let $\pi=
(\nu_{1},\,\nu_{2},\,\dots,\,\nu_{s}\,;\,\ud{\sigma})$ 
be an expression of $\pi \in \BP$. 
We call $\nu_{1} \in \Fh^{\ast}_{\BR}$ 
(resp., $\nu_{s} \in \Fh^{\ast}_{\BR}$) 
the initial (resp., final) direction of $\pi$; 
it is easy to check that 
these elements $\nu_{1},\,\nu_{s} \in \Fh^{\ast}_{\BR}$ 
do not depend on the choice of an expression of $\pi$. 
The initial (resp., final) direction of 
$\pi$ is denoted by $\iota(\pi)$ (resp., $\kappa(\pi)$). 
\end{dfn}
%
%%%%%%%%%%%%%%%
%%% rem:sum %%%
%%%%%%%%%%%%%%%
%
\begin{rem} \label{rem:sum}
Let $\pi_{1},\,\pi_{2} \in \BP$. 
We easily see that 
$\iota(\pi_{1} \pm \pi_{2}) = 
 \iota(\pi_{1}) \pm \iota(\pi_{2})$ and 
$\kappa(\pi_{1} \pm \pi_{2}) = 
 \kappa(\pi_{1}) \pm \kappa(\pi_{2})$.
\end{rem}

Let $\lambda \in P$ be an integral weight.
Recall from \cite[Definition~2.6.1]{NSls} 
(see also \cite[\S4]{L2}) that a Lakshmibai-Seshadri path 
(LS path for short) of shape $\lambda$ is a path $\pi \in \BP$ 
having an expression of the form 
$\pi=(\nu_{1},\,\nu_{2},\,\dots,\,\nu_{s}\,;\,
\sigma_{0},\,\sigma_{1},\,\dots,\,\sigma_{s})$, 
where $\nu_{1},\,\nu_{2},\,\dots,\,\nu_{s} \in W\lambda$, and 
where for each $1 \le u \le s-1$, 
there exists a ``$\sigma_{u}$-chain''
for $(\nu_{u},\,\nu_{u+1})$ 
(see \cite[\S4]{L2} and \cite[Definition~2.4.5]{NSls}
 for the definition of ``$\sigma_{u}$-chain'').
We denote by $\BB(\lambda)$ the set of 
all LS paths of shape $\lambda$. 
%
%%%%%%%%%%%%%%%
%%% rem:str %%%
%%%%%%%%%%%%%%%
%
\begin{rem}[{see \cite[Remark~2.6.2]{NSls}}] \label{rem:str}
Let $\lambda \in P$ be an integral weight. 
The straight line 
$\pi_{\nu}(t):=t\nu$, $t \in [0,1]$, is 
an element of $\BB(\lambda)$ 
for all $\nu \in W\lambda$.
\end{rem}
%
%%%%%%%%%%%%%%%%%
%%% lem:delta %%%
%%%%%%%%%%%%%%%%%
%
\begin{lem} \label{lem:delta}
Let $\lambda \in P^{0}_{+}$ be 
a level-zero dominant integral weight. 
If $\pi \in \BB(\lambda)$, then 
$\pi+\pi_{kd_{\lambda}\delta} \in \BB(\lambda)$ 
for all $k \in \BZ$. 
\end{lem}

\begin{proof}
Let $k \in \BZ$. 
We know from \cite[Lemma~2.7.4]{NSls} 
that $\pi+\pi_{kd_{\lambda}\delta} \in 
\BB(\lambda+kd_{\lambda}\delta)$.
Because 
there exists $w \in W$ such that $w(\lambda)=
\lambda+kd_{\lambda}\delta$ 
by the definition of $d_{\lambda}$, 
we deduce from the definition of LS paths that 
$\BB(\lambda+kd_{\lambda}\delta)=
 \BB(w(\lambda))=\BB(\lambda)$ 
(see \cite[Remark~2.6.3\,(3)]{NSls}). 
This proves the lemma. 
\end{proof}

Following \cite[\S1.2 and \S1.3]{L1} and \cite[\S1]{L2}
(see also \cite[\S5.1]{GL}), we recall the definition of 
the root operators 
$e_{j}$ and $f_{j}$, $j \in I$, for $\BB(\lambda)$. 
Let $\pi \in \BB(\lambda)$, and $j \in I$. We set
%
%%%%%%%%%%%%%
%%% eq:hm %%%
%%%%%%%%%%%%%
%
\begin{equation} \label{eq:hm}
H^{\pi}_{j}(t):=(\pi(t))(h_{j}) 
  \quad \text{for \,} t \in [0,1], \qquad
m^{\pi}_{j}
  :=\min\bigl\{H^{\pi}_{j}(t) \mid t \in [0,1]\bigr\}.
\end{equation}
Then we define $e_{j}\pi$ as follows 
(note that $m^{\pi}_{j} \in \BZ_{\le 0}$ 
 by \cite[Lemma~4.5\,d)]{L2}). 
If $m^{\pi}_{j} = 0$, then $e_{j}\pi:=\bzero$, 
where the $\bzero$ is an additional element 
corresponding to ``$0$''
in the theory of crystals. 
If $m^{\pi}_{j} \le -1$, then 
we define $e_{j}\pi \in \BP$ by:
%
%%%%%%%%%%%%%%%
%%% eq:ro_e %%%
%%%%%%%%%%%%%%%
%
\begin{equation} \label{eq:ro_e}
(e_{j}\pi)(t)=
\begin{cases}
\pi(t) & \text{if \,} 0 \le t \le t_{0}, \\[2mm]
\pi(t_{0})+r_{j}(\pi(t)-\pi(t_{0}))
       & \text{if \,} t_{0} \le t \le t_{1}, \\[2mm]
\pi(t)+\alpha_{j} & \text{if \,} t_{1} \le t \le 1,
\end{cases}
\end{equation}
where we set 
%
%%%%%%%%%%%%
%%% eq:t %%%
%%%%%%%%%%%%
%
\begin{equation} \label{eq:t}
\begin{array}{l}
t_{1}:=\min\bigl\{t \in [0,1] \mid 
       H^{\pi}_{j}(t)=m^{\pi}_{j} \bigr\}, \\[2mm]
t_{0}:=\max\bigl\{t \in [0,t_{1}] \mid
       H^{\pi}_{j}(t) = m^{\pi}_{j}+1 \bigr\}. \\[2mm]
\end{array}
\end{equation}
Similarly, $f_{j}\pi \in \BP \cup \{\bzero\}$ is 
given as follows 
(note that $H^{\pi}_{j}(1)-m^{\pi}_{j} \in \BZ_{\ge 0}$ 
 by \cite[Lemma~4.5\,d)]{L2} and $\pi(1) \in P$).
If $H^{\pi}_{j}(1)-m^{\pi}_{j}=0$, 
then $f_{j}\pi:=\bzero$. 
If $H^{\pi}_{j}(1)-m^{\pi}_{j} \ge 1$, 
then we define $f_{j}\pi \in \BP$ by: 
%
%%%%%%%%%%%%%%%
%%% eq:ro_f %%%
%%%%%%%%%%%%%%%
%
\begin{equation} \label{eq:ro_f}
(f_{j}\pi)(t)=
\begin{cases}
\pi(t) & \text{if \,} 0 \le t \le t_{0}, \\[2mm]
\pi(t_{0})+r_{j}(\pi(t)-\pi(t_{0}))
       & \text{if \,} t_{0} \le t \le t_{1}, \\[2mm]
\pi(t)-\alpha_{j} & \text{if \,} t_{1} \le t \le 1,
\end{cases}
\end{equation}
where we set 
\begin{equation}
\begin{array}{l}
t_{0}:=\max\bigl\{t \in [0,1] \mid 
       H^{\pi}_{j}(t)=m^{\pi}_{j}\bigr\}, \\[2mm]
t_{1}:=\min\bigl\{t \in [t_{0},1] \mid
       H^{\pi}_{j}(t)=m^{\pi}_{j}+1\bigr\}. \\[2mm]
\end{array}
\end{equation}
Using the root operators $e_{j}$ and $f_{j}$, $j \in I$, 
we can endow the set $\BB(\lambda)$ of 
all LS paths of shape $\lambda \in P$ 
with a structure of ($P$-weighted) crystal, i.e.,
a structure of crystal associated to the Cartan datum 
$(A,P,P^{\vee},\Pi,\Pi^{\vee})$
(see \cite[\S\S2 and 4]{L2} and also \cite[Theorems~1.2.3 and 1.4.5]{NSp2}). 
In fact, it follows from \cite[Lemma~2.1\,c)]{L2} that 
%
%%%%%%%%%%%%%%%%%%%%%
%%% eq:ve & eq:vp %%%
%%%%%%%%%%%%%%%%%%%%%
%
\begin{align}
& 
-m^{\pi}_{j} = 
  \max\bigl\{l \in \BZ_{\ge 0} \mid 
  e_{j}^{l}\pi \ne \bzero \bigr\}, 
\label{eq:ve} \\
&
H^{\pi}_{j}(1)-m^{\pi}_{j} = 
  \max\bigl\{l \in \BZ_{\ge 0} \mid 
  f_{j}^{l}\pi \ne \bzero \bigr\},
\label{eq:vp}
\end{align}
and hence that $\ve_{j}(\pi)=-m^{\pi}_{j}$ and 
$\vp_{j}(\pi)=H^{\pi}_{j}(1)-m^{\pi}_{j}$.

Let $\lambda \in P$ be an integral weight. 
For $\pi \in \BB(\lambda)$, define 
a piecewise linear, continuous map 
$\cl(\pi):[0,1] \rightarrow \Fh^{\ast}_{\BR}/\BR\delta$ by:
$(\cl(\pi))(t)=\cl(\pi(t))$ 
for $t \in [0,1]$, where 
$\cl:\Fh^{\ast}_{\BR} \twoheadrightarrow 
\Fh^{\ast}_{\BR}/\BR\delta$ is 
the canonical projection.
We set $\BB(\lambda)_{\cl}:=
 \bigl\{\cl(\pi) \mid \pi \in \BB(\lambda)\bigr\}$.
%
%%%%%%%%%%%%%%%%%%
%%% rem:str-cl %%%
%%%%%%%%%%%%%%%%%%
%
\begin{rem} \label{rem:str-cl}
We see from Remark~\ref{rem:str} that
the straight line $\eta_{\mu}(t)=t\mu$, $t \in [0,1]$, 
is contained in $\BB(\lambda)_{\cl}$ 
for all $\mu \in \cl(W\lambda)=\fin{W}\cl(\lambda)$.
\end{rem}

An expression and a reduced expression of 
$\eta \in \BB(\lambda)_{\cl}$ are defined similarly to 
those for the case of $\BB(\lambda)$. 
In addition, for $\eta \in \BB(\lambda)_{\cl}$, we define 
the initial and final directions of $\eta$ 
(which do not depend on 
 the choice of an expression of $\eta$) 
as in Definition~\ref{dfn:dir}, 
and also denote the initial (resp., final) direction of 
$\eta \in \BB(\lambda)_{\cl}$ by $\iota(\eta)$ 
(resp., $\kappa(\eta)$). 
%
%%%%%%%%%%%%%%%%%%%
%%% rem:cl-expr %%%
%%%%%%%%%%%%%%%%%%%
%
\begin{rem} \label{rem:cl-expr}
Let $\pi=(\nu_{1},\,\nu_{2},\,\dots,\,\nu_{s} \,;\,
 \ud{\sigma})$ be an expression of $\pi \in \BB(\lambda)$. 
Then, $\cl(\pi) \in \BB(\lambda)_{\cl}$ has an expression 
$\cl(\pi)=
 (\cl(\nu_{1}),\,\cl(\nu_{2}),\,\dots,\,\cl(\nu_{s}) \,;\, 
  \ud{\sigma})$.
It follows that 
\begin{equation*}
\begin{cases}
\iota(\cl(\pi))=\cl(\iota(\pi)) \in 
 \cl(W\lambda)=\fin{W}\cl(\lambda) & \\[2mm]
\kappa(\cl(\pi))=\cl(\kappa(\pi)) \in 
 \cl(W\lambda)=\fin{W}\cl(\lambda)
\end{cases}
\quad \text{for every $\pi \in \BB(\lambda)$}.
\end{equation*}
\end{rem}

Let $\eta \in \BB(\lambda)_{\cl}$, and $j \in I$. We set 
$H^{\eta}_{j}(t):=(\eta(t))(h_{j})$, $t \in [0,1]$, and 
define $m^{\eta}_{j}$ to be the minimum of the function 
$H^{\eta}_{j}(t)$ on the interval $[0,1]$. 
It is obvious that 
for every $\pi \in \BB(\lambda)$ and $j \in I$, 
%
%%%%%%%%%%%%
%%% eq:H %%%
%%%%%%%%%%%%
%
\begin{equation} \label{eq:H}
\text{
$H^{\cl(\pi)}_{j}(t) = H^{\pi}_{j}(t)$ 
for all $t \in [0,1]$, and hence 
$m^{\cl(\pi)}_{j} = m^{\pi}_{j}$.}
\end{equation}
%
%%%%%%%%%%%%%%%%%
%%% rem:chain %%%
%%%%%%%%%%%%%%%%%
%
\begin{rem} \label{rem:chain}
Let $\eta \in \BB(\lambda)_{\cl}$, and $j \in I$. 
We see from \cite[Lemma~4.5\,d)]{L2} and 
\eqref{eq:H} that $\eta(1) \in P_{\cl}$, and that 
all local minimums of the function 
$H^{\eta}_{j}(t)$, $t \in [0,1]$, are integers. 
In particular, the minimum $m^{\eta}_{j}$ of 
the function $H^{\eta}_{j}(t)$ on the interval $[0,1]$ is 
a nonpositive integer, and 
$H^{\eta}_{j}(1)-m^{\eta}_{j}$ is 
a nonnegative integer.
\end{rem}

We define $e_{j}\eta, \, f_{j}\eta \in 
\BB(\lambda)_{\cl} \cup \{\bzero\}$ for 
$\eta \in \BB(\lambda)_{\cl}$ and $j \in I$ 
in the same way as in the case of $\BB(\lambda)$. 
We easily see from \eqref{eq:H} that 
%
%%%%%%%%%%%%%%%%
%%% eq:cl-ro %%%
%%%%%%%%%%%%%%%%
%
\begin{equation} \label{eq:cl-ro}
\cl(e_{j}\pi)=e_{j}\cl(\pi), \quad 
\cl(f_{j}\pi)=f_{j}\cl(\pi) \quad 
\text{for \,} \pi \in \BB(\lambda)
\text{\, and \,} j \in I,
\end{equation}
where $\cl(\bzero)$ is understood to be $\bzero$. 
%
%%%%%%%%%%%%%%%%%%%
%%% rem:initial %%%
%%%%%%%%%%%%%%%%%%%
%
\begin{rem} \label{rem:initial}
Let $\eta \in \BB(\lambda)_{\cl}$, and $j \in I$. 
It follows from the definition of 
the root operator $e_{j}$ that 
if $e_{j}\eta \ne \bzero$, then the initial direction 
$\iota(e_{j}\eta)$ is equal
either to $\iota(\eta)$ or to $r_{j}(\iota(\eta))$.
\end{rem}

We know from \cite[Theorem~2.4 and \S3.1]{NSz} that 
the set $\BB(\lambda)_{\cl}$ equipped 
with the root operators $e_{j}$ and $f_{j}$, $j \in I$, 
becomes a $P_{\cl}$-weighted crystal 
($P_{\cl}$-crystal for short), i.e., 
a crystal associated to 
the classical Cartan datum $(A, \cl(\Pi), \Pi^{\vee}, 
P_{\cl}, P_{\cl}^{\vee})$, and 
that the following equations hold:
%
%%%%%%%%%%%%%%%%%%%%%%%%%%%
%%% eq:ve-cl & eq:vp-cl %%%
%%%%%%%%%%%%%%%%%%%%%%%%%%%
%
\begin{align}
& 
-m^{\eta}_{j} = 
  \max\bigl\{l \in \BZ_{\ge 0} \mid 
  e_{j}^{l}\eta \ne \bzero \bigr\}=\ve_{j}(\eta), 
\label{eq:ve-cl} \\
&
H^{\eta}_{j}(1)-m^{\eta}_{j} = 
  \max\bigl\{l \in \BZ_{\ge 0} \mid 
  f_{j}^{l}\eta \ne \bzero \bigr\}=\vp_{j}(\eta). 
\label{eq:vp-cl}
\end{align}

For each $\eta \in \BB(\lambda)_{\cl}$ and $j \in I$, 
we set
$e_{j}^{\max}\eta:=
 e_{j}^{\ve_{j}(\eta)}\eta \in \BB(\lambda)_{\cl}$.
The proof of the next lemma is similar to that 
of \cite[5.3 Lemma]{L1} (see Remark~\ref{rem:chain}). 
%
%
%%%%%%%%%%%%%%%%%
%%% lem:init1 %%%
%%%%%%%%%%%%%%%%%
%
\begin{lem} \label{lem:init1}
Let $\lambda \in P$ be an integral weight, and 
let $\eta \in \BB(\lambda)_{\cl}$, $j \in I$.

\noindent
{\rm (1)} 
If the initial direction 
$\iota(\eta) \in P_{\cl}$ of $\eta$
satisfies $(\iota(\eta))(h_{j}) < 0$, then 
$e_{j}\eta \ne \bzero$. 

\noindent
{\rm (2)} For all $0 \le l \le \ve_{j}(\eta)-1$, 
we have $\iota(e_{j}^{l}\eta)=\iota(\eta)$. 

\noindent
{\rm (3)} If $(\iota(\eta))(h_{j}) \le 0$, then 
we have $\iota(e_{j}^{\max}\eta)=r_{j}(\iota(\eta))$.
\end{lem}

Using Lemma~\ref{lem:init1}\,(3), 
we can show the following lemma by 
induction on $p$. 
%
%%%%%%%%%%%%%%%%%
%%% lem:init2 %%%
%%%%%%%%%%%%%%%%%
%
\begin{lem} \label{lem:init2}
Let $\lambda \in P$ be an integral weight. 
Let $\eta \in \BB(\lambda)_{\cl}$, and 
set $\mu:=\iota(\eta) \in P_{\cl}$. 
If $j_{1},\,j_{2},\,\dots,\,j_{p} \in I$ 
satisfy the condition that
$\bigl(r_{j_{p'}}r_{j_{p'-1}} \cdots r_{j_{1}}(\mu)\bigr)
(h_{j_{p'+1}}) \le 0$ for all $p'=0,\,1,\,\dots,\,p-1$, 
then the initial direction of 
$e_{j_{p}}^{\max}e_{j_{p-1}}^{\max} \cdots 
 e_{j_{1}}^{\max}(\eta) \in \BB(\lambda)_{\cl}$ 
is equal to $r_{j_{p}}r_{j_{p-1}} \cdots r_{j_{1}}(\mu)$. 
\end{lem}
%
%%%%%%%%%%%%%%%%%
%%% lem:final %%%
%%%%%%%%%%%%%%%%%
%
\begin{lem} \label{lem:final}
Let $\lambda \in P$ be an integral weight. 
Let $\eta \in \BB(\lambda)_{\cl}$, and set 
$\mu:=\kappa(\eta) \in P_{\cl}$. 
If $j \in I$ satisfies $\mu(h_{j}) > 0$, 
then $f_{j}\eta \ne \bzero$ holds.
\end{lem}

\begin{proof}
From the assumption of the lemma, we deduce that 
$H^{\pi}_{j}(1) - m^{\pi}_{j} > 0$. 
It follows from Remark~\ref{rem:chain} that 
$H^{\pi}_{j}(1) - m^{\pi}_{j} \ge 1$. 
Hence, by \eqref{eq:vp-cl}, 
we have $\vp_{j}(\eta) \ge 1$, 
which implies that $f_{j}\eta \ne \bzero$. 
This proves the lemma. 
\end{proof}
%
%
%
%==============================%
%     START SUBSECTION 0203    %
%==============================%
%
\subsection{Regular crystals and simple crystals.}
\label{subsec:simple}
For a proper subset $J$ of $I$, we set
$A_{J}:=(a_{ij})_{i,j \in J}$, 
$\Pi_{J}:=
\bigl\{\alpha_{j}\bigr\}_{j \in J} 
\subset \Pi$, 
and 
$\Pi_{J}^{\vee}:=
\bigl\{h_{j}\bigr\}_{j \in J} 
\subset \Pi^{\vee}$. 
When we regard a $P_{\cl}$-crystal $\CB$ 
as a crystal associated to the Cartan datum 
$(A_{J},\cl(\Pi_{J}),\Pi_{J}^{\vee},P_{\cl},P_{\cl}^{\vee})$, 
we denote it by $\res_{J}\CB$. 
Also, we denote by $U_{q}'(\Fg)_{J}$ 
the $\BQ(q)$-subalgebra of $U_{q}'(\Fg)$ 
generated by $x_{j}$, $y_{j}$, $j \in J$, and 
$q^{h}$, $h \in P_{\cl}^{\vee}$. 
Recall that 
a $P_{\cl}$-crystal $\CB$ is said to be regular 
if for every proper subset $J \subsetneq I$, 
$\res_{J} \CB$ is isomorphic to the crystal base of 
an integrable $U_{q}^{\prime}(\Fg)_{J}$-module. 
%
%%%%%%%%%%%%%%%%%%%
%%% rem:regular %%%
%%%%%%%%%%%%%%%%%%%
%
\begin{rem} \label{rem:regular}
Let $\lambda \in P$ be a level-zero integral weight. 
We know from \cite[Proposition~3.13]{NSz} that 
$\BB(\lambda)_{\cl}$ is a regular $P_{\cl}$-crystal 
with finitely many elements.
\end{rem}

If $\CB$ is a regular $P_{\cl}$-crystal 
with Kashiwara operators $e_{j}$ and $f_{j}$, $j \in I$, 
then we set $e_{j}^{\max}b:=e_{j}^{\ve_{j}(b)}b \in \CB$ 
for $b \in \CB$ and $j \in I$, where 
$\ve_{j}(b):=\max \bigl\{l \in \BZ_{\ge 0} \mid 
e_{j}^{l}b \ne \bzero \bigr\}$. 
For regular $P_{\cl}$-crystals 
$\CB_{1}$ and $\CB_{2}$, 
we define the tensor product $P_{\cl}$-crystal 
$\CB_{1} \otimes \CB_{2}$ 
of $\CB_{1}$ and $\CB_{2}$ 
as in \cite[\S7.3]{Kas-onc} and 
\cite[Definition~4.5.3]{HK}; 
note that $\CB_{1} \otimes \CB_{2}$ is also 
a regular $P_{\cl}$-crystal.
The next lemma follows immediately 
from the tensor product rule for crystals. 
%
%%%%%%%%%%%%%%%
%%% lem:max %%%
%%%%%%%%%%%%%%%
%
\begin{lem} \label{lem:max}
Let $\CB_{1}$ and $\CB_{2}$ be 
regular $P_{\cl}$-crystals, and 
$b_{1} \in \CB_{1}$, 
$b_{2} \in \CB_{2}$. 
Let $j \in I$. 

\noindent {\rm (1)} \, We have
$\ve_{j}(b_{1} \otimes b_{2}) \ge \ve_{j}(b_{1})$. 
Therefore, if $e_{j}(b_{1} \otimes b_{2}) = \bzero$, 
then $e_{j}b_{1}=\bzero$. 

\noindent {\rm (2)} \, 
Set $L:=\ve_{j}(b_{1} \otimes b_{2})$. 
Then, for $0 \le l \le L$, 
we have 
\begin{equation*}
e_{j}^{l}(b_{1} \otimes b_{2}) = 
\begin{cases}
b_{1} \otimes e_{j}^{l}b_{2} & 
\text{\rm if $0 \le l \le L-\ve_{j}(b_{1})$}, \\[3mm]
e_{j}^{l-L+\ve_{j}(b_{1})} b_{1} \otimes 
e_{j}^{L-\ve_{j}(b_{1})} b_{2} & 
\text{\rm if $L-\ve_{j}(b_{1}) \le l \le L$}.
\end{cases}
\end{equation*}
In particular, 
$e_{j}^{\max}(b_{1} \otimes b_{2})=
 e_{j}^{\max}b_{1} \otimes b_{2}'$ 
for some $b_{2}' \in \CB_{2}$.
\end{lem}

Let $\CB$ be a regular $P_{\cl}$-crystal.
We define 
%
%%%%%%%%%%%%%%%%%%%
%%% eq:def-norm %%%
%%%%%%%%%%%%%%%%%%%
%
\begin{equation} \label{eq:def-norm}
||b||:=\sqrt{(\wt(b),\wt(b))_{\cl}}
\qquad
\text{for $b \in \CB$}.
\end{equation}
%
%%%%%%%%%%%%%%%%%
%%% lem:e-max %%%
%%%%%%%%%%%%%%%%%
%
\begin{lem} \label{lem:e-max}
Let $\CB$ be a regular $P_{\cl}$-crystal.
For each $b \in \CB$ and $j \in I$, 
we have $||e^{\max}_{j}b|| \ge ||b||$, 
with equality if and only if 
either $e_{j}b=\bzero$ or $f_{j}b=\bzero$ holds. 
\end{lem}

\begin{proof}
Using the equation 
$\vp_{j}(b)=(\wt(b))(h_{j}) + \ve_{j}(b)$, 
we easily see that 
\begin{equation*}
||e^{\max}_{j}b||^{2} = 
||b||^{2} + 
    \ve_{j}(b) \vp_{j}(b)
    (\alpha_{j},\alpha_{j}).
\end{equation*}
The inequality $||e^{\max}_{j}b|| \ge ||b||$ 
follows immediately from the fact that 
$\ve_{j}(b) \ge 0$, $\vp_{j}(b) \ge 0$, and 
$(\alpha_{j},\alpha_{j}) > 0$. 
Also, the equality holds 
if and only if 
$\ve_{j}(b)=0$ or $\vp_{j}(b)=0$, which is equivalent to 
saying that $e_{j}b=\bzero$ or $f_{j}b=\bzero$. 
This proves the lemma. 
\end{proof}

Let $\CB$ be a regular $P_{\cl}$-crystal. 
For each $j \in I$, we define 
$S_{j}:\CB \rightarrow \CB$ by:
%
%%%%%%%%%%%%%
%%% eq:sj %%%
%%%%%%%%%%%%%
%
\begin{equation} \label{eq:sj}
S_{j}b=
\begin{cases}
f_{j}^{l}b & \text{if \ } l:=(\wt b)(h_{j}) \ge 0, \\[2mm]
e_{j}^{-l}b & \text{if \ } l:=(\wt b)(h_{j}) < 0.
\end{cases}
\end{equation}
We know from \cite[\S7]{Kas-mod} that 
there exists a unique action $S:W \rightarrow \Bij(\CB)$, 
$w \mapsto S_{w}$, of the Weyl group $W$ on the set $\CB$ 
such that $S_{r_{j}}=S_{j}$ for all $j \in I$, 
where $\Bij(\CB)$ denotes the group of 
all bijections from the set 
$\CB$ to itself; in fact, if 
$w=r_{j_{1}}r_{j_{2}} \cdots r_{j_{p}} \in W$ for 
$j_{1},\,j_{2},\,\dots,\,j_{p} \in I$, then 
$S_{w}=S_{j_{1}}S_{j_{2}} \cdots S_{j_{p}}$.
Note that $\wt(S_{w}b)=w(\wt(b))$ 
for all $w \in W$ and $b \in \CB$. 
%
%%%%%%%%%%%%%%%
%%% dfn:ext %%%
%%%%%%%%%%%%%%%
%
\begin{dfn}[{\cite[\S1.4]{AK}}] \label{dfn:ext}
Let $\CB$ be a regular $P_{\cl}$-crystal. 
An element $b \in \CB$ is said to be extremal 
if for every $w \in W$, either 
$e_{j}S_{w}b=\bzero$ or 
$f_{j}S_{w}b=\bzero$ holds for each $j \in I$. 
\end{dfn}
%
%%%%%%%%%%%%%%%
%%% rem:ext %%%
%%%%%%%%%%%%%%%
%
\begin{rem} \label{rem:ext}
It follows immediately from the definition above that 
if $b \in \CB$ is an extremal element, 
then $S_{w}b$ is an extremal element 
of weight $w(\wt(b))$ for each $w \in W$. 
\end{rem}
%
%
%%%%%%%%%%%%%%%%
%%% lem:norm %%%
%%%%%%%%%%%%%%%%
%
\begin{lem} \label{lem:norm}
Let $\CB$ be a regular $P_{\cl}$-crystal 
with finitely many elements.
If $b \in \CB$ satisfies the condition that 
$||b||=
 \max 
 \bigl\{
 ||b'|| \mid b' \in \BB(\lambda)_{\cl}
 \bigr\}$, then $b$ is an extremal element. 
\end{lem}

\begin{proof}
Let $w \in W$, and $j \in I$.
Since $(\cdot\,,\,\cdot)_{\cl}$ is 
$W$-invariant, 
it follows that $||S_{w}b|| = ||b||$.
Using this, we deduce that
\begin{align*}
||b||
 & \ge ||e_{j}^{\max}S_{w}b||
   \quad \text{by the maximality of $||b||$} \\
 & \ge ||S_{w}b||
   \quad \text{by Lemma~\ref{lem:e-max}} \\
 & = ||b||,
\end{align*}
and hence that
$||e_{j}^{\max}S_{w}b||=||S_{w}b||$. 
Therefore, by Lemma~\ref{lem:e-max}, 
either $e_{j}S_{w}b=\bzero$ or 
$f_{j}S_{w}b=\bzero$ holds.
This proves the lemma.
\end{proof}
%
%
%%%%%%%%%%%%%%%%%%
%%% dfn:simple %%%
%%%%%%%%%%%%%%%%%%
%
\begin{dfn} \label{dfn:simple}
Let $\CB$ be a regular $P_{\cl}$-crystal with 
finitely many elements. 
The $P_{\cl}$-crystal $\CB$ is said to be simple 
if it satisfies the following conditions: 

\noindent
{\rm (1)} The weights of elements of $\CB$ are all of level zero. 

\noindent 
{\rm (2)} The set of all extremal elements of $\CB$ 
coincides with a $W$-orbit in $\CB$. In addition, 
for each extremal element $b \in \CB$, 
the subset $\CB_{\wt(b)} \subset \CB$ of all elements of 
weight $\wt(b)$ consists of a single element, i.e., 
$\CB_{\wt(b)}=\bigl\{b\bigr\}$.
\end{dfn}
%
%%%%%%%%%%%%%%%%%%
%%% rem:simple %%%
%%%%%%%%%%%%%%%%%%
%
\begin{rem} \label{rem:simple}
Let $\CB$ be a simple $P_{\cl}$-crystal.
Then it follows from Remark~\ref{rem:orbcl} and 
Definition~\ref{dfn:simple}\,(2) that there exists a unique 
extremal element $b$ of $\CB$ such that 
$\wt(b) \in P_{\cl}$ is level-zero dominant.
\end{rem}
%
%%%%%%%%%%%%%%%%%%
%%% lem:simple %%%
%%%%%%%%%%%%%%%%%%
%
\begin{lem} \label{lem:simple}
{\rm (1)} A simple $P_{\cl}$-crystal is connected. 

\noindent 
{\rm (2)} A tensor product of simple $P_{\cl}$-crystals is also 
a simple $P_{\cl}$-crystal.

\noindent
{\rm (3)} 
Let $\CB_{1}$, $\CB_{2}$ be simple $P_{\cl}$-crystals. 
Then there exists at most one isomorphism of 
$P_{\cl}$-crystals from $\CB_{1}$ to $\CB_{2}$. 
In particular, any automorphism of 
a simple $P_{\cl}$-crystal is 
necessarily the identity map.
\end{lem}

\begin{proof}
Parts~(1) and (2) are simply \cite[Lemmas~1.9 and 1.10]{AK}, 
respectively. Let us show part~(3). 
Let $\CB_{1}$, $\CB_{2}$ be simple $P_{\cl}$-crystals, and 
let $\Phi:\CB_{1} \stackrel{\sim}{\rightarrow} \CB_{2}$ be 
an isomorphism of $P_{\cl}$-crystals 
from $\CB_{1}$ to $\CB_{2}$.
We see from Remark~\ref{rem:simple} that 
there exists a unique element 
$b_{1} \in \CB_{1}$ (resp., $b_{2} \in \CB_{2}$) such that 
$b_{1}$ (resp., $b_{2}$) is extremal, and 
$\wt(b_{1})$ (resp., $\wt(b_{2})$) is 
level-zero dominant. 
Since $\Phi:\CB_{1} \stackrel{\sim}{\rightarrow} \CB_{2}$ is 
an isomorphism of $P_{\cl}$-crystals, it follows that 
$\Phi(b_{1})$ is extremal, and 
$\wt(\Phi(b_{1}))$ is level-zero dominant. 
Hence, from the uniqueness of such an element, 
we obtain $\Phi(b_{1})=b_{2}$. 
But, because a simple $P_{\cl}$-crystal is 
connected by part (1), 
an isomorphism of $P_{\cl}$-crystals from 
$\CB_{1}$ to $\CB_{2}$ is determined uniquely
by the requirement that $\Phi(b_{1})=b_{2}$. 
Thus the proof of the lemma is complete. 
\end{proof}
%
%
%
%==============================%
%     START SUBSECTION 0204    %
%==============================%
%
\subsection{
 Tensor product decomposition and combinatorial $R$-matrices.}
\label{subsec:tpd}
We know from \cite[Propositions~3.4.1 and 3.4.2]{NSp2} that 
for each $i \in I_{0}$, 
the $P_{\cl}$-crystal $\BB(\vpi_{i})_{\cl}$ is 
a simple $P_{\cl}$-crystal isomorphic to 
the crystal basis $\CB(\vpi_{i})$ of 
the level-zero fundamental representation $W(\vpi_{i})$, 
which is a finite-dimensional irreducible $U_{q}'(\Fg)$-module 
introduced in \cite[\S5.2]{Kas-lz}. 
Because the $W(\vpi_{i})$, $i \in I_{0}$, are
``good'' $U_{q}'(\Fg)$-modules in the sense 
of \cite[\S8]{Kas-lz}, 
we deduce from \cite[Proposition~10.6]{Kas-lz} that 
for each $i_{1},\,i_{2} \in I_{0}$, 
there exists a unique isomorphism 
(called a combinatorial $R$-matrix)
$R_{\vpi_{i_{1}},\vpi_{i_{2}}}: 
 \BB(\vpi_{i_{1}})_{\cl} \otimes \BB(\vpi_{i_{2}})_{\cl}
 \stackrel{\sim}{\rightarrow}
 \BB(\vpi_{i_{2}})_{\cl} \otimes \BB(\vpi_{i_{1}})_{\cl}$ of 
$P_{\cl}$-crystals  (see also \cite[\S2.3]{O}); 
the uniqueness 
follows from Lemma~\ref{lem:simple}\,(2),(3).
Combining this fact and 
the tensor product decomposition theorem 
\cite[Theorem~3.2]{NSz}, 
we obtain the following theorem. 
%
%%%%%%%%%%%%%%%
%%% thm:tpd %%%
%%%%%%%%%%%%%%%
%
\begin{thm} \label{thm:tpd}
Let $\bi=(i_{1},\,i_{2},\,\dots,\,i_{n})$ be 
an arbitrary sequence of elements of $I_{0}$  
{\rm(}with repetitions allowed\,{\rm)}, and 
set $\lambda=\sum_{k=1}^{n} \vpi_{i_{k}} \in P^{0}_{+}$. 
Then there exists a unique 
isomorphism of $P_{\cl}$-crystals 
\begin{equation}
\Psi_{\bi}:
\BB(\lambda)_{\cl} \stackrel{\sim}{\rightarrow} 
\BB_{\bi}:=
\BB(\vpi_{i_{1}})_{\cl} \otimes 
\BB(\vpi_{i_{2}})_{\cl} \otimes \cdots \otimes 
\BB(\vpi_{i_{n}})_{\cl}.
\end{equation}
\end{thm}
%
%%%%%%%%%%%%%%%%%
%%% rem:tpd01 %%%
%%%%%%%%%%%%%%%%%
%
\begin{rem} \label{rem:tpd01}
Let $\bi=(i_{1},\,i_{2},\,\dots,\,i_{n})$ and 
$\lambda=\sum_{k=1}^{n}\vpi_{i_{k}} \in P^{0}_{+}$ 
be as in Theorem~\ref{thm:tpd}.
It follows from Theorem~\ref{thm:tpd} and 
Lemma~\ref{lem:simple}\,(1), (2) that 
$\BB(\lambda)_{\cl}$ is a simple $P_{\cl}$-crystal 
isomorphic to the crystal basis of 
the tensor product $U_{q}'(\Fg)$-module 
$W_{\bi}:=W(\vpi_{i_{1}}) \otimes 
 W(\vpi_{i_{2}}) \otimes \cdots \otimes 
 W(\vpi_{i_{n}})$ of 
the level-zero fundamental representations
$W(\vpi_{i_{k}})$, $1 \le k \le n$.
\end{rem}
%
%%%%%%%%%%%%%%%%%
%%% rem:tpd02 %%%
%%%%%%%%%%%%%%%%%
%
\begin{rem} \label{rem:tpd02}
Let $\lambda \in P^{0}_{+}$.
We know from \cite[Lemma~3.19\,(1)]{NSz} that 
the straight line 
$\eta_{\cl(\lambda)}$ is an extremal element of 
$\BB(\lambda)_{\cl}$, and that 
$S_{w}\eta_{\cl(\lambda)}=\eta_{w\cl(\lambda)}$ 
for each $w \in W$.
Therefore, from Remark~\ref{rem:tpd01}, 
we deduce (recalling the definition of simple $P_{\cl}$-crystals) 
that each extremal element of $\BB(\lambda)_{\cl}$ 
is a straight line $\eta_{\mu}$ for some 
$\mu \in \cl(W\lambda)=\fin{W}\cl(\lambda)$, and 
that the number of elements of 
weight $\mu$ in 
$\BB(\lambda)_{\cl}$ is equal to $1$ 
for all $\mu \in \cl(W\lambda)=\fin{W}\cl(\lambda)$. 
In particular, 
the straight line $\eta_{\cl(\lambda)}$ is 
the unique extremal element of $\BB(\lambda)_{\cl}$ 
whose weight is level-zero dominant, 
which we mentioned in Remark~\ref{rem:simple}.
\end{rem}
%
%%%%%%%%%%%%%%%%%
%%% rem:tpd03 %%%
%%%%%%%%%%%%%%%%%
%
\begin{rem} \label{rem:tpd03}
Let $\lambda \in P^{0}_{+}$. 
We see from \cite[Lemma~2.6.4]{NSls} that 
the weights of $\BB(\lambda)_{\cl}$ are all contained in 
the set $\cl(\lambda)-a_{0}^{-1}\cl(\fin{Q}_{+})$. 
\end{rem}

Let $\bi=(i_{1},\,i_{2},\,\dots,\,i_{n})$ and 
$\lambda \in P^{0}_{+}$ be as in Theorem~\ref{thm:tpd}.
It is easily seen from 
Remark~\ref{rem:tpd02} and
\cite[Lemma~1.6\,(1)]{AK} that 
$\eta_{\cl(\vpi_{i_{1}})} \otimes 
 \eta_{\cl(\vpi_{i_{2}})} \otimes \cdots \otimes 
 \eta_{\cl(\vpi_{i_{n}})}$ is 
the unique extremal element of the simple $P_{\cl}$-crystal 
$\BB_{\bi}=\BB(\vpi_{i_{1}})_{\cl} \otimes 
 \BB(\vpi_{i_{2}})_{\cl} \otimes \cdots \otimes 
 \BB(\vpi_{i_{n}})_{\cl}$ 
whose weight is level-zero dominant. 
Therefore, we deduce from Remark~\ref{rem:tpd02} and 
the proof of Lemma~\ref{lem:simple}\,(3) that 
%
%%%%%%%%%%%%%%%
%%% eq:map1 %%%
%%%%%%%%%%%%%%%
%
\begin{equation} \label{eq:map1}
\Psi_{\bi}(\eta_{\cl(\lambda)}) = 
 \eta_{\cl(\vpi_{i_{1}})} \otimes 
 \eta_{\cl(\vpi_{i_{2}})} \otimes \cdots \otimes 
 \eta_{\cl(\vpi_{i_{n}})}.
\end{equation}
%
%%%%%%%%%%%%%%%
%%% cor:tpd %%%
%%%%%%%%%%%%%%%
%
\begin{cor} \label{cor:tpd}
Let $\lambda,\,\lambda' \in P^{0}_{+}$ be 
level-zero dominant integral weights.

\noindent 
{\rm (1)} There exists a unique isomorphism 
$\Psi_{\lambda,\lambda'}: 
\BB(\lambda+\lambda')_{\cl}
 \stackrel{\sim}{\rightarrow} 
\BB(\lambda)_{\cl} \otimes \BB(\lambda')_{\cl}$
of $P_{\cl}$-crystals.

\noindent
{\rm (2)} 
There exists a unique isomorphism 
$R_{\lambda,\lambda'}: 
\BB(\lambda)_{\cl} \otimes \BB(\lambda')_{\cl}
 \stackrel{\sim}{\rightarrow} 
\BB(\lambda')_{\cl} \otimes \BB(\lambda)_{\cl}$
of $P_{\cl}$-crystals.
\end{cor}

\begin{proof}
Part~(1) follows from Theorem~\ref{thm:tpd}. 
Part~(2) follows immediately from part~(1).
\end{proof}
Let $\lambda,\,\lambda' \in P^{0}_{+}$ be 
level-zero dominant integral weights.
By the same reasoning as that yielding \eqref{eq:map1}, 
we obtain 
%
%%%%%%%%%%%%%%%%%%%%%%%%
%%% eq:map2 and map3 %%%
%%%%%%%%%%%%%%%%%%%%%%%%
%
\begin{align} 
& \label{eq:map2}
\Psi_{\lambda,\lambda'}(\eta_{\cl(\lambda+\lambda')}) = 
 \eta_{\cl(\lambda)} \otimes 
 \eta_{\cl(\lambda')}, \\
& \label{eq:map3}
R_{\lambda,\lambda'}
(\eta_{\cl(\lambda)} \otimes \eta_{\cl(\lambda')}) = 
 \eta_{\cl(\lambda')} \otimes \eta_{\cl(\lambda)}.
\end{align}
%
%
%
%==============================%
%     START SUBSECTION 0205    %
%==============================%
%
\subsection{Local energy functions.}
\label{subsec:eng}
Let $\lambda,\,\lambda' \in P^{0}_{+}$ be 
level-zero dominant integral weights, and 
let $\bi=(i_{1},\,i_{2},\,\dots,\,i_{n})$, 
$\bi'=(i_{1}',\,i_{2}',\,\dots,\,i_{n'}')$ be 
sequences of elements of $I_{0}$ such that 
$\lambda=\sum_{k=1}^{n}\vpi_{i_{k}}$ and 
$\lambda'=\sum_{k'=1}^{n'}\vpi_{i_{k'}'}$, respectively. 
We define the tensor product $U_{q}'(\Fg)$-modules 
$W_{\bi}$ and $W_{\bi'}$ corresponding to 
$\bi$ and $\bi'$, respectively, as in Remark~\ref{rem:tpd01};
note that both $W_{\bi}$ and $W_{\bi'}$ are 
``good'' $U_{q}'(\Fg)$-modules (in the sense of \cite[\S8]{Kas-lz})
by \cite[Proposition~8.7]{Kas-lz}, and that 
$\BB(\lambda)_{\cl}$ and $\BB(\lambda')_{\cl}$ are 
isomorphic as a $P_{\cl}$-crystal to the crystal bases of 
$W_{\bi}$ and $W_{\bi'}$, respectively.
Therefore, by an argument similar to that in \cite[\S11]{Kas-lz}, 
we obtain the following theorem (see also \cite[\S2.3]{O}). 
%
%
%%%%%%%%%%%%%%%%%%
%%% thm:energy %%%
%%%%%%%%%%%%%%%%%%
%
\begin{thm} \label{thm:energy}
Let $\lambda,\,\lambda' \in P^{0}_{+}$ be 
level-zero dominant integral weights.
Then, there exists a unique $\BZ$-valued function 
{\rm(}called a local energy function{\rm)}
$H_{\lambda,\lambda'}:
 \BB(\lambda)_{\cl} \otimes \BB(\lambda')_{\cl} \rightarrow \BZ$ 
satisfying the conditions\,{\rm:}

\noindent {\rm (H1)} \, 
For each $\eta_{1} \otimes \eta_{2} \in 
\BB(\lambda)_{\cl} \otimes \BB(\lambda')_{\cl}$ and $j \in I$ 
such that $e_{j}(\eta_{1} \otimes \eta_{2}) \ne \bzero$, 
the following equation holds.
%
%%%%%%%%%%%%%%%%%%%
%%% eq:energy-e %%%
%%%%%%%%%%%%%%%%%%%
%
\begin{align}
& H_{\lambda,\lambda'}(e_{j}(\eta_{1} \otimes \eta_{2})) = 
  \nonumber \\[3mm]
& \hspace*{10mm}
\begin{cases}
H_{\lambda,\lambda'}(\eta_{1} \otimes \eta_{2})+1 & \\[1mm]
\hspace*{5mm} \text{\rm if $j=0$, and if
   $e_{0}(\eta_{1} \otimes \eta_{2})=
    e_{0}\eta_{1} \otimes \eta_{2}$,
   $e_{0}(\ti{\eta}_{2} \otimes \ti{\eta}_{1})=
    e_{0}\ti{\eta}_{2} \otimes \ti{\eta}_{1}$,} & \\[3mm]
H_{\lambda,\lambda'}(\eta_{1} \otimes \eta_{2})-1 & \\[1mm]
\hspace*{5mm} \text{\rm if $j=0$, and if
   $e_{0}(\eta_{1} \otimes \eta_{2})=
    \eta_{1} \otimes e_{0}\eta_{2}$, 
   $e_{0}(\ti{\eta}_{2} \otimes \ti{\eta}_{1})=
    \ti{\eta}_{2} \otimes e_{0}\ti{\eta}_{1}$,} & \\[3mm]
H_{\lambda,\lambda'}(\eta_{1} \otimes \eta_{2}) 
\hspace{10mm} \text{\rm otherwise}, & 
\end{cases}
\end{align}
where we set 
$\ti{\eta}_{2} \otimes \ti{\eta}_{1}:=
 R_{\lambda,\lambda'}(\eta_{1} \otimes \eta_{2}) \in 
 \BB(\lambda')_{\cl} \otimes \BB(\lambda)_{\cl}$.

\noindent {\rm (H2)} \, 
$H_{\lambda,\lambda'}
 (\eta_{\cl(\lambda)} \otimes \eta_{\cl(\lambda')})=0$.
\end{thm}

%=========================%
%     START SECTION 03    %
%=========================%
%
\section{Degree functions on LS path crystals.}
\label{sec:degree}

Throughout this section, we fix 
a level-zero dominant integral weight $\lambda \in P$
contained in the set $P^{0}_{+}=
\sum_{i \in I_{0}} \BZ_{\ge 0}\vpi_{i}$.
%

%==============================%
%     START SUBSECTION 0301    %
%==============================%
%
\subsection{Definition of degree functions.}
\label{subsec:dfn-deg}

Let $\lambda \in P^{0}_{+}$.
Recall from Lemma~\ref{lem:orb} that 
every element $\nu$ of $W\lambda$ 
can be written uniquely in the form 
$\nu=\lambda-\beta+kd_{\lambda}\delta$ 
with $\beta \in \fin{Q}_{+}$ and $k \in \BZ$.
%
%
%%%%%%%%%%%%%%%
%%% lem:end %%%
%%%%%%%%%%%%%%%
%
\begin{lem} \label{lem:end}
Let $\pi \in \BB(\lambda)$. If the initial direction 
$\iota(\pi) \in W\lambda$ of $\pi$ is contained 
in $\lambda-\fin{Q}_{+}$, then $\pi(1) \in P$ 
can be written uniquely in the form 
$\pi(1)=\lambda-a_{0}^{-1}\beta+a_{0}^{-1}K\delta$ 
with $\beta \in \fin{Q}_{+}$ and $K \in \BZ_{\ge 0}$.
\end{lem}

\begin{proof}
We know from \cite[Lemma~2.6.4]{NSls}, along with the linear 
independence of $\alpha_{j}$, $j \in I$, and $\delta$, that 
$\pi(1) \in P$ can be written uniquely in the form
$\pi(1)=\lambda-a_{0}^{-1}\beta+a_{0}^{-1}K\delta$ 
with $\beta \in \fin{Q}_{+}$ and $K \in \BZ$.
Let us show that the coefficient $a_{0}^{-1}K$ of 
$\delta$ in this expression of $\pi(1) \in P$ is nonnegative. 
Let $\pi=(\nu_{1},\,\nu_{2},\,\dots,\,\nu_{s} \,;\, \ud{\sigma})$ be 
an expression of $\pi \in \BB(\lambda)$, and 
write each $\nu_{u} \in W\lambda$ for  
$1 \le u \le s$ as 
$\nu_{u}=
 \lambda-\beta_{u}+k_{u}d_{\lambda}\delta$, 
where $\beta_{u} \in \fin{Q}_{+}$ and $k_{u} \in \BZ$ 
(see Lemma~\ref{lem:orb}); 
note that $k_{1}=0$ by the assumption of the lemma. 
It follows from the definition of LS paths that 
$\nu_{1} \ge \nu_{2} \ge \cdots \ge \nu_{s}$, and hence 
from Remark~\ref{rem:bruhat} that 
$0=k_{1} \le k_{2} \le \cdots \le k_{s}$. 
Observe by \eqref{eq:path} that 
the coefficient $a_{0}^{-1}K$ of 
$\delta$ in the expression above of $\pi(1) \in P$ 
is equal to 
$\sum_{u=1}^{s}
(\sigma_{u}-\sigma_{u-1})k_{u}d_{\lambda}$. 
Therefore, we conclude that 
$a_{0}^{-1}K$, and hence $K$ is nonnegative. 
This proves the lemma.
\end{proof}

Let us denote by $\BB_{0}(\lambda) \subset \BB(\lambda)$ 
the connected component of the $P$-crystal $\BB(\lambda)$ 
containing the straight line $\pi_{\lambda}$. 
We know the following lemma from \cite[Lemma~4.2.3]{NSls}.
%
%%%%%%%%%%%%%%%%%%%
%%% lem:inverse %%%
%%%%%%%%%%%%%%%%%%%
%
\begin{lem} \label{lem:inverse}
Let $\eta \in \BB(\lambda)_{\cl}$. Then, the set 
$\cl^{-1}(\eta) \cap \BB_{0}(\lambda)$ is nonempty, 
where we set 
$\cl^{-1}(\eta):=\bigl\{\pi \in \BB(\lambda) \mid 
\cl(\pi)=\eta\bigr\}$.
Furthermore, if we take an arbitrary 
$\pi \in \cl^{-1}(\eta) \cap 
\BB_{0}(\lambda)$, then 
$\cl^{-1}(\eta) \cap \BB_{0}(\lambda) = 
 \bigl\{\pi+\pi_{kd_{\lambda}\delta} \mid k \in \BZ\bigr\}$.
\end{lem}
%
%
%%%%%%%%%%%%%%%%
%%% prop:deg %%%
%%%%%%%%%%%%%%%%
%
\begin{prop} \label{prop:deg}
Let $\eta \in \BB(\lambda)_{\cl}$. 
Then, the set $\cl^{-1}(\eta) \cap \BB_{0}(\lambda)$ contains 
a unique element $\pi_{\eta}$ such that 
$\iota(\pi_{\eta}) \in \lambda-\fin{Q}_{+}$.
\end{prop}

\begin{proof}
Let us take 
$\pi \in \cl^{-1}(\eta) \cap \BB_{0}(\lambda)$, and 
write its initial direction 
$\iota(\pi) \in W\lambda$ as 
$\iota(\pi)=\lambda-\beta+kd_{\lambda}\delta$ 
for $\beta \in \fin{Q}_{+}$ and $k \in \BZ$. 
Then we see from Lemma~\ref{lem:inverse} that 
$\pi_{\eta}:=\pi-\pi_{kd_{\lambda}\delta}$ is also
contained in $\cl^{-1}(\eta) \cap \BB_{0}(\lambda)$. 
In addition, 
it is easy to show (see Remark~\ref{rem:sum}) that 
$\iota(\pi_{\eta})$ is equal to 
$\iota(\pi)-kd_{\lambda}\delta=\lambda-\beta$, and 
hence that $\iota(\pi_{\eta}) \in \lambda-\fin{Q}_{+}$. 
This proves the existence of $\pi_{\eta}$. 
The uniqueness of $\pi_{\eta}$ follows immediately
from Lemma~\ref{lem:inverse}. 
This completes the proof of the proposition. 
\end{proof}

Let $\eta \in \BB(\lambda)_{\cl}$, and take 
$\pi_{\eta} \in \cl^{-1}(\eta) \cap \BB_{0}(\lambda)$ of 
Proposition~\ref{prop:deg}.
Then, by Lemma~\ref{lem:end}, we can 
write $\pi_{\eta}(1) \in P$ in the form 
$\pi_{\eta}(1)=
 \lambda-a_{0}^{-1}\beta+a_{0}^{-1} K \delta$ 
 with $\beta \in \fin{Q}_{+}$ and $K \in \BZ_{\ge 0}$. 
Now, we define 
the degree $\Deg_{\lambda}(\eta) \in \BZ_{\le 0}$ of 
the $\eta \in \BB(\lambda)_{\cl}$ by:
\begin{equation}
\Deg_{\lambda}(\eta)=-K.
\end{equation}
%
%%%%%%%%%%%%%%%%%%%%
%%% prop:deg-min %%%
%%%%%%%%%%%%%%%%%%%%
%
\begin{prop} \label{prop:deg-min}
Let $\eta \in \BB(\lambda)_{\cl}$. 
Take $\pi \in \cl^{-1}(\eta) \cap \BB(\lambda)$ such that 
$\iota(\pi) \in \lambda - \fin{Q}_{+}$ and $\pi \ne \pi_{\eta}$. 

\noindent 
{\rm (1)} \, If we write $\pi(1) \in P$ in the form 
$\pi(1)=\lambda-a_{0}^{-1}\beta'+a_{0}^{-1}K'\delta$ 
with $\beta' \in \fin{Q}_{+}$ and $K' \in \BZ_{\ge 0}$, 
then $-K' < \Deg_{\lambda}(\eta)$. 

\noindent 
{\rm (2)} \, 
If we write the final directions $\kappa(\pi_{\eta})$ and 
$\kappa(\pi)$ of $\pi_{\eta}$ and $\pi$ 
in the form 
$\kappa(\pi_{\eta})=\lambda-\beta+kd_{\lambda}\delta$ and 
$\kappa(\pi)=\lambda-\beta'+k'd_{\lambda}\delta$ with 
$\beta,\,\beta' \in \fin{Q}_{+}$ and $k,\,k' \in \BZ$, 
respectively, then $k < k'$. 
\end{prop}

\begin{rem}
Part (1) of Proposition~\ref{prop:deg-min} characterizes the degree 
$\Deg_{\lambda}(\eta) \in \BZ_{\le 0}$ of $\eta \in \BB(\lambda)_{\cl}$ 
as the maximum of the nonpositive integer $-K$ for which $\pi(1) \in P$ is 
of the form $\pi(1)=\lambda-a_{0}^{-1}\beta+a_{0}^{-1}K\delta$ with 
$\beta \in \fin{Q}_{+}$ and $K \in \BZ_{\ge 0}$, 
where $\pi \in \cl^{-1}(\eta) \cap \BB(\lambda)$ is such that 
$\iota(\pi) \in \lambda-\fin{Q}_{+}$. 
Furthermore, the maximum $\Deg_{\lambda}(\eta)$ is attained only by 
$\pi_{\eta} \in \cl^{-1}(\eta) \cap \BB_{0}(\lambda)$ of 
Proposition~\ref{prop:deg}.
\end{rem}

To prove this proposition, 
we need the following lemma, which can be proved 
by an argument in the proof of \cite[Theorem~3.1.1]{NSls}. 
%
%%%%%%%%%%%%%%%%
%%% lem:comp %%%
%%%%%%%%%%%%%%%%
%
\begin{lem} \label{lem:comp}
Each connected component of 
$\BB(\lambda)$ contains a unique element 
whose reduced expression is of the form\,{\rm:}
%
%%%%%%%%%%%%%%%
%%% eq:comp %%%
%%%%%%%%%%%%%%%
%
\begin{equation} \label{eq:comp}
(
 \lambda,\,
 \lambda+k_{2}d_{\lambda}\delta,\,\dots,\,
 \lambda+k_{s}d_{\lambda}\delta \ ; \ 
 \sigma_{0},\,
 \sigma_{1},\,\dots,\,
 \sigma_{s}
),
\end{equation}
with $k_{2},\,\dots,\,k_{s} \in \BZ$ and 
$0=\sigma_{0} < \sigma_{1} < \cdots < \sigma_{s}=1$. 
\end{lem}
%
%%%%%%%%%%%%%%%%
%%% rem:comp %%%
%%%%%%%%%%%%%%%%
%
\begin{rem} \label{rem:comp}
It follows from 
the definition of LS paths that $\lambda > 
 \lambda+k_{2}d_{\lambda}\delta > \cdots >
 \lambda+k_{s}d_{\lambda}\delta$.
Hence we see from Remark~\ref{rem:bruhat} that 
$0 < k_{2} < \cdots < k_{s}$. 
\end{rem}

\begin{proof}[Proof of Proposition~\ref{prop:deg-min}]
Assume that $\pi \in \BB(\lambda)$ 
lies in a connected component of $\BB(\lambda)$ 
containing an LS path $\pi'$ whose reduced expression is 
of the form \eqref{eq:comp}. 
We see from Proposition~\ref{prop:deg} and the assumption 
of the lemma that $\pi$ does not lie in $\BB_{0}(\lambda)$, 
and hence that $s \ge 2$. 

We set $\psi:=
(0,\,
 k_{2}d_{\lambda}\delta,\,\dots,\,
 k_{s}d_{\lambda}\delta \ ; \ 
 \sigma_{0},\,
 \sigma_{1},\,\dots,\,
 \sigma_{s})$; 
note that $\pi'=\pi_{\lambda}+\psi$. 
Let $X$ be a monomial of $X$ 
in the root operators $e_{j}$, $f_{j}$ for $j \in I$ 
such that $\pi=X \pi'$. Then we deduce that 
\begin{align*}
X \pi_{\lambda} 
 & = X (\pi'-\psi) = X\pi'-\psi
\quad \text{by \cite[Lemma~2.7.1]{NSls}} \\
& = \pi-\psi.
\end{align*}
Since $\cl(\pi)=\eta$, it follows that 
$\cl(X\pi_{\lambda})=\cl(\pi-\psi)=\cl(\pi)=\eta$. 
Hence we have $X\pi_{\lambda} \in 
\cl^{-1}(\eta) \cap \BB_{0}(\lambda)$. 
Also, because $\iota(\pi) \in \lambda-\fin{Q}_{+}$ 
and $\iota(\psi)=0$, we see that $\iota(X \pi_{\lambda})=
\iota(\pi-\psi)=\iota(\pi)-\iota(\psi) \in 
\lambda-\fin{Q}_{+}$. 
It follows from 
Proposition~\ref{prop:deg} that 
$X\pi_{\lambda}=\pi_{\eta}$. 
Thus we obtain $\pi_{\eta}=\pi-\psi$, 
and hence
\begin{equation*}
\Deg(\eta)=-K' + a_{0} \times 
  \text{(the coefficient of $\delta$ in $\psi(1)$)}, 
\quad \text{and} \quad
k=k'-k_{s}.
\end{equation*}
Since $s \ge 2$ as seen above, 
we deduce from Remark~\ref{rem:comp} 
(using \eqref{eq:path}) that 
the coefficient of $\delta$ in $\psi(1)$ is greater than $0$.
Therefore, we conclude that $-K' < \Deg(\eta)$. 
Also, since $k_{s} > 0$ with $s \ge 2$, 
it follows that $k < k'$. 
This completes the proof of the proposition.
\end{proof}
%
%
%
%==============================%
%     START SUBSECTION 0302    %
%==============================%
%
\subsection{Behavior of degree functions under the root operators.}
\label{subsec:deg-ro}
As in \ref{subsec:dfn-deg}, let $\lambda \in P^{0}_{+}$. 
%
%
%%%%%%%%%%%%%%%%%%
%%% lem:deg-ro %%%
%%%%%%%%%%%%%%%%%%
%
\begin{lem} \label{lem:deg-ro}
\noindent {\rm (1)} 
We have $\Deg_{\lambda}(\eta_{\cl(\lambda)})=0$. 

\noindent {\rm (2)} 
Let $\eta \in \BB(\lambda)_{\cl}$, and $j \in I$. 
If $e_{j}\eta \ne \bzero$, then 
{\rm(}see Remark~\ref{rem:initial}\,{\rm)}
%
%%%%%%%%%%%%%%%%%
%%% eq:deg-ro %%%
%%%%%%%%%%%%%%%%%
%
\begin{equation} \label{eq:deg-ro}
\Deg_{\lambda}(e_{j}\eta)=
 \begin{cases}
 \Deg_{\lambda}(\eta)-1 
  & \text{\rm if $j=0$ and 
    $\iota(e_{0}\eta)=\iota(\eta)$}, \\[1.5mm]
 \Deg_{\lambda}(\eta)-(\iota(\eta))(h_{0})-1 
  & \text{\rm if $j=0$ and 
    $\iota(e_{0}\eta)=r_{0}(\iota(\eta))$}, \\[1.5mm]
 \Deg_{\lambda}(\eta) & \text{\rm if $j \ne 0$}.
 \end{cases}
\end{equation}
\end{lem}

\begin{proof}
Part (1) is obvious from the definition of $\Deg_{\lambda}$, 
since $\pi_{\eta_{\cl(\lambda)}}=\pi_{\lambda}$. 
Let us prove part~(2).
It is obvious that 
$e_{j}\pi_{\eta} \in \BB_{0}(\lambda)$ 
since $\pi_{\eta} \in \BB_{0}(\lambda)$ by definition. 
Also, we know from \eqref{eq:cl-ro} that 
$\cl(e_{j}\pi_{\eta})=e_{j}\cl(\pi_{\eta})=e_{j}\eta$. 
Let us write $\pi_{\eta}(1) \in P$ in the form 
$\pi_{\eta}(1)=
 \lambda-a_{0}^{-1} \beta - 
 a_{0}^{-1} \Deg_{\lambda}(\eta)\delta$
with $\beta \in \fin{Q}_{+}$. 

First, assume that $j \ne 0$. 
We deduce from Remark~\ref{rem:initial} along with Lemma~\ref{lem:orb} 
that $\iota(e_{j}\pi_{\eta}) \in \lambda-\fin{Q}_{+}$. 
Because 
$e_{j}\pi_{\eta} \in \BB_{0}(\lambda)$ and 
$\cl(e_{j}\pi_{\eta})=e_{j}\eta$, 
it follows from Proposition~\ref{prop:deg} that 
$\pi_{e_{j}\eta}=e_{j}\pi_{\eta}$. 
Since $j \ne 0$, we have
\begin{equation*}
\pi_{e_{j}\eta}(1)=(e_{j}\pi_{\eta})(1)=
\pi_{\eta}(1)+\alpha_{j}=
 \lambda-
 \underbrace{(a_{0}^{-1} \beta - \alpha_{j})}_{\in a_{0}^{-1}\fin{Q}_{+}} - 
 a_{0}^{-1} \Deg_{\lambda}(\eta)\delta, 
\end{equation*}
and hence 
$\Deg_{\lambda}(e_{j}\eta)=\Deg_{\lambda}(\eta)$. 

Next, assume that $j=0$ and 
$\iota(e_{0}\eta)=\iota(\eta)$. 
Then we deduce (using \eqref{eq:H}) 
from the definitions of the root operator
$e_{0}$ for $\BB(\lambda)$ and 
the one for $\BB(\lambda)_{\cl}$ that 
$\iota(e_{0}\pi_{\eta})=\iota(\pi_{\eta})$, and hence 
$\iota(e_{0}\pi_{\eta}) \in \lambda-\fin{Q}_{+}$. 
Because $e_{0}\pi_{\eta} \in \BB_{0}(\lambda)$ and 
$\cl(e_{0}\pi_{\eta})=e_{0}\eta$, 
it follows from Proposition~\ref{prop:deg} that 
$\pi_{e_{0}\eta}=e_{0}\pi_{\eta}$. 
Now, define $\theta \in \fin{Q}_{+}$ by: 
$\theta=\delta-a_{0}\alpha_{0}$. 
Since $\alpha_{0}=a_{0}^{-1}(\delta-\theta)$, 
we have
\begin{align*}
\pi_{e_{0}\eta}(1) 
 & = (e_{0}\pi_{\eta})(1) 
   = \pi_{\eta}(1)+\alpha_{0} 
   = \lambda-a_{0}^{-1} \beta  - 
     a_{0}^{-1} \Deg_{\lambda}(\eta)\delta + 
     a_{0}^{-1}(\delta-\theta) \\
 & = \lambda-a_{0}^{-1}(\beta +\theta)-
     a_{0}^{-1} (\Deg_{\lambda}(\eta)-1) \delta,
\end{align*}
and hence $\Deg_{\lambda}(e_{0}\eta)=
\Deg_{\lambda}(\eta)-1$. 

Finally, assume that $j=0$ and 
$\iota(e_{0}\eta)=r_{0}(\iota(\eta))$. 
Then we deduce (using \eqref{eq:H}) 
from the definitions of the root operator
$e_{0}$ for $\BB(\lambda)$ and 
the one for $\BB(\lambda)_{\cl}$ that 
%
%%%%%%%%%%%%%%%%%%%
%%% eq:deg-ro01 %%%
%%%%%%%%%%%%%%%%%%%
%
\begin{align}
\iota(e_{0}\pi_{\eta}) & =r_{0}\iota(\pi_{\eta})=
\iota(\pi_{\eta})-(\iota(\pi_{\eta}))(h_{0})\alpha_{0} \nonumber \\
& =
\iota(\pi_{\eta})+a_{0}^{-1}(\iota(\pi_{\eta}))(h_{0})\theta
- a_{0}^{-1}(\iota(\pi_{\eta}))(h_{0})\delta, \label{eq:deg-ro01}
\end{align}
where $\theta=\delta-a_{0}\alpha_{0}$ as above.
Note that $\iota(\pi_{\eta})+a_{0}^{-1}(\iota(\pi_{\eta}))(h_{0})\theta \in 
\lambda-\sum_{j \in I_{0}}\BQ \alpha_{j}$, since 
$\iota(\pi_{\eta}) \in \lambda-\fin{Q}_{+}$. 
Because $\iota(e_{0}\pi_{\eta}) \in W\lambda$, 
it follows from Lemma~\ref{lem:orb} and 
the linear independence of $\alpha_{j}$, $j \in I$, and $\delta$ that 
$\iota(\pi_{\eta})+a_{0}^{-1}(\iota(\pi_{\eta}))(h_{0})\theta \in 
\lambda - \fin{Q}_{+}$ and $a_{0}^{-1}(\iota(\pi_{\eta}))(h_{0}) \in 
 \BZ d_{\lambda}$.
Hence we have 
$a_{0}^{-1}(\iota(\pi_{\eta}))(h_{0})=kd_{\lambda}$
for some $k \in \BZ$.
Because $e_{0}\pi_{\eta} \in \BB_{0}(\lambda)$ and 
$\cl(e_{0}\pi_{\eta})=e_{0}\eta \in \BB_{0}(\lambda)$, 
we deduce from Lemma~\ref{lem:inverse} that 
$e_{0}\pi_{\eta}+\pi_{kd_{\lambda}\delta} \in 
 \cl^{-1}(e_{0}\eta) \cap \BB_{0}(\lambda)$.
In addition,
\begin{align*}
\iota(e_{0}\pi_{\eta}+\pi_{kd_{\lambda}\delta}) 
& = \iota(e_{0}\pi_{\eta})+\iota(\pi_{kd_{\lambda}\delta}) 
\quad \text{by Remark~\ref{rem:sum}} \\
& = \iota(\pi_{\eta})+
  a_{0}^{-1}(\iota(\pi_{\eta}))(h_{0})\theta-
  kd_{\lambda}\delta+ kd_{\lambda}\delta
\quad \text{by \eqref{eq:deg-ro01}} \\
& = \iota(\pi_{\eta})+
  a_{0}^{-1}(\iota(\pi_{\eta}))(h_{0})\theta, 
\end{align*}
which lies in $\lambda-\fin{Q}_{+}$ as seen above. 
Therefore, by Proposition~\ref{prop:deg}, 
we have $\pi_{e_{0}\eta}=
 e_{0}\pi_{\eta}+\pi_{kd_{\lambda}\delta}$. 
From this, we obtain 
\begin{align*}
\pi_{e_{0}\eta}(1) 
 & = (e_{0}\pi_{\eta})(1)+\pi_{kd_{\lambda}\delta}(1)
   = \pi_{\eta}(1)+\alpha_{0}+kd_{\lambda}\delta \\
 & = \lambda-a_{0}^{-1} \beta - 
     a_{0}^{-1} \Deg_{\lambda}(\eta)\delta + 
     a_{0}^{-1}(\delta-\theta)+ kd_{\lambda}\delta \\
 & = \lambda-a_{0}^{-1} (\beta+\theta) - 
     a_{0}^{-1} (\Deg_{\lambda}(\eta)-a_{0}kd_{\lambda}-1) \delta,
\end{align*}
and hence 
\begin{equation*}
\Deg_{\lambda}(e_{0}\eta)=
\Deg_{\lambda}(\eta)-a_{0}kd_{\lambda}-1 = 
\Deg_{\lambda}(\eta)-(\iota(\pi_{\eta}))(h_{0})-1.
\end{equation*}
Note that 
$(\iota(\pi_{\eta}))(h_{0})=(\iota(\eta))(h_{0})$ 
since $\cl(\iota(\pi_{\eta}))=\iota(\cl(\pi_{\eta}))=\iota(\eta)$ 
by Remark~\ref{rem:cl-expr}.
Thus we conclude that 
$\Deg_{\lambda}(e_{0}\eta)=
 \Deg_{\lambda}(\eta)-(\iota(\eta))(h_{0})-1$, as desired.
This completes the proof of the lemma.
\end{proof}
%
%%%%%%%%%%%%%%%%%%%
%%% lem:deg-max %%%
%%%%%%%%%%%%%%%%%%%
%
\begin{lem} \label{lem:deg-max}
Let $\eta \in \BB(\lambda)_{\cl}$, and $j \in I$.
Assume that $e_{j}\eta \ne \bzero$, and that 
$(\iota(\eta))(h_{j}) \le 0$. Then, 
%
%%%%%%%%%%%%%%%%%%
%%% eq:deg-max %%%
%%%%%%%%%%%%%%%%%%
%
\begin{equation} \label{eq:deg-max}
\Deg_{\lambda}(e_{j}^{\max}\eta)=
\begin{cases}
 \Deg_{\lambda}(\eta)-\ve_{0}(\eta)-(\iota(\eta))(h_{0}) 
 & \text{\rm if } j=0, \\[2mm]
 \Deg_{\lambda}(\eta)
 & \text{\rm if } j \ne 0.
\end{cases}
\end{equation}
\end{lem}

\begin{proof}
If $j \ne 0$, then it follows immediately from 
Lemma~\ref{lem:deg-ro}\,(2) that 
$\Deg_{\lambda}(e_{j}^{\max}\eta)=\Deg_{\lambda}(\eta)$. 
Now assume that $j=0$. If $\ve_{0}(\eta)=0$, i.e., 
$e_{0}\eta=\bzero$, then 
we see from Lemma~\ref{lem:init1}\,(1) that 
$(\iota(\eta))(h_{0}) \ge 0$, 
which, when combined with the assumption of the lemma,
implies that $(\iota(\eta))(h_{0})=0$. 
Hence we have
\begin{equation*}
\Deg_{\lambda}(e_{0}^{\max}\eta)=
\Deg_{\lambda}(e_{0}^{0}\eta)=
\Deg_{\lambda}(\eta) = 
\Deg_{\lambda}(\eta)-
\underbrace{\ve_{0}(\eta)}_{= 0}-
\underbrace{(\iota(\eta))(h_{0})}_{= 0}.
\end{equation*}
It remains to consider the case 
$\ve_{0}(\eta) \ge 1$. 
From Lemma~\ref{lem:init1}\,(2) and 
Lemma~\ref{lem:deg-ro}\,(2), it follows that 
%
%%%%%%%%%%%%%%%%%%%
%%% eq:deg-max1 %%%
%%%%%%%%%%%%%%%%%%%
%
\begin{equation} \label{eq:deg-max1}
\Deg_{\lambda}(e_{0}^{\ve_{0}(\eta)-1}\eta)=
\Deg_{\lambda}(\eta)-\ve_{0}(\eta)+1.
\end{equation}
Therefore, 
by using Lemma~\ref{lem:init1}\,(2), (3), 
we deduce from Lemma~\ref{lem:deg-ro}\,(2) that 
\begin{align*}
\Deg_{\lambda}(e_{0}^{\max}\eta) 
& = 
  \Deg_{\lambda}(e_{0}^{\ve_{0}(\eta)}\eta)=
  \Deg_{\lambda}(e_{0}e_{0}^{\ve_{0}(\eta)-1}\eta) \\ 
& =
  \Deg_{\lambda}(e_{0}^{\ve_{0}(\eta)-1}\eta) -
  (\iota(\eta))(h_{0})-1 \\
& = 
  \Deg_{\lambda}(\eta)-\ve_{0}(\eta)+1-
  (\iota(\eta))(h_{0})-1 \qquad \text{by \eqref{eq:deg-max1}} \\
& = 
  \Deg_{\lambda}(\eta)-\ve_{0}(\eta)-(\iota(\eta))(h_{0}).
\end{align*}
This proves the lemma.
\end{proof}

%=========================%
%     START SECTION 04    %
%=========================%
%
\section{Relation between 
 energy functions and degree functions.}
\label{sec:main}

%==============================%
%     START SUBSECTION 0401    %
%==============================%
%
\subsection{Main results.}
\label{subsec:main}

Let $\bi=(i_{1},\,i_{2},\,\dots,\,i_{n})$ be 
an arbitrary sequence of elements of $I_{0}$, and 
define the tensor product $P_{\cl}$-crystal 
$\BB_{\bi}=
 \BB(\vpi_{i_{1}})_{\cl} \otimes 
 \BB(\vpi_{i_{2}})_{\cl} \otimes \cdots \otimes 
 \BB(\vpi_{i_{n}})_{\cl}$.
For an element 
$\eta_{1} \otimes \eta_{2} \otimes \cdots \otimes 
 \eta_{n} \in \BB_{\bi}$, 
we define $\eta_{l}^{(k)} \in \BB(\vpi_{i_{l}})_{\cl}$, 
$1 \le k < l \le n$, 
as follows (see \cite[\S3]{HKOTY} and \cite[\S3.3]{HKOTT}). 
There exists a unique isomorphism
\begin{align*}
&
\BB(\vpi_{i_{k}})_{\cl} \otimes 
\BB(\vpi_{i_{k+1}})_{\cl} \otimes \cdots \otimes 
\BB(\vpi_{i_{l-1}})_{\cl} \otimes \BB(\vpi_{i_{l}})_{\cl} \\
& \hspace*{20mm} \stackrel{\sim}{\rightarrow} 
\BB(\vpi_{i_{l}})_{\cl} \otimes 
\BB(\vpi_{i_{k}})_{\cl} \otimes \cdots \otimes 
\BB(\vpi_{i_{l-2}})_{\cl} \otimes \BB(\vpi_{i_{l-1}})_{\cl}
\end{align*}
of $P_{\cl}$-crystals, 
which is given as the composition 
$R_{\vpi_{i_{k}},\vpi_{i_{l}}} \circ 
 R_{\vpi_{i_{k+1}},\vpi_{i_{l}}} \circ \cdots \circ 
 R_{\vpi_{i_{l-1}},\vpi_{i_{l}}}$ 
of combinatorial $R$-matrices (see \S\ref{subsec:tpd}); 
for uniqueness, see Lemma~\ref{lem:simple}\,(3).
We define $\eta_{l}^{(k)}$ to be 
the first factor
(which lies in $\BB(\vpi_{i_{l}})_{\cl}$) of 
the image of 
$\eta_{k} \otimes \eta_{k+1} \otimes 
\cdots \otimes \eta_{l} \in 
\BB(\vpi_{i_{k}})_{\cl} \otimes 
\BB(\vpi_{i_{k+1}})_{\cl} \otimes \cdots \otimes 
\BB(\vpi_{i_{l}})_{\cl}$ 
under the above isomorphism of $P_{\cl}$-crystals.
For convenience, we set $\eta_{l}^{(l)}:=\eta_{l}$ 
for $1 \le l \le n$. 

For each $1 \le k \le n$, 
take (and fix) an arbitrary element 
$\eta_{k}^{\flat} \in \BB(\vpi_{i_{k}})_{\cl}$ 
such that $f_{j}\eta_{k}^{\flat}=\bzero$ 
for all $j \in I_{0}$.
Note that such an element 
$\eta^{\flat}_{k} \in \BB(\vpi_{i_{k}})_{\cl}$ 
actually exists. 
Indeed, for each $i \in I_{0}$, we know 
from Remark~\ref{rem:str-cl} 
that $\eta_{\ti{\vpi}_{i}} \in \BB(\vpi_{i})_{\cl}$, 
where $\ti{\vpi}_{i}:=w_{0}\cl(\vpi_{i}) \in P_{\cl}$ 
(see also Remark~\ref{rem:orbcl}). 
It follows immediately from the definition of the root operators 
$f_{j}$, $j \in I_{0}$, that $f_{j}\eta_{\ti{\vpi}_{i}}=\bzero$ 
for all $j \in I_{0}$.

Now, following \cite[\S3]{HKOTY} and \cite[\S3.3]{HKOTT}, 
we define the energy function $D_{\bi}:
 \BB_{\bi}=
 \BB(\vpi_{i_{1}})_{\cl} \otimes 
 \BB(\vpi_{i_{2}})_{\cl} \otimes \cdots \otimes 
 \BB(\vpi_{i_{n}})_{\cl} \rightarrow \BZ$ by: 
%
%%%%%%%%%%%%%
%%% eq:Di %%%
%%%%%%%%%%%%%
%
\begin{align}
& 
D_{\bi}(
  \eta_{1} \otimes 
  \eta_{2} \otimes \cdots \otimes 
  \eta_{n}) = \nonumber \\
& \hspace*{15mm}
\sum_{1 \le k < l \le n} 
 H_{\vpi_{i_{k}}, \vpi_{i_{l}}}
 (\eta_{k} \otimes \eta_{l}^{(k+1)})
 + \sum_{k=1}^{n} 
 H_{\vpi_{i_{k}}, \vpi_{i_{k}}}
  (\eta^{\flat}_{k} \otimes \eta_{k}^{(1)}). \label{eq:Di}
\end{align}
Also, we define a constant $D_{\bi}^{\ext} \in \BZ$ by:
%
%%%%%%%%%%%%%%%%
%%% eq:Diext %%%
%%%%%%%%%%%%%%%%
%
\begin{equation}  \label{eq:Diext}
D_{\bi}^{\ext} = 
\sum_{k=1}^{n} 
 H_{\vpi_{i_{k}}, \vpi_{i_{k}}}
  (\eta^{\flat}_{k} \otimes \eta_{\cl(\vpi_{i_{k}})}). 
\end{equation}

The main result of this paper is the following theorem. 
%
%%%%%%%%%%%%%%%%
%%% thm:main %%%
%%%%%%%%%%%%%%%%
%
\begin{thm} \label{thm:main}
Let $\bi=(i_{1},\,i_{2},\,\dots,\,i_{n})$ be 
an arbitrary sequence of elements of $I_{0}$, and set 
$\lambda:=\sum_{k=1}^{n} \vpi_{i_{k}} \in P^{0}_{+}$. 
Then, for every $\eta \in \BB(\lambda)_{\cl}$, 
the following equation holds\,{\rm:}
%
%%%%%%%%%%%%%%%
%%% eq:main %%%
%%%%%%%%%%%%%%%
%
\begin{equation} \label{eq:main}
\Deg_{\lambda}(\eta)=D_{\bi}(\Psi_{\bi}(\eta))-D_{\bi}^{\ext}, 
\end{equation}
where $\Psi_{\bi}:\BB(\lambda)_{\cl} \rightarrow \BB_{\bi}$ 
is the isomorphism of $P_{\cl}$-crystals in Theorem~\ref{thm:tpd}.
\end{thm}

We will establish Theorem~\ref{thm:main} under the following plan. 
First, in \S\ref{subsec:lemmas}, we show some technical lemmas 
needed later. Next, in \S\ref{subsec:key}, using these lemmas, 
we prove Proposition~\ref{prop:key}, which is the key to our proof 
(in \S\ref{subsec:step1}) of Theorem~\ref{thm:step1} below. 
Finally, in \S\ref{subsec:step2}, we prove 
Theorem~\ref{thm:step2} below, which, when combined with 
Theorem~\ref{thm:step1}, establishes Theorem~\ref{thm:main}.
%
%
%%%%%%%%%%%%%%%%%
%%% thm:step1 %%%
%%%%%%%%%%%%%%%%%
%
\begin{thm} \label{thm:step1}
Let $\bi=(i_{1},\,i_{2},\,\dots,\,i_{n})$ be 
an arbitrary sequence of elements of $I_{0}$, and set 
$\lambda:=\sum_{k=1}^{n} \vpi_{i_{k}} \in P^{0}_{+}$.
Let $\eta \in \BB(\lambda)_{\cl}$, and set 
$\Psi_{\bi}(\eta):=
 \eta_{1} \otimes \eta_{2} \otimes \cdots \otimes \eta_{n} \in 
 \BB_{\bi}=
 \BB(\vpi_{i_{1}})_{\cl} \otimes \BB(\vpi_{i_{2}})_{\cl} 
 \otimes \cdots \otimes \BB(\vpi_{i_{n}})_{\cl}$.
Then, the following equation holds\,{\rm:}
%
%%%%%%%%%%%%%%%%
%%% eq:step1 %%%
%%%%%%%%%%%%%%%%
%
\begin{equation} \label{eq:step1}
\Deg_{\lambda}(\eta)=
\sum_{1 \le k < l \le n} 
 H_{\vpi_{i_{k}},\vpi_{i_{l}}}
 (\eta_{k} \otimes \eta_{l}^{(k+1)}) +
 \sum_{k=1}^{n} \Deg_{\vpi_{i_{k}}} (\eta_{k}^{(1)}).
\end{equation}
\end{thm}
%
%%%%%%%%%%%%%%%%%
%%% thm:step2 %%%
%%%%%%%%%%%%%%%%%
%
\begin{thm} \label{thm:step2}
Let $i \in I_{0}$, and 
let $\eta^{\flat} \in \BB(\vpi_{i})_{\cl}$ be an element of 
$\BB(\vpi_{i})_{\cl}$ such that $f_{j}\eta^{\flat} = \bzero$ 
for all $j \in I_{0}$. 
Then, for every $\eta \in \BB(\vpi_{i})_{\cl}$, 
the the following equation holds\,{\rm:}
%
%%%%%%%%%%%%%%%%
%%% eq:step2 %%%
%%%%%%%%%%%%%%%%
%
\begin{equation} \label{eq:step2}
\Deg_{\vpi_{i}}(\eta)= 
H_{\vpi_{i}, \vpi_{i}}(\eta^{\flat} \otimes \eta) - 
H_{\vpi_{i}, \vpi_{i}}(\eta^{\flat} \otimes \eta_{\cl(\vpi_{i})}).
\end{equation}
\end{thm}
%
%
%
%==============================%
%     START SUBSECTION 0402    %
%==============================%
%
\subsection{Some technical lemmas.}
\label{subsec:lemmas}
Recall from \cite[Proposition~6.3]{Kac} that 
a real root of $\Fg$ lies either 
in $a_{0}^{-1}\fin{Q}_{+}+a_{0}^{-1}\BZ\delta$ or 
in $-a_{0}^{-1} \fin{Q}_{+}+ a_{0}^{-1}\BZ\delta$.
%
%
%%%%%%%%%%%%%%
%%% lem:fp %%%
%%%%%%%%%%%%%%
%
\begin{lem} \label{lem:fp}
Let $\lambda \in P^{0}_{+}$ 
be a level-zero dominant integral weight, and 
let $w \in W$, $j \in I$. 

\noindent
{\rm (1)} 
If $\bigl(w(\cl(\lambda))\bigr)(h_{j}) < 0$, then 
the real root $w^{-1}(\alpha_{j})$ lies 
in $-a_{0}^{-1}\fin{Q}_{+}+a_{0}^{-1}\BZ\delta$. 

\noindent
{\rm (2)} 
If the real root $w^{-1}(\alpha_{j})$ lies in 
$-a_{0}^{-1}\fin{Q}_{+}+a_{0}^{-1}\BZ\delta$, 
then $\bigl(w(\cl(\lambda))\bigr)(h_{j}) \le 0$. 

\noindent
{\rm (3)} Assume that $\lambda$ is 
strictly level-zero dominant. Then, 
$\bigl(w(\cl(\lambda))\bigr)(h_{j}) < 0$ 
if and only if 
the real root $w^{-1}(\alpha_{j})$ lies in 
$-a_{0}^{-1}\fin{Q}_{+}+a_{0}^{-1}\BZ\delta$. 
\end{lem}
\begin{proof}
Since $(\alpha_{j},\alpha_{j}) \in \BZ_{> 0}$ for all $j \in I$, and 
since 
\begin{equation*}
\bigl(w(\cl(\lambda))\bigr)(h_{j}) = 
\bigl(w(\lambda)\bigr)(h_{j})=
\frac{2(w(\lambda),\alpha_{j})}{(\alpha_{j},\alpha_{j})}=
\frac{2(\lambda,w^{-1}(\alpha_{j}))}{(\alpha_{j},\alpha_{j})},
\end{equation*}
it follows immediately that 
$\bigl(w(\cl(\lambda))\bigr)(h_{j}) < 0$ if and only if
$(\lambda,w^{-1}(\alpha_{j})) < 0$, and that 
$\bigl(w(\cl(\lambda))\bigr)(h_{j})=0$ if and only if
$(\lambda,w^{-1}(\alpha_{j}))=0$. 
In addition, because $\lambda$ is level-zero dominant and 
$(\alpha_{j},\alpha_{j}) \in \BZ_{> 0}$ for all $j \in I_{0}$, 
we have
\begin{equation*}
\bigl(
 \lambda,\, \pm a_{0}^{-1}\fin{Q}_{+}+a_{0}^{-1}\BZ\delta
\bigr)
 =
\bigl(
 \lambda,\, \pm a_{0}^{-1}\fin{Q}_{+}
\bigr)
\subset \pm \BQ_{\ge 0}.
\end{equation*}
All the assertions of the lemma follows immediately 
from the discussion above. 
\end{proof}
%
%
%%%%%%%%%%%%%%
%%% lem:s1 %%%
%%%%%%%%%%%%%%
%
\begin{lem} \label{lem:s1}
Let $\lambda \in P^{0}_{+}$ 
be a level-zero dominant integral weight, and 
let $\eta \in \BB(\lambda)_{\cl}$, $j \in I$. 
Assume that 
$\eta$ has an expression of the form 
$\eta=(\mu_{1},\,\mu_{2} \,;\, 0,\,\sigma,\,1)$, with 
$\mu_{1},\,\mu_{2} \in \cl(W\lambda)=\fin{W}\cl(\lambda)$ 
and $0 < \sigma < 1$.
If $\mu_{1}(h_{j}) < 0$, then 
$e_{j}^{\max}\eta=
 (r_{j}(\mu_{1}),\,\mu_{2}' \,;\, 0,\,\sigma,\,1)$, where 
\begin{equation}
\mu_{2}':=
\begin{cases}
\mu_{2} & \text{\rm if } \mu_{2}(h_{j}) \ge 0, \\[1.5mm]
r_{j}\mu_{2} & \text{\rm if } \mu_{2}(h_{j}) < 0.
\end{cases}
\end{equation}
\end{lem}

\begin{proof}
First, assume that $\mu_{2}(h_{j}) \ge 0$. 
Then, since $\mu_{1}(h_{j}) < 0$ by the assumption of the lemma, 
it follows that the function 
$H^{\eta}_{j}(t)$ is strictly decreasing 
on the interval $[0,\sigma]$, and 
$m^{\eta}_{j}=H^{\eta}_{j}(\sigma) < 0$; 
note that $\ve_{j}(\eta)=-m^{\eta}_{j}$ by \eqref{eq:ve-cl}.
For $0 \le l \le \ve_{j}(\eta)=-m^{\eta}_{j}$, 
let $\sigma^{(l)}$ be the unique point in $[0,\sigma]$ 
such that $H^{\eta}_{j}(\sigma^{(l)})=m^{\eta}_{j}+l$; 
observe that $0=\sigma^{(\ve_{j}(\eta))} < 
\sigma^{(\ve_{j}(\eta)-1)} < \cdots < \sigma^{(0)}=\sigma$.
Now it is easily shown by induction on $l$ that 
for $0 \le l \le \ve_{j}(\eta)$, 
\begin{equation*}
(e_{j}^{l}\eta)(t) = 
\begin{cases}
\eta(t) & 
  \text{if } 0 \le t \le \sigma^{(l)}, \\[2mm]
\eta(\sigma^{(l)})+r_{j}\bigl(\eta(t)-\eta(\sigma^{(l)})\bigr) & 
  \text{if } \sigma^{(l)} \le t \le \sigma, \\[2mm]
\eta(t)+l\alpha_{j} & 
  \text{if } \sigma \le t \le 1.
\end{cases}
\end{equation*}
In particular, by taking $l=\ve_{j}(\eta)$, we have
\begin{equation*}
(e_{j}^{\max}\eta)(t) = 
\begin{cases}
r_{j}\bigl(\eta(t)\bigr) & 
  \text{if } 0 \le t \le \sigma, \\[2mm]
\eta(t)+\ve_{j}(\eta)\alpha_{j} & 
  \text{if } \sigma \le t \le 1, 
\end{cases}
\end{equation*}
which implies that $e_{j}^{\max}\eta=
 (r_{j}(\mu_{1}),\,\mu_{2} \,;\, 0,\,\sigma,\,1)$, 
since $\eta=(\mu_{1},\,\mu_{2} \,;\, 0,\,\sigma,\,1)$. 

Next, assume that $\mu_{2}(h_{j}) < 0$. 
Then, since $\mu_{1}(h_{j}) < 0$ by the assumption of the lemma, 
it follows that the function 
$H^{\eta}_{j}(t)$ is strictly decreasing on the interval
$[0,1]$, and $m^{\eta}_{j}=H^{\eta}_{j}(1) < 0$; 
note that $\ve_{j}(\eta)=-m^{\eta}_{j}$ by \eqref{eq:ve-cl}.
For $0 \le l \le \ve_{j}(\eta)=-m^{\eta}_{j}$, 
let $\sigma^{(l)}$ be the unique point in $[0,1]$ 
such that $H^{\eta}_{j}(\sigma^{(l)})=m^{\eta}_{j}+l$; 
observe that $0=\sigma^{(\ve_{j}(\eta))} < 
\sigma^{(\ve_{j}(\eta)-1)} < \cdots < \sigma^{(0)}=1$.
Now it is easily shown by induction on $l$ that 
for $0 \le l \le \ve_{j}(\eta)$, 
\begin{equation*}
(e_{j}^{l}\eta)(t) = 
\begin{cases}
\eta(t) & 
  \text{if } 0 \le t \le \sigma^{(l)}, \\[2mm]
\eta(\sigma^{(l)})+r_{j}\bigl(\eta(t)-\eta(\sigma^{(l)})\bigr) & 
  \text{if } \sigma^{(l)} \le t \le 1.
\end{cases}
\end{equation*}
In particular, by taking $l=\ve_{j}(\eta)$, we have
\begin{equation*}
(e_{j}^{\max}\eta)(t) = 
r_{j}\bigl(\eta(t)\bigr) 
\qquad 
\text{for $t \in [0,1]$}, 
\end{equation*}
which implies that 
$e_{j}^{\max}\eta=
 (r_{j}(\mu_{1}),\,r_{j}(\mu_{2}) \,;\, 0,\,\sigma,\,1)$, 
since $\eta=(\mu_{1},\,\mu_{2} \,;\, 0,\,\sigma,\,1)$. 
This completes the proof of the lemma.
\end{proof}
%
%
%%%%%%%%%%%%%%
%%% lem:s2 %%%
%%%%%%%%%%%%%%
%
\begin{lem} \label{lem:s2}
Let $\lambda \in P^{0}_{+}$. 
Let 
$\eta=
 (\mu_{1},\,\mu_{2},\,\dots,\,\mu_{s} \, ; \, 
  \sigma_{0},\,\sigma_{1},\,\sigma_{2},\,\dots,\,\sigma_{s})$ 
be an expression of $\eta \in \BB(\lambda)_{\cl}$, and 
assume that $s \ge 2$. 

\noindent {\rm (1)} 
$\eta':=
 (\mu_{1},\,\mu_{2} \,;\, 
  \sigma_{0},\,\sigma_{1},\,\sigma_{s})$ 
is contained in $\BB(\lambda)_{\cl}$. 

\noindent {\rm (2)} 
If $j \in I$ satisfies $\mu_{1}(h_{j}) < 0$, 
then 
$(e_{j}^{\max}\eta)(t)=(e_{j}^{\max}\eta')(t)$
for all $t \in [0,\sigma_{2}]$. Hence, by Lemma~\ref{lem:s1}, 
$e_{j}^{\max}\eta \in \BB(\lambda)_{\cl}$ has an expression of the form\,{\rm:}
\begin{equation}
e_{j}^{\max}\eta=
\bigl(r_{j}\mu_{1},\,\mu_{2}',\,\mu_{3}'\dots,\,\mu_{s'}' \, ; \, 
  \sigma_{0},\,\sigma_{1},\,\sigma_{2},\,
  \sigma_{3}',\,\dots,\,\sigma_{s'}'\bigr),
\end{equation}
where 
\begin{equation}
\mu_{2}'=
\begin{cases}
\mu_{2} & \text{\rm if } \mu_{2}(h_{j}) \ge 0, \\[1.5mm]
r_{j}\mu_{2} & \text{\rm if } \mu_{2}(h_{j}) < 0.
\end{cases}
\end{equation}
\end{lem}

%%%%
\vsp
%%%%

\begin{proof}
(1) Let $\pi \in \BB(\lambda)$ be such that $\cl(\pi)=\eta$, 
and let 
\begin{equation*}
\pi=(\nu_{1},\,\nu_{2},\,\dots,\,\nu_{s''} \, ; \, 
  \sigma_{0}'',\,\sigma_{1}'',\,
  \sigma_{2}'',\,\dots,\,
  \sigma_{s''}'')
\end{equation*}
be an expression of $\pi$. 
By ``inserting'' (see \cite[Remark~2.5.2\,(2)]{NSls}) 
$\sigma_{1}$ (resp., $\sigma_{2}$) between 
$\sigma_{k}''$ and $\sigma_{k+1}''$ such that 
$\sigma_{k}'' < \sigma_{1} \, 
\text{(resp., $\sigma_{2}$)} < \sigma_{k+1}''$ 
if necessary, 
we may assume that 
there exists $1 \le u_{1} < u_{2} \le s''$
such that $\sigma_{u_{1}}''=\sigma_{1}$ and 
$\sigma_{u_{2}}''=\sigma_{2}$.
By Remark~\ref{rem:cl-expr} and 
the condition that $\cl(\pi)=\eta$, we have 
%
%%%%%%%%%%%%%%%
%%% eq:s2-a %%%
%%%%%%%%%%%%%%%
%
\begin{equation} \label{eq:s2-a}
\begin{array}{l}
\cl(\nu_{u})=\mu_{1} \quad \text{for all $1 \le u \le u_{1}$}, \\[2mm]
\cl(\nu_{u})=\mu_{2} \quad \text{for all $u_{1}+1 \le u \le u_{2}$}.
\end{array}
\end{equation} 
Set $\pi':=
(\nu_{1},\,\nu_{2},\,\dots,\,
 \nu_{u_{2}} \ ; \ 
 \sigma_{0}'',\,\sigma_{1}'',\,\dots,\,
 \sigma_{u_{2}-1}'',\,\sigma_{s''}'')$. 
Then we can easily deduce from the definition of LS paths 
(see also \cite[Lemma~4.5\,b)]{L2}) that $\pi' \in \BB(\lambda)$. 
Furthermore, it is clear from Remark~\ref{rem:cl-expr} and 
\eqref{eq:s2-a} that $\cl(\pi')=\eta'$. 
Thus we have proved that $\eta' \in \BB(\lambda)_{\cl}$. 

\vsp

\noindent 
(2) 
For $0 \le l \le \ve_{j}(\eta)-1$, 
we set $\eta_{l}:=e_{j}^{l}\eta \in \BB(\lambda)_{\cl}$, and 
\begin{equation*}
t_{1}^{(l)}:=\min\bigl\{t \in [0,1] \mid 
     H^{\eta_{l}}_{j}(t)=m^{\eta_{l}}_{j} \bigr\}, 
\quad
t_{0}^{(l)}:=\max\bigl\{t \in [0,t_{1}^{l}] \mid
     H^{\eta_{l}}_{j}(t) = m^{\eta_{l}}_{j}+1\bigr\}.
\end{equation*}
We note that by \eqref{eq:ve-cl}, 
%
%%%%%%%%%%%%%%%
%%% eq:ve01 %%%
%%%%%%%%%%%%%%%
%
\begin{equation} \label{eq:ve01}
m^{\eta_{l}}_{j}=
-\ve_{j}(\eta_{l})=
-\ve_{j}(e_{j}^{l}\eta)=
-\ve_{j}(\eta)+l=m^{\eta}_{j}+l
\quad
\text{for $0 \le l \le \ve_{j}(\eta)-1$}.
\end{equation}

Now, let us assume that $\mu_{2}(h_{j}) \ge 0$. 
Since $\mu_{1}(h_{j}) < 0$ by the assumption of the lemma, 
it follows that $k_{1}:=
 H^{\eta}_{j}(\sigma_{1}) < 0$. 
In addition, since $\mu_{2}(h_{j}) \ge 0$, 
we see that 
the function $H^{\eta}_{j}(t)$, $t \in [0,1]$, 
attains a local minimum at $t=\sigma_{1}$. 
Therefore, we obtain 
$k_{1} \in \BZ_{< 0}$ by Remark~\ref{rem:chain}. 
Observe that 
%
%%%%%%%%%%%%%%%
%%% eq:gek1 %%%
%%%%%%%%%%%%%%%
%
\begin{equation} \label{eq:gek1}
H^{\eta}_{j}(t) \ge k_{1}
\quad
\text{for all }  t \in [0,\sigma_{2}].
\end{equation}
We set $l_{1}:=k_{1}-m^{\eta}_{j} \in \BZ_{\ge 0}$; 
note that $l_{1} < -m^{\eta}_{j}=\ve_{j}(\eta)$ 
since $k_{1} \in \BZ_{< 0}$.

%%%%%%%%%%%%%
\begin{claim}
For all $0 \le l \le l_{1}$, we have 
$\eta_{l}(t)=\eta(t)$ for all $t \in [0,\sigma_{2}]$. 
\end{claim}
%%%%%%%%%%%

\noindent {\it Proof of Claim.} 
We show the assertion by induction on $l$. 
When $l=0$, the assertion obviously holds. 
Assume that $0 < l \le l_{1}$ and that 
$\eta_{l-1}(t)=\eta(t)$ for all $t \in [0,\sigma_{2}]$. 
Note that $\eta_{l}=e_{j}\eta_{l-1}$. 
Therefore, by the definition of the root operator $e_{j}$, 
it suffices to show that $\sigma_{2} \le t_{0}^{(l-1)}$. 
We see from \eqref{eq:ve01} that 
%
%%%%%%%%%%%%%%%%%
%%% eq:gek1-1 %%%
%%%%%%%%%%%%%%%%%
%
\begin{equation} \label{eq:gek1-1}
H^{\eta_{l-1}}_{j}(t_{1}^{(l-1)}) = 
m^{\eta_{l-1}}_{j}=m^{\eta}_{j}+l-1 \le
m^{\eta}_{j}+l_{1}-1=k_{1}-1.
\end{equation}
Also, 
since $\eta_{l-1}(t)=\eta(t)$ for all $t \in [0,\sigma_{2}]$ 
by the inductive assumption, it follows from \eqref{eq:gek1} that
$H^{\eta_{l-1}}_{j}(t)=H^{\eta}_{j}(t) \ge k_{1}$ 
for all $t \in [0,\sigma_{2}]$. 
Hence we deduce from \eqref{eq:gek1-1} that 
$t_{1}^{(l-1)} \notin [0,\sigma_{2}]$, i.e., 
that $\sigma_{2} < t_{1}^{(l-1)}$.

We note that 
$H^{\eta_{l-1}}_{j}(\sigma_{2}) = 
 H^{\eta}_{j}(\sigma_{2}) \ge k_{1}=
 m^{\eta}_{j}+l_{1} \ge 
 m^{\eta}_{j}+l=m^{\eta_{l-1}}_{j}+1$
as seen above, and that 
$H^{\eta_{l-1}}_{j}(t_{1}^{(l-1)})=
 m^{\eta_{l-1}}_{j} < 
 m^{\eta_{l-1}}_{j}+1$. 
Therefore, from the continuity of the function 
$H^{\eta_{l-1}}_{j}(t)$ on the interval 
$\sigma_{2} \le t \le t_{1}^{(l-1)}$, we conclude that 
there exists $\sigma_{2} \le t' < t_{1}^{(l-1)}$ such that 
$H^{\eta_{l-1}}_{j}(t')=m^{\eta_{l-1}}_{j}+1$. 
It follows from the definition of $t_{0}^{(l-1)}$ that 
$t_{0}^{(l-1)} \ge t' \ge \sigma_{2}$. 
This proves the claim.

%%%%
\vsp
%%%%

From the claim above, by taking $l=l_{1}$, we obtain 
$\eta_{l_{1}}(t)=\eta(t)$ for all 
$t \in [0,\sigma_{2}]$. 
Consequently, we see 
from the definition of $t_{1}^{(l_{1})}$ 
that $t_{1}^{(l_{1})}=\sigma_{1}$, 
since $m^{\eta_{l_{1}}}_{j}=m^{\eta}_{j}+l_{1}=k_{1}$ and 
$H^{\eta}_{j}(\sigma_{1})=k_{1}$. 
Therefore, as in the proof of Lemma~\ref{lem:s1}, 
we can show (using $\eta(t)=\eta_{l_{1}}(t)$ for 
$t \in [0,\sigma_{2}]$) that 
\begin{align*}
(e_{j}^{\max}\eta)(t) & = 
(e_{j}^{-m^{\eta}_{j}}\eta)(t) = 
(e_{j}^{l_{1}-k_{1}}\eta)(t) \\[3mm]
& = 
(e_{j}^{-k_{1}}\eta_{l_{1}})(t)
=
\begin{cases}
r_{j}\bigl(\eta(t)\bigr) & 
  \text{if } 0 \le t \le \sigma_{1}, \\[2mm]
\eta(t)-k_{1}\alpha_{j} & 
  \text{if } \sigma_{1} \le t \le \sigma_{2}, \\[2mm]
\eta_{l_{1}}(t)-k_{1}\alpha_{j} & 
  \text{if } \sigma_{2} \le t \le 1.
\end{cases}
\end{align*}
From this and Lemma~\ref{lem:s1}, 
we conclude that 
$(e_{j}^{\max}\eta)(t)=(e_{j}^{\max}\eta')(t)$ 
for $t \in [0,\sigma_{2}]$.

The proof for the case $\mu_{2}(h_{j}) < 0$ 
is similar; we give only a sketch of the proof. 
Take the largest 
$u \in \bigl\{2,\,3,\,\dots,\,s\bigr\}$ such that 
$\mu_{u'}(h_{j}) < 0$ for all $1 \le u' \le u$. 
Then we see that the function 
$H^{\eta}_{j}(t)$, $t \in [0,1]$, attains 
a local minimum at $t=\sigma_{u}$, and hence that
$k_{u}:=H^{\eta}_{j}(\sigma_{u}) \in \BZ_{< 0}$ 
by Remark~\ref{rem:chain}. 
We set $l_{u}:=k_{u}-m^{\eta}_{j} \in \BZ_{\ge 0}$.
In exactly the same way as above, 
we can show that $\eta(t)=\eta_{l_{u}}(t)$ for all 
$t \in [0,\sigma_{u}]$.
Consequently, we have $t_{1}^{(l_{u})}=\sigma_{u}$, 
since $\mu_{u'}(h_{j}) < 0$ for all $1 \le u' \le u$. 
Therefore, 
as in the proof of Lemma~\ref{lem:s1}, 
we can show (using $\eta(t)=\eta_{l_{u}}(t)$ for 
$t \in [0,\sigma_{u}]$) that 
\begin{equation*}
(e_{j}^{\max}\eta)(t) = 
(e_{j}^{-k_{u}}\eta_{l_{u}})(t) = 
\begin{cases}
r_{j}\bigl(\eta(t)\bigr) & 
  \text{if } 0 \le t \le \sigma_{u}, \\[2mm]
\eta_{l_{u}}(t)-k_{u}\alpha_{j} & 
  \text{if } \sigma_{u} \le t \le 1. 
\end{cases}
\end{equation*}
From this and Lemma~\ref{lem:s1}, 
we conclude that 
$(e_{j}^{\max}\eta)(t)=(e_{j}^{\max}\eta')(t)$ 
for $t \in [0,\sigma_{2}]$.
This completes the proof of the lemma.
\end{proof}
%
%
%
%==============================%
%     START SUBSECTION 0403    %
%==============================%
%
\subsection{
  Key proposition to 
  the proof of Theorem~\ref{thm:step1}.
  }
\label{subsec:key}
%
%%%%%%%%%%%%%%%%
%%% prop:key %%%
%%%%%%%%%%%%%%%%
%
\begin{prop} \label{prop:key}
Let $\lambda \in P^{0}_{+}$ be 
a level-zero dominant integral weight, 
and let $\eta \in \BB(\lambda)_{\cl}$ be an element of 
$\BB(\lambda)_{\cl}$ such that 
$e_{j}\eta=\bzero$ for all $j \in I_{0}$. 
Then, there exists a sequence 
$j_{1},\,j_{2},\,\dots,\,j_{N} \in I$ such that 

\noindent
{\rm (A)} \, 
    $e_{j_{N}}^{\max}e_{j_{N-1}}^{\max} \cdots 
     e_{j_{1}}^{\max}\eta=\eta_{\cl(\lambda)}$, 

\noindent
{\rm (B)} \, 
    $r_{j_{1}}r_{j_{2}} \dots r_{j_{p}}(\alpha_{j_{p+1}}) \in 
     -a_{0}^{-1}\fin{Q}_{+} + a_{0}^{-1}\BZ \delta$ 
     for $p=0,\,1,\,\dots,\,N-1$.
\end{prop}

In order to prove 
Proposition~\ref{prop:key}, 
we need Lemmas~\ref{lem:strict1} 
and \ref{lem:strict2} below. 
%
%
%%%%%%%%%%%%%%%%%%%
%%% lem:strict1 %%%
%%%%%%%%%%%%%%%%%%%
%
\begin{lem} \label{lem:strict1}
Let $\lambda \in \sum_{i \in I_{0}}\BZ_{> 0} \vpi_{i}$ 
be a strictly level-zero dominant integral weight.
Let $\eta$ be an element of $\BB(\lambda)_{\cl}$ 
such that $\iota(\eta)=\cl(\lambda)$. 
Then, there exists a sequence 
$j_{1},\,j_{2},\,\dots,\,j_{N} \in I$ such that 

\noindent
{\rm (A)} \, 
  $e_{j_{N}}^{\max}e_{j_{N-1}}^{\max} \cdots 
   e_{j_{1}}^{\max}\eta=\eta_{\cl(\lambda)}$, 

\noindent
{\rm (B)'} \, 
  $\bigl(r_{j_{p}}r_{j_{p-1}} \cdots r_{j_{1}}(\cl(\lambda))\bigr)
   (h_{j_{p+1}}) < 0$ for $p=0,\,1,\,\dots,\,N-1$.
\end{lem}

\begin{proof}
First, we show the next sublemma, 
which is a special case of Lemma~\ref{lem:strict1}.
%
%%%%%%%%%%%%%%%%
%%% sublemma %%%
%%%%%%%%%%%%%%%%
%
\begin{sublem}
Keep the notation and assumption of 
Lemma~\ref{lem:strict1}. 
Further, we assume that 
$\eta \in \BB(\lambda)_{\cl}$ is of the form
$\eta=(\cl(\lambda),\,\mu \,;\, 0,\sigma,1)$, 
with $\mu \in \cl(W\lambda)=\fin{W}\cl(\lambda)$ and 
$0 < \sigma < 1$.
Then the assertion of Lemma~\ref{lem:strict1} holds. 
\end{sublem}

%%%%%%%%%%%%%%%%%%%%%%%%%%
\noindent {\it Proof of Sublemma.} 
Since $\BB(\lambda)_{\cl}$ is a finite set 
(see Remark~\ref{rem:regular}), it follows that 
$\max \bigl\{||\eta'|| \mid 
 \eta' \in \BB(\lambda)_{\cl}\bigr\} < \infty$, 
where $||\eta||$ is defined as in \eqref{eq:def-norm}; 
recall that $\wt(\eta)$ is defined by: 
$\wt(\eta)=\eta(1) \in P_{\cl}$.
Now we show the assertion of the sublemma
by descending induction on $||\eta||$.
When $||\eta||=\max \bigl\{||\eta'|| \mid 
 \eta' \in \BB(\lambda)_{\cl}\bigr\}$, 
it follows from Lemma~\ref{lem:norm} that 
$\eta$ is an extremal element. 
Since the initial direction $\iota(\eta)$ of $\eta$ 
is equal to $\cl(\lambda)$ by assumption, 
we deduce from Remark~\ref{rem:tpd02}
that $\eta=\eta_{\cl(\lambda)}$ 
(hence there is nothing to prove). 

Assume that 
$||\eta|| < \max \bigl\{||\eta'|| \mid 
 \eta' \in \BB(\lambda)_{\cl}\bigr\}$.
Set $\Lambda:=w_{0}\lambda$, 
where $w_{0} \in \fin{W}$ is 
the longest element of $\fin{W}$.
Then it is easy to check that 
$\Lambda(h_{j}) \in \BZ_{< 0}$ 
for all $j \in I_{0}$. 
Note that $\mu \in \cl(W\lambda)=\fin{W}\cl(\lambda)$ 
satisfies the condition that 
$\mu(h_{j}) \le 0$ for all $j \in I_{0}$ if and only if 
$\mu=\cl(\Lambda)$ (see Remark~\ref{rem:orbcl}). 

We see from [AK, Lemma~1.4] that there exists 
a sequence $j_{1},\,j_{2},\,\dots,\,j_{N'} \in I$
such that 

\vspace{1.5mm}

\noindent
(a) \, 
$r_{j_{N'}}r_{j_{N'-1}} \cdots r_{j_{1}}(\cl(\lambda))=\cl(\Lambda)$, 

\noindent
(b) \, 
$\bigl(r_{j_{p}}r_{j_{p-1}} \cdots r_{j_{1}}(\cl(\lambda))\bigr)
 (h_{j_{p+1}}) < 0$ for $p=0,\,1,\,\dots,\,N'-1$.

\vspace{1.5mm}

\noindent Set $\eta':=
 e_{j_{N'}}^{\max}e_{j_{N'-1}}^{\max} \cdots 
 e_{j_{1}}^{\max}\eta \in \BB(\lambda)_{\cl}$.
Then it follows from Lemma~\ref{lem:e-max} that 
%
%%%%%%%%%%%%%%
%%% eq:cs2 %%%
%%%%%%%%%%%%%%
%
\begin{equation} \label{eq:cs2}
||\eta'|| \ge ||\eta||.
\end{equation}
Furthermore, from (b), 
we deduce by repeated application of 
Lemma~\ref{lem:s1} that 
\begin{equation*}
e_{j_{N'}}^{\max}
e_{j_{N'-1}}^{\max} \cdots 
e_{j_{1}}^{\max}\eta 
= 
(r_{j_{N'}}r_{j_{N'-1}} \cdots r_{j_{1}}(\cl(\lambda)),\,
 \mu' \,;\, 0,\,\sigma,\,1) 
\qquad \text{for some $\mu' \in \cl(W\lambda)$},
\end{equation*}
and hence from (a) that 
\begin{align*}
\eta'
& =
e_{j_{N'}}^{\max}e_{j_{N'-1}}^{\max} \cdots 
e_{j_{1}}^{\max}\eta 
= 
(r_{j_{N'}}r_{j_{N'-1}} \cdots r_{j_{1}}(\cl(\lambda)),\,
 \mu' \,;\, 0,\,\sigma,\,1) \\
& =
(\cl(\Lambda),\,\mu' \, ; \, 0,\,\sigma,\,1).
\end{align*}

\paragraph{Case 1.} $\mu'= \cl(\Lambda)$. 

\noindent 
In this case, we have 
$\eta'= 
 (\cl(\Lambda),\,\cl(\Lambda) \, ; \, 0,\,\sigma,\,1)= 
 \eta_{\cl(\Lambda)}$. 
Let $w_{0}=r_{j_{N'+1}}r_{j_{N'+2}} \cdots r_{j_{N-1}}r_{j_{N}}$ be 
a reduced expression of $w_{0} \in \fin{W}$; note that 
$j_{N'+1},\,j_{N'+2},\,\dots,\,j_{N-1},\,j_{N} \in I_{0}$. 
Then, since $\Lambda=w_{0}\lambda$, it follows that 

\vspace{1.5mm}

\noindent
(c) \, 
$r_{j_{N}}r_{j_{N-1}} \cdots r_{j_{N'+1}}
 (\cl(\Lambda))=\cl(\lambda)$. 

\vspace{1.5mm}

\noindent 
In addition, using \cite[Lemma~3.11\,b)]{Kac}, 
we obtain

\vspace{1.5mm}

\noindent
(d) \, 
$\bigl(
   r_{j_{p}}r_{j_{p-1}} \cdots r_{j_{N'+1}}(\cl(\Lambda))
 \bigr)
 (h_{j_{p+1}}) < 0$ 
 for all $p=N',\,N'+1,\,\dots,\,N-1$.

\vspace{1.5mm}

\noindent 
From (c) and (d), 
we deduce by repeated use of 
Lemma~\ref{lem:s1} that 
\begin{align*}
& 
e_{j_{N}}^{\max}e_{j_{N-1}}^{\max} \cdots 
e_{j_{N'+1}}^{\max}
e_{j_{N'}}^{\max}e_{j_{N'-1}}^{\max} \cdots 
e_{j_{1}}^{\max}\eta = 
e_{j_{N}}^{\max}e_{j_{N-1}}^{\max} \cdots 
e_{j_{N'+1}}^{\max}\eta' \\
& \hspace*{5mm} = 
e_{j_{N}}^{\max}e_{j_{N-1}}^{\max} \cdots 
e_{j_{N'+1}}^{\max}
(\cl(\Lambda),\cl(\Lambda) \,;\, 0,\sigma,1) = 
(\cl(\lambda),\cl(\lambda) \,;\, 0,\sigma,1) =
\eta_{\cl(\lambda)}.
\quad 
\end{align*}
Therefore, 
the sequence $j_{1},\,j_{2},\,\dots,\,j_{N} \in I$ 
satisfies condition (A). Also, it follows 
from (a), (b), and (d) that
$\bigl(
 r_{j_{p}}r_{j_{p-1}} \cdots r_{j_{1}}(\cl(\lambda))
\bigr)(h_{j_{p+1}}) < 0$ for all $p=0,\,1,\,\dots,\,N-1$,
and hence that 
the sequence $j_{1},\,j_{2},\,\dots,\,j_{N} \in I$ 
satisfies condition (B)'.

\paragraph{Case 2.} $\mu' \ne \cl(\Lambda)$. 

\noindent
In this case, 
we can take (and fix) $j_{N'+1} \in I_{0}$ 
such that $\mu'(h_{j_{N'+1}}) > 0$. 
It is well-known that 
there exists a reduced expression 
of $w_{0} \in \fin{W}$ of the form 
$w_{0}=r_{j_{N'+1}}r_{j_{N'+2}} \cdots r_{j_{N''}}$; 
note that 
$j_{N'+1},\,j_{N'+2},\,\dots,\,j_{N''} \in I_{0}$. 
Since $\Lambda=w_{0}\lambda$, it follows that 

\vspace{1.5mm}

\noindent
(e) \, 
$r_{j_{N''}}r_{j_{N''-1}} \cdots r_{j_{N'+1}}
 (\cl(\Lambda))=\cl(\lambda)$. 

\vspace{1.5mm}

\noindent 
In addition, using \cite[Lemma~3.11\,b)]{Kac}, 
we obtain 

\vspace{1.5mm}

\noindent
(f) \, 
$\bigl(
   r_{j_{p}}r_{j_{p-1}} \cdots r_{j_{N'+1}}(\cl(\Lambda))
 \bigr)
 (h_{j_{p+1}}) < 0$ 
all for $p=N',\,N'+1,\,\dots,\,N''-1$.

\vspace*{1.5mm}

\noindent 
If we set $\eta'':=
e_{j_{N''}}^{\max}e_{j_{N''-1}}^{\max} \cdots 
e_{j_{N'+1}}^{\max}\eta'$, then $||\eta''|| > ||\eta||$.
Indeed, since $(\cl(\Lambda))(h_{j_{N'+1}}) < 0$ 
by (f) with $p=N'$, it follows 
from Lemma~\ref{lem:init1}\,(1) that 
$e_{j_{N'+1}}\eta' \ne \bzero$. 
Also, since $\mu'(h_{j_{N'+1}}) > 0$, it follows 
from Lemma~\ref{lem:final} that 
$f_{j_{N'+1}}\eta' \ne \bzero$. 
Therefore, 
we see from Lemma~\ref{lem:e-max} that
%
%%%%%%%%%%%%%%
%%% eq:cs3 %%%
%%%%%%%%%%%%%%
%
\begin{equation} \label{eq:cs3}
||e_{j_{N'+1}}^{\max}\eta'|| > ||\eta'||.
\end{equation}
Hence we have
\begin{align*}
||\eta''|| & \ge 
||e_{j_{N'+1}}^{\max}\eta'||
\qquad \text{by Lemma~\ref{lem:e-max}} \\
& > ||\eta'|| \qquad \text{by \eqref{eq:cs3}} \\
& \ge ||\eta|| \qquad \text{by \eqref{eq:cs2}}. 
\end{align*}
Furthermore, from (e) and (f), we deduce, 
by applying Lemma~\ref{lem:s1} successively, that 
\begin{align*}
\eta'' & = 
e_{j_{N''}}^{\max}e_{j_{N''-1}}^{\max} \cdots 
e_{j_{N'+1}}^{\max}\eta'
= 
(r_{j_{N''}}r_{j_{N''-1}} \cdots r_{j_{N'+1}}
 (\cl(\Lambda)),\,
 \mu'' \,;\, 0,\,\sigma,\,1) \\
& =(\cl(\lambda),\mu'' \,;\, 0,\,\sigma,\,1)
\end{align*}
for some $\mu'' \in \cl(W\lambda)$. 
Therefore, by the inductive assumption, 
there exists a sequence 
$j_{N''+1},\,j_{N''+2},\,\dots,\,j_{N} \in I$ satisfying 
conditions (A) and (B)' for $\eta''$. 
It is easily checked by (a),\,(b),\,(e),\,(f), and 
the inductive assumption that the sequence
\begin{equation*}
j_{1},\,j_{2},\,\dots,\,j_{N'},\,
 j_{N'+1},\,j_{N'+2},\,\dots,\,j_{N''},\,
 j_{N''+1},\,j_{N''+2},\,\dots,\,j_{N} \in I
\end{equation*}
satisfies conditions (A) and (B)' for $\eta$. 
This proves the sublemma. 
%%%%%%%%%%%%%%%%%%%%%%%%%%

\vsp

Now we turn to the proof of Lemma~\ref{lem:strict1}.
By Remark~\ref{rem:regular}, 
the set of all elements of $[0,1]$ appearing as $\sigma_{u}$'s 
in the reduced expression 
$(\ud{\nu} \,;\, 
  \sigma_{0},\,\sigma_{1},\,\dots,\,\sigma_{s})$ of 
some element of $\BB(\lambda)_{\cl}$, is a finite set. 
Let $\eta =
 (\cl(\lambda),\,\mu_{2},\,\dots,\,\mu_{s} \, ; \, 
  \sigma_{0},\,\sigma_{1},\,\sigma_{2},\,\dots,\,\sigma_{s})$
be the reduced expression of $\eta$. 
Using this fact, 
we show the assertion of 
Lemma~\ref{lem:strict1} by 
descending induction on $\sigma_{1}$ 
in the reduced expression 
$(\cl(\lambda),\,\mu_{2},\,\dots,\,\mu_{s} \, ; \, 
  \sigma_{0},\,\sigma_{1},\,\sigma_{2},\,\dots,\,\sigma_{s})$
of $\eta \in \BB(\lambda)_{\cl}$.
When $\sigma_{1}=1$, $\eta=\eta_{\cl(\lambda)}$, 
since $\iota(\eta)=\cl(\lambda)$ by the assumption of the lemma. 
Thus the assertion obviously holds.
Assume that $\sigma_{1} < 1$, 
or equivalently, $s \ge 2$. 
Set $\eta':=
(\cl(\lambda),\mu_{2} \,;\, 
 \sigma_{0},\,\sigma_{1},\,\sigma_{s})$; 
note that $\eta' \in \BB(\lambda)_{\cl}$ 
by Lemma~\ref{lem:s2}\,(1).
It follows from the sublemma above
that there exists a sequence 
$j_{1},\,j_{2},\,\dots,\,j_{N'} \in I$ 
satisfying conditions (A) and (B)' for $\eta'$. 
We set 
$\eta'':=
 e_{j_{N'}}^{\max}e_{j_{N'-1}}^{\max} \cdots 
 e_{j_{1}}^{\max}\eta$.
Then, repeated use of Lemma~\ref{lem:s2}\,(2) shows that 
%
%%%%%%%%%%%%%%%%%%%%
%%% eq:strict1-1 %%%
%%%%%%%%%%%%%%%%%%%%
%
\begin{equation} \label{eq:strict1-1}
\eta''(t)=
(e_{j_{N'}}^{\max}e_{j_{N'-1}}^{\max} \cdots 
 e_{j_{1}}^{\max}\eta')(t)=\eta_{\cl(\lambda)}(t)=t\cl(\lambda)
\qquad 
\text{for $t \in [0,\sigma_{2}]$}.
\end{equation}
Hence the initial direction $\iota(\eta'')$, 
which equals 
$r_{j_{N'}}r_{j_{N'-1}} \cdots r_{j_{1}}(\cl(\lambda))$ 
by Lemma~\ref{lem:s2}\,(2), must be equal to 
$\cl(\lambda)$. 
Therefore, if 
$\eta''=
 (\cl(\lambda),\,\mu_{2}',\,\dots,\,\mu_{s'}' \, ; \, 
  \sigma_{0}',\,\sigma_{1}',\,\sigma_{2}',\,\dots,\,\sigma_{s'}')$
is the reduced expression of $\eta''$, 
then we see from \eqref{eq:strict1-1} that 
$\sigma_{1}' \ge \sigma_{2} > \sigma_{1}$. 
Hence, by the inductive assumption, 
there exists a sequence 
$j_{N'+1},\,j_{N'+2},\,\dots,\,j_{N} \in I$ 
satisfying conditions (A) and (B)' for $\eta''$.
Thus we obtain a sequence
\begin{equation*}
j_{1},\,j_{2},\,\dots,\,j_{N'},\,
j_{N'+1},\,j_{N'+2},\,\dots,\,j_{N} \in I
\end{equation*}
satisfying conditions (A) and (B)' for $\eta$ 
(note that 
 $r_{j_{N'}}r_{j_{N'-1}} \cdots r_{j_{1}}(\cl(\lambda))=
 \iota(\eta'')=\cl(\lambda)$). 
This completes the proof of Lemma~\ref{lem:strict1}. 
\end{proof}

Since $\lambda$ is assumed to be 
strictly level-zero dominant, 
we see from Lemma~\ref{lem:fp}\,(3) that 
the condition~(B)' of Lemma~\ref{lem:strict1} can be replaced 
with condition~(B) of Proposition~\ref{prop:key}. Namely, we have 
%
%%%%%%%%%%%%%%%%%%%
%%% lem:strict2 %%%
%%%%%%%%%%%%%%%%%%%
%
\begin{lem} \label{lem:strict2}
Let $\lambda \in \sum_{i \in I_{0}}\BZ_{> 0} \vpi_{i}$ 
be a strictly level-zero dominant integral weight.
Let $\eta$ be an element of $\BB(\lambda)_{\cl}$ 
such that $\iota(\eta)=\cl(\lambda)$. 
Then, there exists a sequence 
$j_{1},\,j_{2},\,\dots,\,j_{N} \in I$ satisfying 
conditions {\rm (A)} and {\rm (B)} of 
Proposition~\ref{prop:key}.
\end{lem}

Finally, let us give 
a proof of Proposition~\ref{prop:key}.
\begin{proof}[Proof of Proposition~\ref{prop:key}]
We set $\rho:=\sum_{i \in I_{0}} \vpi_{i} \in P^{0}_{+}$; 
note that $\rho$ is strictly level-zero dominant, 
and hence so is $\lambda+\rho \in 
\sum_{i \in I} \BZ_{> 0} \vpi_{i}$.
Since $e_{j}\eta = \bzero$ for all $j \in I_{0}$ by assumption 
and $e_{j} \eta_{\cl(\rho)} = \bzero$ for all $j \in I_{0}$ 
by the definition of the root operators $e_{j}$, 
we see from the tensor product rule for crystals that 
$\eta \otimes \eta_{\cl(\rho)} \in 
 \BB(\lambda)_{\cl} \otimes \BB(\rho)_{\cl}$ 
also satisfies the condition that 
$e_{j}(\eta \otimes \eta_{\cl(\rho)}) = \bzero$ 
for all $j \in I_{0}$. 
Recall from Corollary~\ref{cor:tpd}\,(1) that 
there exists an isomorphism $\Psi_{\lambda,\rho}:
\BB(\lambda+\rho)_{\cl} \stackrel{\sim}{\rightarrow} 
\BB(\lambda)_{\cl} \otimes \BB(\rho)_{\cl}$ of 
$P_{\cl}$-crystals. 
Set $\eta':=
 \Psi_{\lambda,\rho}^{-1}(\eta \otimes \eta_{\cl(\rho)}) \in 
 \BB(\lambda+\rho)_{\cl}$. 
Then, clearly $e_{j}\eta' = \bzero$ for all $j \in I_{0}$.
Hence it follows from Lemma~\ref{lem:init1}\,(1) 
(see also Remark~\ref{rem:cl-expr}) 
that the initial direction $\iota(\eta') \in 
\cl(W(\lambda+\rho))=\fin{W}\cl(\lambda+\rho)$ of 
$\eta'$ is level-zero dominant, and hence 
from Remark~\ref{rem:orbcl} that 
$\iota(\eta')=\cl(\lambda+\rho)$.
Since $\lambda+\rho \in \sum_{i \in I} \BZ_{> 0} \vpi_{i}$ 
is strictly level-zero dominant, 
we know from Lemma~\ref{lem:strict2} that 
there exists 
a sequence $j_{1},\,j_{2},\,\dots,\,j_{N} \in I$ 
satisfying conditions (A) and (B) for $\eta'$.
It remains to show that 
$e_{j_{N}}^{\max}e_{j_{N-1}}^{\max} \cdots
 e_{j_{1}}^{\max}\eta = \eta_{\cl(\lambda)}$, 
i.e., that 
the sequence 
$j_{1},\,j_{2},\,\dots,\,j_{N} \in I$ also 
satisfies condition (A) for $\eta$.
We deduce that 
\begin{align*}
& 
e_{j_{N}}^{\max}e_{j_{N-1}}^{\max} \cdots 
     e_{j_{1}}^{\max}(\eta \otimes \eta_{\cl(\rho)}) 
= 
e_{j_{N}}^{\max}e_{j_{N-1}}^{\max} \cdots 
e_{j_{1}}^{\max} \bigl(\Psi_{\lambda,\rho}(\eta')\bigr) \\
& \hspace*{5mm}
= 
\Psi_{\lambda,\rho}
\bigl(e_{j_{N}}^{\max}e_{j_{N-1}}^{\max} \cdots 
      e_{j_{1}}^{\max} \eta'\bigr)
=
\Psi_{\lambda,\rho}(\eta_{\cl(\lambda+\rho)})
\quad \text{by condition (A) for $\eta'$} \\
& \hspace*{5mm}
= 
\eta_{\cl(\lambda)} \otimes \eta_{\cl(\rho)}
\qquad \text{by \eqref{eq:map2}}.
\end{align*}
Also, we see from Lemma~\ref{lem:max}\,(2) that 
\begin{equation*}
e_{j_{N}}^{\max}e_{j_{N-1}}^{\max} \cdots 
     e_{j_{1}}^{\max}(\eta \otimes \eta_{\cl(\rho)}) = 
(e_{j_{N}}^{\max}e_{j_{N-1}}^{\max} \cdots 
     e_{j_{1}}^{\max}\eta) \otimes \eta''
\quad \text{for some $\eta'' \in \BB(\rho)_{\cl}$}.
\end{equation*}
Thus, we obtain 
$e_{j_{N}}^{\max}e_{j_{N-1}}^{\max} \cdots
 e_{j_{1}}^{\max}\eta = \eta_{\cl(\lambda)}$.
This establishes Proposition~\ref{prop:key}.
\end{proof}
%
%
%
%==============================%
%     START SUBSECTION 0404    %
%==============================%
%
\subsection{Proof of Theorem~\ref{thm:step1}.}
\label{subsec:step1}

Let $\lambda,\,\lambda' \in P^{0}_{+}$ be 
level-zero dominant integral weights.
For $\eta_{1} \otimes \eta_{2} \in 
 \BB(\lambda)_{\cl} \otimes \BB(\lambda')_{\cl}$, 
we define (see Theorem~\ref{thm:energy})
%
%%%%%%%%%%%%%%%%
%%% eq:dfn-D %%%
%%%%%%%%%%%%%%%%
%
\begin{equation} \label{eq:dfn-D}
D_{\lambda,\lambda'}(\eta_{1} \otimes \eta_{2}):=
 H_{\lambda,\lambda'}(\eta_{1} \otimes \eta_{2})+
 \Deg_{\lambda}(\eta_{1})+\Deg_{\lambda'}(\ti{\eta}_{2}),
\end{equation}
where $\ti{\eta}_{2} \in \BB(\lambda')_{\cl}$ is defined by: 
$\ti{\eta}_{2} \otimes \ti{\eta}_{1}=
 R_{\lambda,\lambda'}(\eta_{1} \otimes \eta_{2})$ 
 (see Corollary~\ref{cor:tpd}\,(2)). 
From Theorem~\ref{thm:energy}\,(H2) and 
Lemma~\ref{lem:deg-ro}\,(1), 
we deduce by use of \eqref{eq:map3} that 
%
%
%%%%%%%%%%%%%%%%
%%% eq:D-ext %%%
%%%%%%%%%%%%%%%%
%
\begin{equation} \label{eq:D-ext}
D_{\lambda,\lambda'}
 (\eta_{\cl(\lambda)} \otimes \eta_{\cl(\lambda')})=
H_{\lambda,\lambda'}
(\eta_{\cl(\lambda)} \otimes \eta_{\cl(\lambda')})+
\Deg_{\lambda}(\eta_{\cl(\lambda)})+
\Deg_{\lambda'}(\eta_{\cl(\lambda')})=0.
\end{equation}
%
%%%%%%%%%%%%%%%%%%%
%%% lem:eng-max %%%
%%%%%%%%%%%%%%%%%%%
%
\begin{lem} \label{lem:eng-max}
Let $\lambda,\,\lambda' \in P^{0}_{+}$ 
be level-zero dominant integral weights, and 
let $\eta_{1} \otimes \eta_{2} \in 
 \BB(\lambda)_{\cl} \otimes \BB(\lambda')_{\cl}$.

\noindent {\rm (1)} 
$D_{\lambda,\lambda'}(e_{j}^{\max}(\eta_{1} \otimes \eta_{2})) = 
 D_{\lambda,\lambda'}(\eta_{1} \otimes \eta_{2})$ for all $j \in I_{0}$.

\noindent {\rm (2)} 
If $(\iota(\eta_{1}))(h_{0}) \le 0$ and 
$(\iota(\ti{\eta}_{2}))(h_{0}) \le 0$, then 
%
%%%%%%%%%%%%%%%%%%
%%% eq:eng-max %%%
%%%%%%%%%%%%%%%%%%
%
\begin{equation} \label{eq:eng-max}
D_{\lambda,\lambda'}
   (e_{0}^{\max}(\eta_{1} \otimes \eta_{2})) = 
D_{\lambda,\lambda'}
 (\eta_{1} \otimes \eta_{2})
  -(\iota(\eta_{1})+\iota(\ti{\eta}_{2}))(h_{0})
  -\ve_{0}(\eta_{1} \otimes \eta_{2}). 
\end{equation}
\end{lem}

\begin{proof}
Let $j \in I$. 
We see from Lemma~\ref{lem:max}\,(2) that
\begin{equation*}
e_{j}^{\max}(\eta_{1} \otimes \eta_{2}) = 
e_{j}^{\max}\eta_{1} \otimes \eta_{2}'
\quad
\text{for some $\eta_{2}' \in \BB(\lambda')_{\cl}$}, 
\end{equation*}
and hence that
\begin{align*}
R_{\lambda,\lambda'}(e_{j}^{\max}(\eta_{1} \otimes \eta_{2})) 
& =
  e_{j}^{\max}R_{\lambda,\lambda'}(\eta_{1} \otimes \eta_{2}) 
  =
  e_{j}^{\max}(\ti{\eta}_{2} \otimes \ti{\eta}_{1}) \\
& = 
  e_{j}^{\max}\ti{\eta}_{2} \otimes \eta_{1}'
\quad
  \text{for some $\eta_{1}' \in \BB(\lambda)_{\cl}$.}
\end{align*}
Therefore, from the definition \eqref{eq:dfn-D} 
of $D_{\lambda,\lambda'}$, we obtain
%
%%%%%%%%%%%%%%%%%%%
%%% eq:eng-max0 %%%
%%%%%%%%%%%%%%%%%%%
%
\begin{equation} \label{eq:eng-max0}
D_{\lambda,\lambda'}
 (e_{j}^{\max}(\eta_{1} \otimes \eta_{2}))=
 H_{\lambda,\lambda'}
 (e_{j}^{\max}(\eta_{1} \otimes \eta_{2}))+
 \Deg_{\lambda}(e_{j}^{\max}\eta_{1})+
 \Deg_{\lambda'}(e_{j}^{\max}\ti{\eta}_{2}).
\end{equation}
From \eqref{eq:eng-max0}, 
part (1) follows immediately 
by Theorem~\ref{thm:energy}\,(H1) and 
Lemma~\ref{lem:deg-ro}\,(2).
Let us prove part~(2).
We give a proof only for the case
$\ve_{0}(\eta_{1}) \ge 
 \ve_{0}(\ti{\eta}_{2})$; 
the proof for the case
$\ve_{0}(\eta_{1}) < 
 \ve_{0}(\ti{\eta}_{2})$ is similar. 
By Lemma~\ref{lem:deg-max}, 
%
%%%%%%%%%%%%%%%%%%%
%%% eq:eng-max1 %%%
%%%%%%%%%%%%%%%%%%%
%
\begin{equation} \label{eq:eng-max1}
\begin{cases}
\Deg_{\lambda}(e_{0}^{\max}\eta_{1})=
  \Deg_{\lambda}(\eta_{1})-
  \ve_{0}(\eta_{1})-(\iota(\eta_{1}))(h_{0}), & \\[3mm]
\Deg_{\lambda'}(e_{0}^{\max}\ti{\eta}_{2})=
  \Deg_{\lambda'}(\ti{\eta}_{2})-
  \ve_{0}(\ti{\eta}_{2})-(\iota(\ti{\eta}_{2}))(h_{0}). &
\end{cases}
\end{equation}
For simplicity of notation,
we set $L:=\ve_{0}(\eta_{1} \otimes \eta_{2})$; 
note that $L=\ve_{0}(\ti{\eta}_{2} \otimes \ti{\eta}_{1})$ 
since $\ti{\eta}_{2} \otimes \ti{\eta}_{1}=
R_{\lambda,\lambda'}(\eta_{1} \otimes \eta_{2})$, and that 
$L \ge \ve_{0}(\eta_{1}),\,\ve_{0}(\ti{\eta}_{2})$ 
by Lemma~\ref{lem:max}\,(1).
It follows from Lemma~\ref{lem:max}\,(2) that 
for $0 \le l \le L$, 
\begin{equation*}
e_{0}^{l}(\eta_{1} \otimes \eta_{2}) = 
\begin{cases}
\eta_{1} \otimes e_{0}^{l}\eta_{2} & 
\text{if $0 \le l \le L-\ve_{0}(\eta_{1})$}, \\[3mm]
e_{0}^{l-L+\ve_{0}(\eta_{1})} \eta_{1} \otimes 
e_{0}^{L-\ve_{0}(\eta_{1})}\eta_{2} & 
\text{if $L-\ve_{0}(\eta_{1}) \le l \le L$},
\end{cases}
\end{equation*}
and 
\begin{equation*}
e_{0}^{l}(\ti{\eta}_{2} \otimes \ti{\eta}_{1}) = 
\begin{cases}
\ti{\eta}_{2} \otimes e_{0}^{l}\ti{\eta}_{1} & 
\text{if $0 \le l \le L-\ve_{0}(\ti{\eta}_{2})$}, \\[3mm]
e_{0}^{l-L+\ve_{0}(\ti{\eta}_{2})} \ti{\eta}_{2} \otimes 
e_{0}^{L-\ve_{0}(\ti{\eta}_{2})}\ti{\eta}_{1} & 
\text{if $L-\ve_{0}(\ti{\eta}_{2}) \le l \le L$}.
\end{cases}
\end{equation*}
Therefore, 
we deduce from Theorem~\ref{thm:energy}\,(H1) 
that for $0 \le l' \le L-\ve_{0}(\eta_{1})$, 
%
%%%%%%%%%%%%%%%
%%% eq:le01 %%%
%%%%%%%%%%%%%%%
%
\begin{equation} \label{eq:le01}
H_{\lambda,\lambda'}(e_{0}^{l'}(\eta_{1} \otimes \eta_{2}))= 
H_{\lambda,\lambda'}(\eta_{1} \otimes \eta_{2})-l'. 
\end{equation}
Similarly, 
we deduce from Theorem~\ref{thm:energy}\,(H1) that 
for $0 \le l' \le \ve_{0}(\eta_{1})-\ve_{0}(\ti{\eta}_{2})$, 
%
%%%%%%%%%%%%%%%
%%% eq:le02 %%%
%%%%%%%%%%%%%%%
%
\begin{align}
& 
H_{\lambda,\lambda'}
 (e_{0}^{L-\ve_{0}(\eta_{1})+l'}(\eta_{1} \otimes \eta_{2}))= 
H_{\lambda,\lambda'}
 (e_{0}^{L-\ve_{0}(\eta_{1})}(\eta_{1} \otimes \eta_{2})) \nonumber \\
& \hspace*{5mm} =
H_{\lambda,\lambda'}
 (\eta_{1} \otimes \eta_{2})-L+\ve_{0}(\eta_{1})
\quad
\text{by \eqref{eq:le01} with $l'=L-\ve_{0}(\eta_{1})$}, \label{eq:le02}
\end{align}
and that for $0 \le l' \le \ve_{0}(\ti{\eta}_{2})$, 
%
%%%%%%%%%%%%%%%
%%% eq:le03 %%%
%%%%%%%%%%%%%%%
%
\begin{align}
& 
H_{\lambda,\lambda'}
 (e_{0}^{L-\ve_{0}(\ti{\eta}_{2})+l'}(\eta_{1} \otimes \eta_{2}))=
H_{\lambda,\lambda'}
 (e_{0}^{L-\ve_{0}(\ti{\eta}_{2})}(\eta_{1} \otimes \eta_{2}))+l' \nonumber \\
& \hspace*{5mm} 
= 
H_{\lambda,\lambda'}
 (\eta_{1} \otimes \eta_{2})-L+\ve_{0}(\eta_{1})+l'
\quad 
\text{by \eqref{eq:le02} with 
$l'=\ve_{0}(\eta_{1})-\ve_{0}(\ti{\eta}_{2})$}. \label{eq:le03}
\end{align}
Finally, by taking 
$l'=\ve_{0}(\ti{\eta}_{2})$ in \eqref{eq:le03}, 
we conclude that 
%
%%%%%%%%%%%%%%%%%%%
%%% eq:eng-max2 %%%
%%%%%%%%%%%%%%%%%%%
%
\begin{align} \label{eq:eng-max2}
H_{\lambda,\lambda'}
 (e_{0}^{\max}(\eta_{1} \otimes \eta_{2})) & = 
H_{\lambda,\lambda'}
 (e_{0}^{L}(\eta_{1} \otimes \eta_{2})) \nonumber \\
& = 
H_{\lambda,\lambda'}
 (\eta_{1} \otimes \eta_{2})-L+
  \ve_{0}(\eta_{1})+\ve_{0}(\ti{\eta}_{2}) \nonumber \\
& = 
H_{\lambda,\lambda'}
 (\eta_{1} \otimes \eta_{2})-\ve_{0}(\eta_{1} \otimes \eta_{2})+
  \ve_{0}(\eta_{1})+\ve_{0}(\ti{\eta}_{2}).
\end{align}
By substituting \eqref{eq:eng-max1} and \eqref{eq:eng-max2} 
into \eqref{eq:eng-max0}, and then comparing 
the resulting equation with \eqref{eq:dfn-D}, 
we obtain the desired equation 
\eqref{eq:eng-max}. This proves the lemma.
\end{proof}
%
%%%%%%%%%%%%%%%%
%%% prop:eng %%%
%%%%%%%%%%%%%%%%
%
\begin{prop} \label{prop:eng}
Let 
$\lambda,\,\lambda' \in P^{0}_{+}$ be 
level-zero dominant integral weights, and let 
$\eta \in \BB(\lambda+\lambda')_{\cl}$. 
Then, 
%
%%%%%%%%%%%%%%
%%% eq:eng %%%
%%%%%%%%%%%%%%
%
\begin{equation} \label{eq:eng}
\Deg_{\lambda+\lambda'}(\eta)=
 D_{\lambda,\lambda'}(\Psi_{\lambda,\lambda'}(\eta)), 
\end{equation}
where $\Psi_{\lambda,\lambda'}:\BB(\lambda+\lambda')_{\cl} 
\stackrel{\sim}{\rightarrow} 
\BB(\lambda)_{\cl} \otimes \BB(\lambda')_{\cl}$ is 
the isomorphism of $P_{\cl}$-crystals 
in Corollary~\ref{cor:tpd}\,{\rm(1)}.
\end{prop}

\begin{proof}
If $j \in I_{0}$, then it follows from 
Lemmas~\ref{lem:deg-ro}\,(2) and 
\ref{lem:eng-max}\,(1) that 
$\Deg_{\lambda+\lambda'}(e_{j}^{\max}\eta)=
 \Deg_{\lambda+\lambda'}(\eta)$ and 
$D_{\lambda, \lambda'}(\Psi_{\lambda,\lambda'}(e_{j}^{\max}\eta))=
 D_{\lambda, \lambda'}\bigl(e_{j}^{\max}(\Psi_{\lambda,\lambda'}(\eta))\bigr)=
 D_{\lambda, \lambda'}(\Psi_{\lambda,\lambda'}(\eta))$. 
Therefore, we may assume that 
$\eta \in \BB(\lambda+\lambda')_{\cl}$ satisfies 
the condition that $e_{j}\eta = \bzero$ 
for all $j \in I_{0}$, since the $P_{\cl}$-crystal 
$\BB(\lambda+\lambda')_{\cl}$ is regular. 
Then it follows from 
Lemma~\ref{lem:init1}\,(1) 
(see also Remark~\ref{rem:cl-expr}) that
$\iota(\eta) \in \cl(W\lambda)=\fin{W}\cl(\lambda)$ is 
level-zero dominant, and hence 
from Remark~\ref{rem:orbcl} that 
$\iota(\eta)=\cl(\lambda+\lambda')$.
Set $\eta_{1} \otimes \eta_{2}:=
 \Psi_{\lambda,\lambda'}(\eta) \in 
 \BB(\lambda)_{\cl} \otimes \BB(\lambda')_{\cl}$, and 
$\ti{\eta}_{2} \otimes \ti{\eta}_{1}:=
 R_{\lambda,\lambda'}(\eta_{1} \otimes \eta_{2})$. 
Since $e_{j}\eta = \bzero$ for all $j \in I_{0}$, 
we have 
$e_{j}(\eta_{1} \otimes \eta_{2})=\bzero$ and 
$e_{j}(\ti{\eta}_{2} \otimes \ti{\eta}_{1})=\bzero$
for all $j \in I_{0}$, which implies that  
$e_{j}\eta_{1}=\bzero$ and 
$e_{j}\ti{\eta}_{2}=\bzero$
for all $j \in I_{0}$ by Lemma~\ref{lem:max}\,(1). 
Hence an argument similar to the above shows that 
$\iota(\eta_{1})=\cl(\lambda)$ and 
$\iota(\ti{\eta}_{2})=\cl(\lambda')$.

By Proposition~\ref{prop:key}, there exists a sequence 
$j_{1},\,j_{2},\,\dots,\,j_{N} \in I$ satisfying 
conditions~(A) and (B) 
for $\eta \in \BB(\lambda+\lambda')_{\cl}$. 
Condition (B) implies that 
%
%%%%%%%%%%%%%%%%%
%%% eq:eng1-1 %%%
%%%%%%%%%%%%%%%%%
%
\begin{equation} \label{eq:eng1-1}
\bigl(
 w^{(p)}(\cl(\lambda+\lambda'))
\bigr)
  (h_{j_{p+1}}) \le 0
\quad
\text{for $p=0,\,1,\,\dots,\,N-1$},
\end{equation}
by Lemma~\ref{lem:fp}\,(2), where we set 
$w^{(p)}:=r_{j_{p}}r_{j_{p-1}} \cdots r_{j_{1}}$
for $0 \le p \le N$.
Therefore, we see from Lemma~\ref{lem:init2} that 
%
%%%%%%%%%%%%%%%%%
%%% eq:eng1-2 %%%
%%%%%%%%%%%%%%%%%
%
\begin{equation} \label{eq:eng1-2}
\iota (E^{(p)}\eta)
  =
w^{(p)}(\cl(\lambda+\lambda'))
\quad
\text{for $p=0,\,1,\,\dots,\,N$},
\end{equation}
where we set 
$E^{(p)}:=
 e_{j_{p}}^{\max}
 e_{j_{p-1}}^{\max} \cdots 
 e_{j_{1}}^{\max}$ 
for $0 \le p \le N$.
Because $\iota(\eta_{1})=\cl(\lambda)$ and 
$\iota(\ti{\eta}_{2})=\cl(\lambda')$, 
an argument similar to the above shows that
%
%%%%%%%%%%%%%%%
%%% eq:eng2 %%%
%%%%%%%%%%%%%%%
%
\begin{equation} \label{eq:eng2}
\begin{cases}
  \iota (E^{(p)}\eta_{1})
  =
  w^{(p)}(\cl(\lambda))
  & \text{for $p=0,\,1,\,\dots,\,N$}, \\[3mm]
  \bigl(
  w^{(p)}(\cl(\lambda))
  \bigr)
  (h_{j_{p+1}}) \le 0
  & \text{for $p=0,\,1,\,\dots,\,N-1$},
\end{cases}
\end{equation}
and 
%
%%%%%%%%%%%%%%%
%%% eq:eng3 %%%
%%%%%%%%%%%%%%%
%
\begin{equation} \label{eq:eng3}
\begin{cases}
  \iota (E^{(p)}\ti{\eta}_{2})
  =
  w^{(p)}(\cl(\lambda'))
  & \text{for $p=0,\,1,\,\dots,\,N$}, \\[3mm]
  \bigl(
  w^{(p)}(\cl(\lambda'))
  \bigr)
  (h_{j_{p+1}}) \le 0
  & \text{for $p=0,\,1,\,\dots,\,N-1$}.
\end{cases}
\end{equation}

Now we deduce that
\begin{align*}
0 & 
 =\Deg_{\lambda+\lambda'}(\eta_{\cl(\lambda+\lambda')})
\quad \text{by Lemma~\ref{lem:deg-ro}\,(1)} \\[1mm]
& =\Deg_{\lambda+\lambda'}
  (e_{j_{N}}^{\max}e_{j_{N-1}}^{\max} \cdots e_{j_{1}}^{\max}\eta)
\quad \text{by condition (A) for $\eta$} \\[1mm]
& =\Deg_{\lambda+\lambda'}(\eta)-
  \sum_{1 \le p \le N \,;\,j_{p}=0}
  \left\{\ve_{0}(E^{(p-1)}\eta)+
  \bigl(w^{(p-1)}(\cl(\lambda+\lambda'))\bigr)(h_{j_{p}})
  \right\} \\
& \hspace*{40mm} 
  \text{by Lemma~\ref{lem:deg-max} 
  along with \eqref{eq:eng1-1}, \eqref{eq:eng1-2}},
\end{align*}
and hence that 
%
%%%%%%%%%%%%%%%
%%% eq:eng4 %%%
%%%%%%%%%%%%%%%
%
\begin{equation} \label{eq:eng4}
\Deg_{\lambda+\lambda'}(\eta) = 
  \sum_{1 \le p \le N \,;\,j_{p}=0}
  \left\{\ve_{0}(E^{(p-1)}\eta)+
  \bigl(w^{(p-1)}(\cl(\lambda+\lambda'))\bigr)(h_{j_{p}})
  \right\}.
\end{equation}
Here, by Lemma~\ref{lem:max}\,(2), we have
\begin{align*}
& E^{(p)}(\eta_{1} \otimes \eta_{2})=
  E^{(p)}\eta_{1} \otimes \eta_{2}'
\quad \text{for some $\eta_{2}' \in \BB(\lambda')_{\cl}$}, \\
& R_{\lambda,\lambda'}(E^{(p)}(\eta_{1} \otimes \eta_{2}))=
  E^{(p)}(R_{\lambda,\lambda'}(\eta_{1} \otimes \eta_{2})) \\
& \hspace*{10mm} =
  E^{(p)}(\ti{\eta}_{2} \otimes \ti{\eta}_{1})=
  E^{(p)}\ti{\eta}_{2} \otimes \eta_{1}'
\quad \text{for some $\eta_{1}' \in \BB(\lambda)_{\cl}$}.
\end{align*}
Therefore, we deduce that 
\begin{align*}
0 & 
 =D_{\lambda,\lambda'}
    (\eta_{\cl(\lambda)} \otimes \eta_{\cl(\lambda')})
\quad \text{by \eqref{eq:D-ext}} \\[1mm]
& 
 =D_{\lambda,\lambda'}
  (e_{j_{N}}^{\max}e_{j_{N-1}}^{\max} \cdots e_{j_{1}}^{\max}
   (\eta_{1} \otimes \eta_{2})) 
\quad 
\text{by condition (A) for 
  $\eta$ and \eqref{eq:map2}} \\[1mm]
&
 =D_{\lambda,\lambda'}(\eta_{1} \otimes \eta_{2})-
  \sum_{1 \le p \le N \,;\,j_{p}=0}
  \left\{\ve_{0}(E^{(p-1)}(\eta_{1} \otimes \eta_{2}))+
  \bigl(w^{(p-1)}(\cl(\lambda+\lambda'))\bigr)(h_{j_{p}})
  \right\} \\
& \hspace*{50mm}  
  \text{by Lemma~\ref{lem:eng-max} along with 
  \eqref{eq:eng2}, \eqref{eq:eng3}}.
\end{align*}
Since 
$\Psi_{\lambda,\lambda'}(E^{(p)}\eta)=
 E^{(p)}\Psi_{\lambda,\lambda'}(\eta)=
 E^{(p)}(\eta_{1} \otimes \eta_{2})$, 
it follows that 
$\ve_{0}(E^{(p)}\eta)=
 \ve_{0}(E^{(p)}(\eta_{1} \otimes \eta_{2}))$ for all 
$p=0,\,1,\,\dots,\,N$. 
Thus we obtain 
%
%%%%%%%%%%%%%%%
%%% eq:eng5 %%%
%%%%%%%%%%%%%%%
%
\begin{equation} \label{eq:eng5}
D_{\lambda,\lambda'}(\eta_{1} \otimes \eta_{2}) = 
  \sum_{1 \le p \le N \,;\,j_{p}=0}
  \left\{\ve_{0}(E^{(p-1)}\eta)+
  \bigl(w^{(p-1)}(\cl(\lambda+\lambda'))\bigr)(h_{j_{p}})
  \right\}.
\end{equation}
Equation~\eqref{eq:eng} follows immediately 
from \eqref{eq:eng4} and \eqref{eq:eng5}. 
This completes the proof of the proposition.
\end{proof}

Now we are ready to prove Theorem~\ref{thm:step1}. 

\begin{proof}[Proof of Theorem~\ref{thm:step1}]
We proceed by induction on the length $n$ of 
the sequence $\bi=(i_{1},\,i_{2},\,\dots,\,i_{n})$. 
When $n=1$, the assertion obviously holds.
Assume that $n > 1$, and set 
$\bi':=(i_{1},\,i_{2},\,\dots,\,i_{n-1})$, 
$\lambda':=\lambda-\vpi_{i_{n}} \in P^{0}_{+}$.
Recall from Theorem~\ref{thm:tpd} and 
Corollary~\ref{cor:tpd}\,(1)
that there exist isomorphisms 
$\Psi_{\bi'}:
 \BB(\lambda')_{\cl} 
 \stackrel{\sim}{\rightarrow} 
 \BB_{\bi'}=
 \BB(\vpi_{i_{1}})_{\cl} \otimes 
 \BB(\vpi_{i_{2}})_{\cl} \otimes \cdots \otimes 
 \BB(\vpi_{i_{n-1}})_{\cl}$
and 
$\Psi_{\lambda',\vpi_{i_{n}}}:
 \BB(\lambda)_{\cl} 
 \stackrel{\sim}{\rightarrow}
 \BB(\lambda')_{\cl} \otimes \BB(\vpi_{i_{n}})_{\cl}$ 
of $P_{\cl}$-crystals. 
Let $\eta \in \BB(\lambda)_{\cl}$, and set
$\eta' \otimes \eta'':=\Psi_{\lambda',\vpi_{i_{n}}}(\eta) \in 
 \BB(\lambda')_{\cl} \otimes \BB(\vpi_{i_{n}})_{\cl}$, 
$\eta_{1} \otimes \eta_{2} \otimes \cdots \otimes \eta_{n}:=
 \Psi_{\bi}(\eta) \in \BB_{\bi}=
 \BB(\vpi_{i_{1}})_{\cl} \otimes 
 \BB(\vpi_{i_{2}})_{\cl} \otimes \cdots \otimes 
 \BB(\vpi_{i_{n}})_{\cl}$. 
Note that both of $\Psi_{\bi}$ and 
$(\Psi_{\bi'} \otimes \id) \circ 
 \Psi_{\lambda',\vpi_{i_{n}}}$ are
isomorphisms of $P_{\cl}$-crystals from 
$\BB(\lambda)_{\cl}$ to 
$\BB_{\bi}$.
Since the $P_{\cl}$-crystals $\BB(\lambda)_{\cl}$ and $\BB_{\bi}$ 
are both simple, it follows from Lemma~\ref{lem:simple}\,(3) that 
$\Psi_{\bi}=(\Psi_{\bi'} \otimes \id) \circ 
 \Psi_{\lambda',\vpi_{i_{n}}}$, and hence that 
%
%%%%%%%%%%%%%%%
%%% eq:s1-1 %%%
%%%%%%%%%%%%%%%
%
\begin{equation} \label{eq:s1-1}
\eta''=\eta_{n}
\quad
\text{and}
\quad
\Psi_{\bi'}(\eta')=
\eta_{1} \otimes \eta_{2} \otimes \cdots \otimes \eta_{n-1}.
\end{equation}

We see from Proposition~\ref{prop:eng} that 
%
%%%%%%%%%%%%%%%%
%%% eq:s1-15 %%%
%%%%%%%%%%%%%%%%
%
\begin{equation} \label{eq:s1-15}
\Deg_{\lambda}(\eta)=
 D_{\lambda',\vpi_{i_{n}}}(\eta' \otimes \eta'')=
 H_{\lambda',\vpi_{i_{n}}}(\eta' \otimes \eta'') + 
 \Deg_{\lambda'}(\eta')+\Deg_{\vpi_{i_{n}}}(\ti{\eta}''),
\end{equation}
where we set 
$\ti{\eta}'' \otimes \ti{\eta}':=
 R_{\lambda',\vpi_{i_{n}}}(\eta' \otimes \eta'')$.
Here we remark 
(see the definition of $\eta_{n}^{(1)}$ in \S\ref{subsec:main}) 
that the element $\eta_{n}^{(1)} \in \BB(\vpi_{i_{n}})_{\cl}$ 
is the first factor of the image of 
$\eta_{1} \otimes 
 \eta_{2} \otimes \cdots \otimes 
 \eta_{n} \in \BB_{\bi}$ 
under the (unique) isomorphism 
$(\id \otimes \Psi_{\bi'}) \circ 
R_{\lambda',\vpi_{i_{n}}} \circ 
\Psi_{\lambda',\vpi_{i_{n}}} \circ 
\Psi_{\bi}^{-1}$ of 
$P_{\cl}$-crystals, which is obtained as follows: 
\begin{align*}
& 
\begin{CD}
 \BB(\vpi_{i_{1}})_{\cl} \otimes 
 \BB(\vpi_{i_{2}})_{\cl} \otimes \cdots \otimes 
 \BB(\vpi_{i_{n}})_{\cl}=\BB_{\bi}
& @>{\Psi_{\bi}^{-1}}>{\sim}> & 
\BB(\lambda)_{\cl}
& @>{\Psi_{\lambda',\vpi_{i_{n}}}}>{\sim}> & 
\BB(\lambda')_{\cl} \otimes \BB(\vpi_{i_{n}})_{\cl}
\end{CD}
\\[2mm]
& 
\begin{CD}
@>{R_{\lambda',\vpi_{i_{n}}}}>{\sim}> &
\BB(\vpi_{i_{n}})_{\cl} \otimes \BB(\lambda')_{\cl} & 
@>{\id \otimes \Psi_{\bi'}}>{\sim}> & 
\BB(\vpi_{i_{n}})_{\cl} \otimes 
\BB(\vpi_{i_{1}})_{\cl} \otimes \cdots \otimes 
\BB(\vpi_{i_{n-1}})_{\cl}.
\end{CD}
\end{align*}
Also, it is easy to check that 
this first factor is equal to 
$\ti{\eta}'' \in \BB(\vpi_{i_{n}})_{\cl}$. 
Thus we have $\eta_{n}^{(1)}=\ti{\eta}''$. 
Furthermore, using \eqref{eq:s1-1}, 
we see from \cite[Lemma~5.2]{O} that 
%
%%%%%%%%%%%%%%%
%%% eq:s1-2 %%%
%%%%%%%%%%%%%%%
%
\begin{equation} \label{eq:s1-2}
H_{\lambda',\vpi_{i_{n}}}(\eta' \otimes \eta'') = 
\sum_{1 \le k \le n-1} 
H_{\vpi_{i_{k}},\vpi_{i_{n}}}(\eta_{k} \otimes \eta_{n}^{(k+1)}).
\end{equation}
Now, the inductive assumption 
along with \eqref{eq:s1-1} implies that
%
%%%%%%%%%%%%%%%
%%% eq:s1-3 %%%
%%%%%%%%%%%%%%%
%
\begin{equation} \label{eq:s1-3}
\Deg_{\lambda'}(\eta') = 
\sum_{1 \le k < l \le n-1} 
 H_{\vpi_{i_{k}},\vpi_{i_{l}}}
 (\eta_{k} \otimes \eta_{l}^{(k+1)}) + 
 \sum_{k=1}^{n-1} \Deg_{\vpi_{i_{k}}} (\eta_{k}^{(1)}).
\end{equation}
By substituting \eqref{eq:s1-2} and \eqref{eq:s1-3} 
into \eqref{eq:s1-15}, and 
using the fact that 
$\eta_{n}^{(1)}=\ti{\eta}''$, 
we obtain
\begin{align*}
\Deg_{\lambda}(\eta) & =
H_{\lambda',\vpi_{i_{n}}}(\eta' \otimes \eta'') + 
\Deg_{\lambda'}(\eta')+\Deg_{\vpi_{i_{n}}}(\ti{\eta}'') \\
& = 
 \sum_{1 \le k \le n-1} 
 H_{\vpi_{i_{k}},\vpi_{i_{n}}}(\eta_{k} \otimes \eta_{n}^{(k+1)}) \\
& \hspace*{10mm} + 
 \sum_{1 \le k < l \le n-1} 
 H_{\vpi_{i_{k}},\vpi_{i_{l}}}(\eta_{k} \otimes \eta_{l}^{(k+1)}) +
 \sum_{k=1}^{n-1} \Deg_{\vpi_{i_{k}}} (\eta_{k}^{(1)})+
 \Deg_{\vpi_{i_{n}}}(\eta_{n}^{(1)}) \\
& = 
 \sum_{1 \le k < l \le n} 
 H_{\vpi_{i_{k}},\vpi_{i_{l}}}(\eta_{k} \otimes \eta_{l}^{(k+1)}) +
 \sum_{k=1}^{n} \Deg_{\vpi_{i_{k}}} (\eta_{k}^{(1)}).
\end{align*}
This completes the proof of equation \eqref{eq:step1},
thereby establishing Theorem~\ref{thm:step1}.
\end{proof}
%
%
%
%==============================%
%     START SUBSECTION 0405    %
%==============================%
%
\subsection{Proof of Theorem~\ref{thm:step2}.}
\label{subsec:step2}
Fix an arbitrary $i \in I_{0}$. 
Recall from Corollary~\ref{cor:tpd}\,(1) that 
there exists an isomorphism 
$\Psi_{\vpi_{i},\vpi_{i}}:\BB(2\vpi_{i})_{\cl} 
 \stackrel{\sim}{\rightarrow} 
 \BB(\vpi_{i})_{\cl} \otimes \BB(\vpi_{i})_{\cl}$
of $P_{\cl}$-crystals. 
The next lemma follows 
from the proof of \cite[Proposition~3.4.4]{NSp2}.
%
%%%%%%%%%%%%%%%%%
%%% lem:2picl %%%
%%%%%%%%%%%%%%%%%
%
\begin{lem} \label{lem:2picl}
Let $\eta_{1} \otimes \eta_{2} \in 
 \BB(\vpi_{i})_{\cl} \otimes \BB(\vpi_{i})_{\cl}$. 
Then, the inverse image
$\Psi_{\vpi_{i},\vpi_{i}}^{-1}(\eta_{1} \otimes \eta_{2}) 
 \in \BB(2\vpi_{i})_{\cl}$ of $\eta_{1} \otimes \eta_{2} \in 
 \BB(\vpi_{i})_{\cl} \otimes \BB(\vpi_{i})_{\cl}$ 
under the isomorphism $\Psi_{\vpi_{i},\vpi_{i}}$ is equal to
the concatenation $\eta_{1} \ast \eta_{2}$ of 
$\eta_{1}$ and $\eta_{2}$ defined by\,{\rm:}
%
%%%%%%%%%%%%%%
%%% eq:cat %%%
%%%%%%%%%%%%%%
%
\begin{equation} \label{eq:cat}
(\eta_{1} \ast \eta_{2})(t)=
\begin{cases}
 \eta_{1}(2t) 
   & \text{\rm if \ } 0 \le t \le 1/2, \\[2mm]
 \eta_{1}(1)+\eta_{2}(2t-1) 
   & \text{\rm if \ } 1/2 \le t \le 1.
\end{cases}
\end{equation}
\end{lem}

In addition, 
from \cite[Proposition~3.2.2]{NSp2}, we have
%
%%%%%%%%%%%%%%%
%%% lem:2pi %%%
%%%%%%%%%%%%%%%
%
\begin{lem} \label{lem:2pi}
The set $\BB(2\vpi_{i})$ 
coincides with the set of all concatenations $\pi_{1} \ast \pi_{2}$ 
of LS paths $\pi_{1},\,\pi_{2} \in \BB(\vpi_{i})$ 
such that $\kappa(\pi_{1}) \ge \iota(\pi_{2})$. 
Here, the concatenation $\pi_{1} \ast \pi_{2}$ is defined 
by the same formula as in \eqref{eq:cat}.
\end{lem}
%
%
%%%%%%%%%%%%%%%%
%%% lem:s2-1 %%%
%%%%%%%%%%%%%%%%
%
\begin{lem} \label{lem:s2-1}
Let $\eta \in \BB(2\vpi_{i})_{\cl}$, and set 
$\eta_{1} \otimes \eta_{2}:=
 \Psi_{\vpi_{i},\vpi_{i}}(\eta) \in 
 \BB(\vpi_{i})_{\cl} \otimes \BB(\vpi_{i})_{\cl}$. 
Then,
%
%%%%%%%%%%%%%%%
%%% eq:s2-1 %%%
%%%%%%%%%%%%%%%
%
\begin{equation} \label{eq:s2-1}
\Deg_{2\vpi_{i}}(\eta)=
 H_{\vpi_{i},\vpi_{i}}
 (\eta_{1} \otimes \eta_{2}) +
 2\Deg_{\vpi_{i}} (\eta_{1}).
\end{equation}
\end{lem}

\begin{proof}
Applying Theorem~\ref{thm:step1} to the case in which 
$\bi=(i,\,i)$ and hence $\lambda=2\vpi_{i}$, we have
\begin{align*}
\Deg_{2\vpi_{i}}(\eta) 
& =
 H_{\vpi_{i},\vpi_{i}}
 (\eta_{1} \otimes \eta_{2}^{(2)}) +
 \Deg_{\vpi_{i}} (\eta_{1}^{(1)}) + 
 \Deg_{\vpi_{i}} (\eta_{2}^{(1)}) \\
& =
 H_{\vpi_{i},\vpi_{i}}
 (\eta_{1} \otimes \eta_{2}) +
 \Deg_{\vpi_{i}} (\eta_{1}) + 
 \Deg_{\vpi_{i}} (\eta_{2}^{(1)}), 
\end{align*}
since $\eta_{1}^{(1)}=\eta_{1}$ and 
$\eta_{2}^{(2)}=\eta_{2}$ by definition.
From Lemma~\ref{lem:simple}\,(3) and 
the definition of $\eta_{2}^{(1)}$, we deduce that 
$\eta_{2}^{(1)}=\eta_{1}$, since the $P_{\cl}$-crystal 
$\BB(\vpi_{i})_{\cl} \otimes \BB(\vpi_{i})_{\cl}$ is simple. 
Thus we obtain \eqref{eq:s2-1}. This proves the lemma.
\end{proof}
Let us fix an (arbitrary) 
$\eta^{\flat} \in \BB(\vpi_{i})_{\cl}$ such that 
$f_{j}\eta^{\flat}=\bzero$ for all $j \in I_{0}$. 
Then it follows from \eqref{eq:cl-ro} that 
$f_{j}\pi_{\eta^{\flat}}=\bzero$ for all $j \in I_{0}$. 
Therefore, by Lemma~\ref{lem:final}, 
the final direction $\nu^{\flat}:=
\kappa(\pi_{\eta^{\flat}}) \in W\vpi_{i}$ of 
$\pi_{\eta^{\flat}}$ 
satisfies the condition that 
$\nu^{\flat}(h_{j}) \le 0$ for all $j \in I_{0}$.
Hence $\cl(\nu^{\flat})=w_{0}(\cl(\vpi_{i}))$ 
by Remark~\ref{rem:orbcl}. From this, we deduce 
(using Lemma~\ref{lem:orb}) that 
$\nu^{\flat} \in W\vpi_{i}$ can be written as 
$\nu^{\flat}=w_{0}\vpi_{i}+k^{\flat}d_{\vpi_{i}}\delta$ 
for some $k^{\flat} \in \BZ$.
%
%
%%%%%%%%%%%%%%%%
%%% lem:s2-2 %%%
%%%%%%%%%%%%%%%%
%
\begin{lem} \label{lem:s2-2}
Let $\eta \in \BB(\vpi_{i})_{\cl}$, and set 
$\eta':=
 \Psi_{\vpi_{i},\vpi_{i}}^{-1}(\eta^{\flat} \otimes \eta) \in 
 \BB(2\vpi_{i})_{\cl}$. 

\noindent {\rm (1)} 
The path 
$\pi_{\eta^{\flat}} \ast 
 (\pi_{\eta}+\pi_{k^{\flat}d_{\vpi_{i}}\delta})$ 
lies in $\cl^{-1}(\eta') \cap \BB(2\vpi_{i})$. 

\noindent {\rm (2)}
$\pi_{\eta'}=\pi_{\eta^{\flat}} \ast 
 (\pi_{\eta}+\pi_{k^{\flat}d_{\vpi_{i}}\delta})$.
\end{lem}

\begin{proof}
(1) First, note that $\eta'$ is equal to the concatenation 
$\eta^{\flat} \ast \eta$ by Lemma~\ref{lem:2picl}. 
In addition, by the definition of 
concatenations (see \eqref{eq:cat}), 
$\cl\bigl(
 \pi_{\eta^{\flat}} \ast 
 (\pi_{\eta}+\pi_{k^{\flat}d_{\vpi_{i}}\delta})
\bigr) = 
\cl(\pi_{\eta^{\flat}}) \ast 
\cl(\pi_{\eta}+\pi_{k^{\flat}d_{\vpi_{i}}\delta})$.
Since $\cl(\pi_{\eta^{\flat}})=\eta^{\flat}$ and 
$\cl(\pi_{\eta}+\pi_{k^{\flat}d_{\vpi_{i}}\delta})=
\cl(\pi_{\eta})=\eta$, we deduce that
\begin{equation*}
\cl\bigl(
 \pi_{\eta^{\flat}} \ast 
 (\pi_{\eta}+\pi_{k^{\flat}d_{\vpi_{i}}\delta})
\bigr) = 
\cl(\pi_{\eta^{\flat}}) \ast 
\cl(\pi_{\eta}+\pi_{k^{\flat}d_{\vpi_{i}}\delta}) = 
\eta^{\flat} \ast \eta = \eta', 
\end{equation*}
and hence that $\pi_{\eta^{\flat}} \ast 
 (\pi_{\eta}+\pi_{k^{\flat}d_{\vpi_{i}}\delta}) 
\in \cl^{-1}(\eta')$.

Next, we show that
$\pi_{\eta^{\flat}} \ast 
 (\pi_{\eta}+\pi_{k^{\flat}d_{\vpi_{i}}\delta}) \in \BB(2\vpi_{i})$. 
Note that $\pi_{\eta}+\pi_{k^{\flat}d_{\vpi_{i}}\delta} \in \BB(\vpi_{i})$ 
by Lemma~\ref{lem:delta}. Therefore, by means of Lemma~\ref{lem:2pi}, 
it suffices to show that 
$\iota(\pi_{\eta}+\pi_{k^{\flat}d_{\vpi_{i}}\delta}) \le 
 \nu^{\flat}=\kappa(\pi_{\eta^{\flat}})$.
Let us write 
$\iota(\pi_{\eta}) \in \vpi_{i}-\fin{Q}_{+}$ as 
$\iota(\pi_{\eta}) = w\vpi_{i}$ 
for some $w \in \fin{W}$ (see Lemma~\ref{lem:orb}). 
Since $w_{0} \in \fin{W}$ is greater than or equal to $w \in \fin{W}$ 
with respect to the usual Bruhat ordering on the (finite) 
Weyl group $\fin{W}$ of $\Fg_{I_{0}}$, 
it follows (see \cite[Remark~4.2]{L2}) that 
$w_{0}\vpi_{i} \ge w\vpi_{i}$. 
From this, using the definition of 
the ordering $\ge$ on $W\vpi_{i}$, 
we deduce that 
%
%%%%%%%%%%%%%%
%%% s2-2-1 %%%
%%%%%%%%%%%%%%
%
\begin{equation} \label{eq:s2-2-1}
\nu^{\flat}=
w_{0}\vpi_{i}+k^{\flat}d_{\vpi_{i}}\delta \ge 
w\vpi_{i} +k^{\flat}d_{\vpi_{i}}\delta=
\iota(\pi_{\eta}) +k^{\flat}d_{\vpi_{i}}\delta.
\end{equation}
Hence, by Remark~\ref{rem:sum}, 
\begin{equation*}
\iota(\pi_{\eta}+\pi_{k^{\flat}d_{\vpi_{i}}\delta}) =
\iota(\pi_{\eta})+\iota(\pi_{k^{\flat}d_{\vpi_{i}}\delta})=
\iota(\pi_{\eta}) +k^{\flat}d_{\vpi_{i}}\delta 
\le w_{0}\vpi_{i}+k^{\flat}d_{\vpi_{i}}\delta=\nu^{\flat}.
\end{equation*}
This proves part (1). 

\vsp

\noindent 
(2) First, by the definitions of 
$\Deg_{\vpi_{i}}(\eta^{\flat})$ and 
$\Deg_{\vpi_{i}}(\eta)$, we can write 
$\pi_{\eta^{\flat}}(1) \in P$ and $\pi_{\eta}(1) \in P$ as 
$\pi_{\eta^{\flat}}(1)=\vpi_{i}-a_{0}^{-1}\beta^{\flat}-
  a_{0}^{-1} \Deg_{\vpi_{i}}(\eta^{\flat})$ and 
$\pi_{\eta}(1)=\vpi_{i}-a_{0}^{-1}\beta-
  a_{0}^{-1} \Deg_{\vpi_{i}}(\eta)$ 
for some $\beta^{\flat},\,\beta \in \fin{Q}_{+}$, 
respectively.
Hence we have
%
%%%%%%%%%%%%%%
%%% eq:end %%%
%%%%%%%%%%%%%%
%
\begin{align} 
& 
(\pi_{\eta^{\flat}} \ast 
(\pi_{\eta}+\pi_{k^{\flat}d_{\vpi_{i}}\delta}))(1) = 
\pi_{\eta^{\flat}}(1) + 
(\pi_{\eta}+\pi_{k^{\flat}d_{\vpi_{i}}\delta})(1)= 
\pi_{\eta^{\flat}}(1) + 
\pi_{\eta}(1)+\pi_{k^{\flat}d_{\vpi_{i}}\delta}(1) \nonumber \\
& \hspace*{5mm}
= \bigl\{\vpi_{i}-a_{0}^{-1}\beta^{\flat}-
  a_{0}^{-1} \Deg_{\vpi_{i}}(\eta^{\flat})\bigr\} + 
  \bigl\{\vpi_{i}-a_{0}^{-1}\beta-
  a_{0}^{-1} \Deg_{\vpi_{i}}(\eta)\bigr\} + k^{\flat}d_{\vpi_{i}}\delta
\nonumber \\
& \hspace*{5mm}
= 
2\vpi_{i}-a_{0}^{-1}(\beta^{\flat}+\beta)-a_{0}^{-1}
(\Deg_{\vpi_{i}}(\eta^{\flat})+
\Deg_{\vpi_{i}}(\eta)-a_{0}k^{\flat}d_{\vpi_{i}})\delta. 
\label{eq:end} 
\end{align}
Now, in view of Proposition~\ref{prop:deg-min}\,(1), 
it suffices to show the following: 
%
%%%%%%%%%%%%%%%%%
%%% eq:s2-2-2 %%%
%%%%%%%%%%%%%%%%%
%
\begin{equation} \label{eq:s2-2-2}\iota
 \bigl(\pi_{\eta^{\flat}} \ast 
 (\pi_{\eta}+\pi_{k^{\flat}d_{\vpi_{i}}\delta})
 \bigr) \in 2\vpi_{i}-\fin{Q}_{+}, 
\end{equation}
%
%%%%%%%%%%%%%%%%%
%%% eq:s2-2-3 %%%
%%%%%%%%%%%%%%%%%
%
\begin{equation} \label{eq:s2-2-3}
\Deg_{2\vpi_{i}}(\eta') \le \Deg_{\vpi_{i}}(\eta^{\flat})+
 \Deg_{\vpi_{i}}(\eta)-a_{0}k^{\flat}d_{\vpi_{i}}. 
\end{equation}
But, we easily deduce from the definition of 
concatenations that 
$\iota
 \bigl(\pi_{\eta^{\flat}} \ast 
 (\pi_{\eta}+\pi_{k^{\flat}d_{\vpi_{i}}\delta})
 \bigr)
 =
 2\iota(\pi_{\eta^{\flat}})$,
which lies in $2\vpi_{i}-\fin{Q}_{+}$ 
since $\iota(\pi_{\eta^{\flat}}) \in \vpi_{i}-\fin{Q}_{+}$ 
by definition. 

Let us show \eqref{eq:s2-2-3}. 
Set $\pi_{1}(t):=\pi_{\eta'}(\frac{t}{2})$ and 
$\pi_{2}(t):=\pi_{\eta'}(\frac{t+1}{2})-
 \pi_{\eta'}(\frac{1}{2})$ for $t \in [0,1]$.
Then it is obvious that 
$\pi_{\eta'}=\pi_{1} \ast \pi_{2}$. 
Also, we see from the proof of 
\cite[Proposition~3.2.2]{NSp2} that 
$\pi_{1},\,\pi_{2} \in \BB(\vpi_{i})$. 
Furthermore, since 
$\eta^{\flat} \ast \eta=\eta'=
 \cl(\pi_{\eta'})=
 \cl(\pi_{1}) \ast \cl(\pi_{2})$, 
we see that 
$\cl(\pi_{1})=\eta^{\flat}$ and 
$\cl(\pi_{2})=\eta$.

Since $\iota(\pi_{\eta'}) \in 2\vpi_{i}-\fin{Q}_{+}$, and 
since $\iota(\pi_{1})=\frac{1}{2}\iota(\pi_{\eta'})$ 
by the definition of $\pi_{1} \in \BB(\vpi_{i})$,
it follows from Lemma~\ref{lem:orb} that 
$\iota(\pi_{1}) \in W\vpi_{i}$ lies in $\vpi_{i}-\fin{Q}_{+}$. 
If we write $\pi_{1}(1) \in P$ in the form
$\pi_{1}(1)=\vpi_{i}-a_{0}^{-1}\beta_{1} + 
 a_{0}^{-1}K_{1}\delta$ 
with $\beta_{1} \in \fin{Q}_{+}$ and 
$K_{1} \in \BZ_{\ge 0}$ (see Lemma~\ref{lem:end}), 
then by Proposition~\ref{prop:deg-min}\,(1) applied to 
$\cl(\pi_{1})=\eta^{\flat} \in \BB(\vpi_{i})_{\cl}$, 
%
%%%%%%%%%%%%%%%
%%% eq:s2-2 %%%
%%%%%%%%%%%%%%%
%
\begin{equation} \label{eq:s2-2}
-K_{1} \le \Deg_{\vpi_{i}}(\eta^{\flat}).
\end{equation}
Also, if we write $\kappa(\pi_{1}) \in W\vpi_{i}$ 
in the form 
$\vpi_{i}-\beta_{1}'+k_{1}'d_{\vpi_{i}}\delta$ with 
$\beta_{1}' \in \fin{Q}_{+}$ and $k_{1}' \in \BZ$, 
then by Proposition~\ref{prop:deg-min}\,(2) 
applied to $\cl(\pi_{1})=\eta^{\flat} \in \BB(\vpi_{i})_{\cl}$, 
%
%%%%%%%%%%%%%%%%
%%% eq:s2-21 %%%
%%%%%%%%%%%%%%%%
%
\begin{equation} \label{eq:s2-21}
k_{1}' \ge k^{\flat},
\end{equation}
where we recall from the discussion preceding this lemma that 
$\nu^{\flat}=\kappa(\pi_{\eta^{\flat}}) \in W\vpi_{i}$ equals 
$w_{0}\vpi_{i}+k^{\flat}d_{\vpi_{i}}\delta$.

Since $\pi_{\eta'}=
\pi_{1} \ast \pi_{2} \in \BB(2\vpi_{i})$, 
it follows from Lemma~\ref{lem:2pi} that 
$\iota(\pi_{2}) \le \kappa(\pi_{1})=
 \vpi_{i}-\beta_{1}'+k_{1}'d_{\vpi_{i}}\delta$.
Hence, by Lemma~\ref{lem:orb} and 
Remark~\ref{rem:bruhat}, 
$\iota(\pi_{2}) \in W\vpi_{i}$ is of the form 
$\iota(\pi_{2})=\vpi_{i}-\beta_{2}'+k_{2}'d_{\vpi_{i}}\delta$, 
with $\beta_{2}' \in \fin{Q}_{+}$ and $k_{2}' \in \BZ$ such that
%
%%%%%%%%%%%%%%%%
%%% eq:s2-22 %%%
%%%%%%%%%%%%%%%%
%
\begin{equation} \label{eq:s2-22}
k_{2}' \ge k_{1}'.
\end{equation}
If we set 
$\pi_{2}':=\pi_{2}-\pi_{k_{2}'d_{\vpi_{i}}\delta}$,
then it follows from Lemma~\ref{lem:delta} that 
$\pi_{2}' \in \BB(\vpi_{i})$.
In addition, we have
$\cl(\pi_{2}')=\cl(\pi_{2})=\eta$, and 
by Remark~\ref{rem:sum}, 
$\iota(\pi_{2}')=
 \iota(\pi_{2})-k_{2}'d_{\vpi_{i}}\delta=
 \vpi_{i}-\beta_{2}' 
 \in \vpi_{i}-\fin{Q}_{+}$. 
Hence, if we write $\pi_{2}'(1) \in P$ in the form
$\pi_{2}'(1)=
\vpi_{i}-a_{0}^{-1}\beta_{2} + a_{0}^{-1}K_{2}\delta$
with $\beta_{2} \in \fin{Q}_{+}$ and 
$K_{2} \in \BZ_{\ge 0}$ (see Lemma~\ref{lem:end}), 
then by Proposition~\ref{prop:deg-min}\,(1) 
applied to $\cl(\pi_{2}')=\eta \in \BB(\vpi_{i})_{\cl}$, 
%
%%%%%%%%%%%%%%%
%%% eq:s2-3 %%%
%%%%%%%%%%%%%%%
%
\begin{equation} \label{eq:s2-3}
-K_{2} \le \Deg_{\vpi_{i}}(\eta).
\end{equation}

From the above, we deduce that 
\begin{align*}
\pi_{\eta'}(1) & 
 = (\pi_{1} \ast \pi_{2})(1) 
 = \pi_{1}(1)+\pi_{2}(1) 
 = \pi_{1}(1)+\pi_{2}'(1)+k_{2}'d_{\vpi_{i}}\delta \\
& =
(\vpi_{i}-a_{0}^{-1}\beta_{1} + a_{0}^{-1}K_{1}\delta)+
(\vpi_{i}-a_{0}^{-1}\beta_{2} + a_{0}^{-1}K_{2}\delta)+
 k_{2}'d_{\vpi_{i}}\delta \\
& =
2\vpi_{i}-a_{0}^{-1}(\beta_{1}+\beta_{2}) + 
a_{0}^{-1}(K_{1}+K_{2}+a_{0}k_{2}'d_{\vpi_{i}})\delta, 
\end{align*}
and hence that 
$\Deg_{2\vpi_{i}}(\eta')=
 -K_{1}-K_{2}-a_{0}k_{2}'d_{\vpi_{i}}$.
By \eqref{eq:s2-2}, \eqref{eq:s2-3}, and 
the inequalities \eqref{eq:s2-21} and \eqref{eq:s2-22}, 
we obtain 
\begin{equation*}
\Deg_{2\vpi_{i}}(\eta')=
-K_{1}-K_{2}-a_{0}k_{2}'d_{\vpi_{i}} \le 
\Deg_{\vpi_{i}}(\eta^{\flat})+
\Deg_{\vpi_{i}}(\eta) - 
a_{0}k^{\flat}d_{\vpi_{i}},
\end{equation*}
which is the desired inequality \eqref{eq:s2-2-3}. 
This completes the proof of part (2). 
\end{proof}

Now we are ready to prove Theorem~\ref{thm:step2}. 

\begin{proof}[Proof of Theorem~\ref{thm:step2}]
Let $\eta \in \BB(\vpi_{i})_{\cl}$, and set 
$\eta':=
 \Psi_{\vpi_{i},\vpi_{i}}^{-1}(\eta^{\flat} \otimes \eta) \in 
 \BB(2\vpi_{i})_{\cl}$. 
It follows from Lemma~\ref{lem:s2-1} that 
%
%%%%%%%%%%%%%%%
%%% eq:s2-0 %%%
%%%%%%%%%%%%%%%
%
\begin{equation} \label{eq:s2-0}
\Deg_{2\vpi_{i}}(\eta')=
 H_{\vpi_{i},\vpi_{i}}
 (\eta^{\flat} \otimes \eta) +
 2\Deg_{\vpi_{i_{1}}} (\eta^{\flat}).
\end{equation}
In addition, we deduce from Lemma~\ref{lem:s2-2}\,(2) and 
the definition of $\Deg_{2\vpi_{i}}(\eta')$, 
by using \eqref{eq:end}, that 
%
%%%%%%%%%%%%%%%
%%% eq:s2-4 %%%
%%%%%%%%%%%%%%%
%
\begin{equation} \label{eq:s2-4}
\Deg_{2\vpi_{i}}(\eta')=
 \Deg_{\vpi_{i}}(\eta^{\flat})+
 \Deg_{\vpi_{i}}(\eta)-a_{0}k^{\flat}d_{\vpi_{i}}.
\end{equation}
Hence, by combining \eqref{eq:s2-0} and \eqref{eq:s2-4}, 
we obtain 
%
%%%%%%%%%%%%%%%
%%% eq:s2-5 %%%
%%%%%%%%%%%%%%%
%
\begin{equation} \label{eq:s2-5}
 H_{\vpi_{i},\vpi_{i}}
 (\eta^{\flat} \otimes \eta)= 
 \Deg_{\vpi_{i}}(\eta)-
 \Deg_{\vpi_{i}}(\eta^{\flat})-a_{0}k^{\flat}d_{\vpi_{i}}.
\end{equation}
In particular, by taking 
$\eta=\eta_{\cl(\vpi_{i})} \in \BB(\vpi_{i})_{\cl}$ 
in \eqref{eq:s2-5}, we obtain $H_{\vpi_{i},\vpi_{i}}
 (\eta^{\flat} \otimes \eta_{\cl(\vpi_{i})})= 
 \Deg_{\vpi_{i}}(\eta_{\cl(\vpi_{i})})
 -\Deg_{\vpi_{i}}(\eta^{\flat})-a_{0}k^{\flat}d_{\vpi_{i}}$. 
But, since $\Deg_{\vpi_{i}}(\eta_{\cl(\vpi_{i})})=0$ 
by Lemma~\ref{lem:deg-ro}\,(1), it follows that 
%
%%%%%%%%%%%%%%%
%%% eq:s2-6 %%%
%%%%%%%%%%%%%%%
%
\begin{equation} \label{eq:s2-6}
 H_{\vpi_{i},\vpi_{i}}
 (\eta^{\flat} \otimes \eta_{\cl(\vpi_{i})})= 
 -\Deg_{\vpi_{i}}(\eta^{\flat})-a_{0}k^{\flat}d_{\vpi_{i}}.
\end{equation}
Equation \eqref{eq:step2} follows immediately from 
\eqref{eq:s2-5} and \eqref{eq:s2-6}. 
This establishes Theorem~\ref{thm:step2}.
\end{proof}
%
%
%
%==============================%
%     START SUBSECTION 0405    %
%==============================%
%
\subsection{Relation to one-dimensional sums.}
\label{subsec:1dsum}

For each $i \in I_{0}$, $s \in \BZ_{\ge 0}$, 
we denote by $W_{s}^{(i)}$ the so-called 
Kirillov-Reshetikhin module (KR module for short)
over the quantum affine algebra $U_{q}'(\Fg)$, 
which is a finite-dimensional irreducible 
$U_{q}'(\Fg)$-module with specific Drinfeld polynomials
(for details, see the comment preceding 
 \cite[Conjecture~2.1]{HKOTT}). 
In \cite[Conjecture~2.1]{HKOTT} (see also \cite[Conjecture~2.1]{HKOTY}), 
they conjectured that the KR modules 
$W_{s}^{(i)}$, $i \in I_{0}$, $s \in \BZ_{\ge 0}$, have 
simple crystal bases $\CB^{i,s}$, called KR crystals. 
Further, assuming the existence of the KR crystals $\CB^{i,s}$ 
for $i \in I_{0}$ and $s \in \BZ_{\ge 1}$, they introduced a special kind of 
classically restricted one-dimensional sums (1dsums for short) associated to 
tensor products of $\CB^{i,s}$, $i \in I_{0}$, $s \in \BZ_{\ge 1}$. 
Here we should mention that 
it is confirmed through enough evidence 
(see, for example, \cite{Kas-rims} and \cite{FL}) that 
the level-zero fundamental representation 
$W(\vpi_{i})$ is indeed 
the KR module $W_{1}^{(i)}$, 
and hence 
the crystal basis $\CB(\vpi_{i})$ of $W(\vpi_{i})$ is 
the KR crystal $\CB^{i,1}$ for every $i \in I_{0}$ 
(see also \cite[Remark~2.3]{HKOTT}). 

Now, following the definition in 
\cite[\S3]{HKOTY} and \cite[\S3.3]{HKOTT} of 
classically restricted 1dsums, 
we define classically restricted 1dsums associated to tensor 
products of the simple crystals $\BB(\vpi_{i})_{\cl} \, 
(\simeq \CB(\vpi_{i}))$, $i \in I_{0}$, in place of $\CB^{i,1}$, 
$i \in I_{0}$, as follows. 
Let $\bi = (i_{1}, i_{2}, \ldots, i_{n})$ be 
a sequence of elements of $I_{0}$, and 
$\BB_{\bi}= 
 \BB(\vpi_{i_{1}})_{\cl} \otimes 
 \BB(\vpi_{i_{2}})_{\cl} \otimes \cdots \otimes
 \BB(\vpi_{i_{n}})_{\cl}$.
Then, for an element
$\mu \in \cl(P^{0}_{+})=\sum_{i \in I_{0}} \BZ_{\ge 0} \cl(\vpi_{i})$, 
the classically restricted 1dsum 
$X(\BB_{\bi},\mu; q)$ is defined by:
\begin{equation*}
X(\BB_{\bi},\mu; q)=\sum_{
 \begin{subarray}{c}
 b \in \BB_{\bi} \\[1mm]
 e_{j}b=\bzero \ (j \in I_{0}) \\[1mm]
 \wt b = \mu
 \end{subarray}} 
 q^{D_{\bi}(b)},
\end{equation*}
where $D_{\bi}:\BB_{\bi} \rightarrow \BZ$ 
is as defined in \eqref{eq:Di}. 
Because the root operators $e_{j}$, $j \in I_{0}$, 
on $\BB(\lambda)_{\cl}$ and those on $\BB_{\bi}$ are 
compatible with the isomorphism 
$\Psi_{\bi}:\BB(\lambda)_{\cl} \stackrel{\sim}{\rightarrow} \BB_{\bi}$ 
of $P_{\cl}$-crystals, we obtain the following corollary of 
Theorem~\ref{thm:main}.
%
%%%%%%%%%%%%%%%%%
%%% cor:1dsum %%%
%%%%%%%%%%%%%%%%%
%
\begin{cor} \label{cor:1dsum}
Let $\bi = (i_{1}, i_{2}, \ldots, i_{n})$ be 
a sequence of elements of $I_{0}$, and set 
$\lambda := \sum_{k=1}^{n} \vpi_{i_{k}} \in P^{0}_{+}$.
For every $\mu \in \cl(P^{0}_{+})=
\sum_{i \in I_{0}}\BZ_{\ge 0} \cl(\vpi_{i})$, 
the following equation holds\,{\rm:}
%
%%%%%%%%%%%%%%%%
%%% eq:1dsum %%%
%%%%%%%%%%%%%%%%
%
\begin{equation} \label{eq:1dsum}
\sum_{
 \begin{subarray}{c}
 \eta \in \BB(\lambda)_{\cl} \\[1mm]
 e_{j}\eta=\bzero \ (j \in I_{0}) \\[1mm]
 \eta(1) = \mu
 \end{subarray}} 
 q^{\Deg_{\lambda}(\eta)}
= q^{-D_{\bi}^{\ext}} X(\BB_{\bi}, \mu; q).
\end{equation}
\end{cor}
%
%%%%%%%%%%%%%%%%%
%%% rem:1dsum %%%
%%%%%%%%%%%%%%%%%
%
\begin{rem}[{see also \cite[Proposition~3.9]{HKOTT}}] \label{rem:1dsum}
Let $\lambda \in P^{0}_{+}$, and let $\mu \in \cl(P^{0}_{+})$. 
We see from Corollary~\ref{cor:1dsum} that 
$q^{-D_{\bi}^{\ext}} X(\BB_{\bi}, \mu; q)$ does not 
depend on the choice of a sequence 
$\bi=(i_{1}, i_{2}, \ldots, i_{n})$ of elements of $I_{0}$ 
such that $\lambda = \sum_{k=1}^{n} \vpi_{i_{k}}$. 
\end{rem}
%
%
%==============================%
%     START SUBSECTION 0407    %
%==============================%
%
\subsection{Relation to the Kostka-Foulkes polynomials.}
\label{subsec:kostka}
Let $i \in I_{0}$. 
Recall the element $\eta_{\ti{\vpi}_{i}} \in \BB(\vpi_{i})_{\cl}$ 
introduced in \S\ref{subsec:main}, 
where $\ti{\vpi}_{i}=w_{0}\cl(\vpi_{i}) \in P_{\cl}$. 
We note that 
if $\res_{I_{0}} \BB(\vpi_{i})_{\cl}$ is connected 
(see \S\ref{subsec:simple} for 
the definition of $\res_{I_{0}} \BB(\vpi_{i})_{\cl}$), 
then the element $\eta_{\ti{\vpi}_{i}}$ is 
the unique element of $\BB(\vpi_{i})_{\cl}$ such that 
$f_{j}\eta_{\ti{\vpi}_{i}}=\bzero$ for all $j \in I_{0}$, 
since the $P_{\cl}$-crystal $\BB(\vpi_{i})_{\cl}$ is regular. 
%
%%%%%%%%%%%%%%%%%
%%% lem:flat0 %%%
%%%%%%%%%%%%%%%%%
%
\begin{lem} \label{lem:flat0}
Let $i \in I_{0}$, and assume that 
$\res_{I_{0}} \BB(\vpi_{i})_{\cl}$ is connected.
Then, the equation 
$H_{\vpi_{i},\vpi_{i}}(\eta_{\ti{\vpi}_{i}} \otimes \eta)=0$ 
holds for all $\eta \in \BB(\vpi_{i})_{\cl}$.
\end{lem}
%
%%%%%%%%%%%%%%%%%
%%% rem:flat0 %%%
%%%%%%%%%%%%%%%%%
%
\begin{rem} \label{rem:flat0}
Since the $P_{\cl}$-crystal $\BB(\vpi_{i})_{\cl}$ is regular, 
we see that
$\res_{I_{0}} \BB(\vpi_{i})_{\cl}$ is connected if and only if 
$W(\vpi_{i})$ is irreducible when regarded as 
a $U'_{q}(\Fg)_{I_{0}}$-module by restriction. 
Therefore, 
we deduce from \cite[Lemma~4.3\,(ii)]{Kas-rims} that 
if $a_{i}^{\vee}=1$, then 
$\res_{I_{0}} \BB(\vpi_{i})_{\cl}$ is connected. 
\end{rem}

\begin{proof}[Proof of Lemma~\ref{lem:flat0}]
Since $\res_{I_{0}} \BB(\vpi_{i})_{\cl}$ is connected 
by the assumption of the lemma, 
there exists a monomial $X$ in the root operators 
$f_{j}$ for $j \in I_{0}$ such that $X\eta=\eta_{\ti{\vpi}_{i}}$. 
Since $f_{j}\eta_{\ti{\vpi}_{i}}=\bzero$ for all $j \in I_{0}$, 
it follows from the tensor product rule for crystals that
\begin{equation*}
X(\eta_{\ti{\vpi}_{i}} \otimes \eta)= 
\eta_{\ti{\vpi}_{i}} \otimes X\eta = 
\eta_{\ti{\vpi}_{i}} \otimes \eta_{\ti{\vpi}_{i}}.
\end{equation*}
Furthermore, we see from 
the proof of \cite[Lemma~1.6\,(1)]{AK} that 
\begin{equation*}
S_{w_{0}}X(\eta_{\ti{\vpi}_{i}} \otimes \eta)=
S_{w_{0}}
(\eta_{\ti{\vpi}_{i}} \otimes \eta_{\ti{\vpi}_{i}})=
(S_{w_{0}} \eta_{\ti{\vpi}_{i}}) \otimes 
(S_{w_{0}}\eta_{\ti{\vpi}_{i}}),
\end{equation*}
and from Remark~\ref{rem:tpd02} that 
$S_{w_{0}} \eta_{\ti{\vpi}_{i}}=\eta_{\cl(\vpi_{i})}$. 
Since $w_{0} \in \fin{W}$, it follows from the definition of $S_{w_{0}}$ 
(see \eqref{eq:sj}) that there exists 
a monomial $X'$ in the root operators 
$e_{j},\,f_{j}$ for $j \in I_{0}$ such that 
$S_{w_{0}}X(\eta_{\ti{\vpi}_{i}} \otimes \eta)=
X'X(\eta_{\ti{\vpi}_{i}} \otimes \eta)$. 
Therefore, we deduce that 
\begin{align*}
0 & = 
H_{\vpi_{i},\vpi_{i}}(\eta_{\cl(\vpi_{i})} \otimes \eta_{\cl(\vpi_{i})})
\quad
\text{by Theorem~\ref{thm:energy} (H2)} \\
& 
= H_{\vpi_{i},\vpi_{i}}
  ((S_{w_{0}}\eta_{\ti{\vpi}_{i}}) \otimes (S_{w_{0}}\eta_{\ti{\vpi}_{i}})) \\
& 
= H_{\vpi_{i},\vpi_{i}}(S_{w_{0}}X(\eta_{\ti{\vpi}_{i}} \otimes \eta)) \\
& 
= H_{\vpi_{i},\vpi_{i}}(X'X(\eta_{\ti{\vpi}_{i}} \otimes \eta)) \\
& = H_{\vpi_{i},\vpi_{i}}(\eta_{\ti{\vpi}_{i}} \otimes \eta)
\quad
\text{by Theorem~\ref{thm:energy} (H1)}.
\end{align*}
This proves the lemma. 
\end{proof}
%
%%%%%%%%%%%%%%%%%
%%% rem:flat1 %%%
%%%%%%%%%%%%%%%%%
%
\begin{rem} \label{rem:flat1}
Let $\bi=(i_{1},\,i_{2},\,\dots,\,i_{n})$ be a sequence of 
elements of $I_{0}$ such that 
$a_{i_{k}}^{\vee}=1$ for all $1 \le k \le n$. 
Then we know from Remark~\ref{rem:flat0} that 
$\res_{I_{0}} \BB(\vpi_{i_{k}})_{\cl}$ is connected 
for all $1 \le k \le n$. 
Therefore, we see from Lemma~\ref{lem:flat0} and 
the definitions of $D_{\bi}$ and $D_{\bi}^{\ext}$ that 
for every 
$\eta_{1} \otimes 
 \eta_{2} \otimes \cdots \otimes 
 \eta_{n} \in 
\BB_{\bi}=\BB(\vpi_{i_{1}})_{\cl} \otimes 
\BB(\vpi_{i_{2}})_{\cl} \otimes \cdots \otimes 
\BB(\vpi_{i_{n}})_{\cl}$, 
%
%%%%%%%%%%%%%%
%%% eq:a-1 %%%
%%%%%%%%%%%%%%
%
\begin{equation} \label{eq:a-1}
D_{\bi}(
  \eta_{1} \otimes 
  \eta_{2} \otimes \cdots \otimes 
  \eta_{n}) = 
\sum_{1 \le k < l \le n} 
 H_{\vpi_{i_{k}}, \vpi_{i_{l}}}
 (\eta_{k} \otimes \eta_{l}^{(k+1)}), 
\quad \text{and} \quad
D_{\bi}^{\ext}=0.
\end{equation}
\end{rem}

Now we restrict our attention to the case 
in which $\Fg$ is of type $A_{\ell-1}^{(1)}$.
Because $a_{i}^{\vee}=1$ for 
all $i \in I_{0}$ 
in the case of type $A_{\ell-1}^{(1)}$, 
we know from Remark~\ref{rem:flat0} that 
$\res_{I_{0}} \BB(\vpi_{i})_{\cl}$ is 
connected for all $i \in I_{0}$. 
In fact, we can check by direct calculation 
that every $\vpi_{i}$ is ``minuscule'', i.e., that 
$\vpi_{i}(\beta^{\vee}) \in \bigl\{0,\pm 1\bigr\}$ 
for all $\beta \in \fin{\Delta}:=\fin{W}\Pi_{I_{0}} \subset \rr$, 
where $\Pi_{I_{0}}=\bigl\{\alpha_{j}\bigr\}_{j \in I_{0}}$, and hence 
that $\vpi_{i}(\xi^{\vee}) \in \bigl\{0,\pm 1\bigr\}$ 
for all $\xi \in \prr$. 
Using this fact, we deduce from 
the definition of LS paths that 
%
%%%%%%%%%%%%%%%
%%% eq:mini %%%
%%%%%%%%%%%%%%%
%
\begin{equation} \label{eq:mini}
\BB(\vpi_{i})_{\cl}=
 \bigl\{\eta_{\mu} \mid \mu \in \cl(W\vpi_{i})=
 \fin{W}\cl(\vpi_{i})\bigr\}
\end{equation}
for every $i \in I_{0}$. 
Also, it follows from the definition of the root operators
$e_{j}$ (resp., $f_{j}$), $j \in I$, that 
for $\mu \in \fin{W}\cl(\vpi_{i})$ and $j \in I$: 
%
%%%%%%%%%%%%%%%%%%
%%% eq:mini-ro %%%
%%%%%%%%%%%%%%%%%%
%
\begin{equation} \label{eq:mini-ro}
\begin{array}{ll}
\text{(1)} & \text{
$e_{j}\eta_{\mu} \ne \bzero$ 
(resp., $f_{j}\eta_{\mu} \ne \bzero$) 
if and only if 
$\mu(h_{j})=-1$ (resp., $=1$);} \\[3mm]
\text{(2)} & \text{
if $e_{j}\eta_{\mu} \ne \bzero$ 
(resp., $f_{j}\eta_{\mu} \ne \bzero$), then 
$e_{j}\eta_{\mu} \ (\text{resp., } f_{j}\eta_{\mu}) = 
 \eta_{r_{j}\mu}$,} \\[3mm]
 & \text{with $r_{j}(\mu)=\mu+\cl(\alpha_{j})$ 
   (resp., $=\mu-\cl(\alpha_{j})$).}
\end{array}
\end{equation}

Here we recall from \cite[\S\S2.1 and 2.4]{NY} 
the $P_{\cl}$-crystals $B_{\cl(\vpi_{i})}$ (adapted to our notation), 
$i \in I_{0}=\bigl\{1,\,2,\,\dots,\,\ell-1\bigr\}$. 
Let $i \in I_{0}=\bigl\{1,\,2,\,\dots,\,\ell-1\bigr\}$. 
Then, the crystal $B_{\cl(\vpi_{i})}$ consists of 
all Young diagrams of shape $(1^{\ell})$ having $i$ dots (in all)
with at most one dot in each box. 
The weight $\wt(b) \in P_{\cl}$ of 
the element $b \in B_{\cl(\vpi_{i})}$ having dots exactly 
in the $j_{1}$-th box, $j_{2}$-th box,\,\dots,\,$j_{i}$-th box 
with $1 \le j_{1} < j_{2} < \dots < j_{i} \le \ell$, 
is given by: $\gamma_{j_{1}} + \gamma_{j_{2}} + \cdots + 
\gamma_{j_{i}} \in P_{\cl}$, where we set 
$\gamma_{j}:=\cl(\vpi_{j})-\cl(\vpi_{j-1})$ for $1 \le j \le \ell$, 
with $\vpi_{0}=\vpi_{\ell}=0 \in P$. 

Now, for each $i \in I_{0}=\bigl\{1,\,2,\,\dots,\,\ell-1\bigr\}$, 
we define a map 
$\Phi_{i}:B_{\cl(\vpi_{i})} \rightarrow \BB(\vpi_{i})_{\cl}$ 
by: $\Phi_{i}(b)=\eta_{\wt(b)}$ for $b \in B_{\cl(\vpi_{i})}$, 
which is easily seen to be a (well-defined) bijection that
preserves weights, since the subgroup $\fin{W} \subset W$ is 
isomorphic to the symmetric group $S_{\ell}$ permuting 
the $\gamma_{j}$, $1 \le j \le \ell$, and since 
$\cl(\vpi_{i})=\gamma_{1}+\gamma_{2}+\cdots+\gamma_{i}$.
Actually, the map 
$\Phi_{i}:B_{\cl(\vpi_{i})} \rightarrow \BB(\vpi_{i})_{\cl}$ 
is an isomorphism of $P_{\cl}$-crystals. 
Indeed, from the definition in \cite[\S2.4]{NY} of 
the Kashiwara operators $e_{j}$ (resp., $f_{j}$), $j \in I$, on 
$B_{\cl(\vpi_{i})}$, we see that for $b \in B_{\cl(\vpi_{i})}$ 
and $j \in I$, $e_{j}b \ne \bzero$ (resp., $f_{j}b \ne \bzero$) 
if and only if $(\wt(b))(h_{j})=-1$ (resp., $=1$), 
since 
$\gamma_{j}=\cl(\vpi_{j})-\cl(\vpi_{j-1})=
 \cl(\Lambda_{j})-\cl(\Lambda_{j-1})$ for $2 \le j \le \ell-1$, and 
$\gamma_{1}=\cl(\vpi_{1})=
 \cl(\Lambda_{1})-\cl(\Lambda_{0})$, 
$\gamma_{\ell}=-\cl(\vpi_{\ell-1})=
 -\cl(\Lambda_{\ell-1})+\cl(\Lambda_{0})$.
Also, if $e_{j}b \ne \bzero$ (resp., $f_{j}b \ne \bzero$), 
then we have $\wt(e_{j}b)=\wt(b)+\cl(\alpha_{j})$ 
(resp., $\wt(f_{j}b)=\wt(b)-\cl(\alpha_{j})$). 
Therefore, from \eqref{eq:mini-ro}, we conclude that 
$\Phi_{i}:B_{\cl(\vpi_{i})} \rightarrow \BB(\vpi_{i})_{\cl}$ 
is an isomorphism of $P_{\cl}$-crystals. 

Let $\bi=(i_{1},\,i_{2},\,\dots,\,i_{n})$ 
be an arbitrary sequence of elements of 
$I_{0}=\bigl\{1,\,2,\,\dots,\,\ell-1\bigr\}$ 
such that $i_{1} \ge i_{2} \ge \cdots \ge i_{n}$, 
and set $\lambda:=\sum_{k=1}^{n} \vpi_{i_{k}} \in P^{0}_{+}$; 
we denote this sequence 
$\bi=(i_{1},\,i_{2},\,\dots,\,i_{n})$ by 
$\lambda^{\dagger}$ when we regard it as a partition 
(or a Young diagram).
In the following, we identify an element 
$\mu=\sum_{i=1}^{\ell-1} \mu^{(i)} \cl(\vpi_{i}) 
\in \cl(P^{0}_{+})$ with the partition
\begin{equation*}
\left(
 \sum_{i=1}^{\ell-1}\mu^{(i)}, \ 
 \sum_{i=2}^{\ell-1}\mu^{(i)}, \ \dots, \ 
 \mu^{(\ell-1)}
\right);
\end{equation*}
note that the elements $\cl(\vpi_{i})$, $1 \le i \le \ell-1$, 
of $\Fh^{\ast}/\BQ\delta$ are linearly independent over $\BQ$. 
Because $\BB(\vpi_{i})_{\cl}$ is isomorphic 
as a $P_{\cl}$-crystal to $B_{\cl(\vpi_{i})}$ through 
the map $\Phi_{i}:
 B_{\cl(\vpi_{i})} \stackrel{\sim}{\rightarrow} 
 \BB(\vpi_{i})_{\cl}$, and, in addition, 
$\res_{I_{0}} \BB(\vpi_{i})_{\cl}$ is connected 
for every $i \in I_{0}=\bigl\{1,\,2,\,\dots,\,\ell-1\bigr\}$, 
we deduce from \cite[Corollary~4.3]{NY} 
along with \eqref{eq:a-1} that 
for every $\mu \in \cl(P^{0}_{+})$, 
the Kostka-Foulkes polynomial $K_{\mu^{t},\,\lambda^{\dagger}}(q)$
(defined in \cite[Chap.\,III, \S6]{M}) associated to 
the conjugate (or transpose) $\mu^{t}$
of the partition $\mu$ and 
the partition $\lambda^{\dagger}$ is 
equal to the following: 
%
%%%%%%%%%%%%%
%%% eq:ny %%%
%%%%%%%%%%%%%
%
\begin{equation} \label{eq:ny}
X(\BB_{\bi},\mu; q^{-1})=\sum_{
 \begin{subarray}{c}
 b \in \BB_{\bi} \\[1mm]
 e_{j}b=\bzero \ (j \in I_{0}) \\[1mm]
 \wt b = \mu
 \end{subarray}} 
 q^{-D_{\bi}(b)}.
\end{equation}
By combining this fact 
and Corollary~\ref{cor:1dsum} (along with \eqref{eq:a-1}), 
we obtain the following expression for 
the Kostka-Foulkes polynomials 
in terms of LS paths. 
%
%%%%%%%%%%%%%%%%%%
%%% cor:kostka %%%
%%%%%%%%%%%%%%%%%%
%
\begin{cor} \label{cor:kostka}
Assume that $\Fg$ is of type $A_{\ell-1}^{(1)}$, and 
keep the notation above. 
Let $\mu \in \cl(P^{0}_{+})=\sum_{i \in I_{0}} \BZ_{\ge 0}\cl(\vpi_{i})$. 
Then, the following equation holds\,{\rm:}
\begin{equation*}
K_{\mu^{t},\,\lambda^{\dagger}}(q)= 
\sum_{
 \begin{subarray}{c}
 \eta \in \BB(\lambda)_{\cl} \\[1mm]
 e_{j}\eta=\bzero \, (j \in I_{0}) \\[1mm]
 \eta(1)=\mu
 \end{subarray}} 
 q^{-\Deg_{\lambda}(\eta)}.
\end{equation*}
\end{cor}
%
%======================%
%     BIBLIOGRAPHY     %
%======================%

{\small
\setlength{\baselineskip}{15pt}
\renewcommand{\refname}{References}

}

\end{document}